\documentclass[12pt]{amsart}
\usepackage{amsbsy}
\usepackage{graphicx,epsfig,subfigure}
\begin{document}

%
%
\newtheorem{theorem}{Theorem}
\newtheorem{proposition}[theorem]{Proposition}
\newtheorem{lemma}[theorem]{Lemma}
\newtheorem{corollary}[theorem]{Corollary}
\newtheorem{definition}[theorem]{Definition}
\newtheorem{remark}[theorem]{Remark}
\numberwithin{equation}{section}
\numberwithin{theorem}{section}
\newcommand{\be}{\begin{equation}}
\newcommand{\ee}{\end{equation}}
\newcommand{\re}{{\mathbb R}}
\newcommand{\n}{\nabla}
\newcommand{\ren}{{\mathbb R}^N}
\newcommand{\iy}{\infty}
\newcommand{\pa}{\partial}
\newcommand{\ms}{\medskip\vskip-.1cm}
\newcommand{\mpb}{\medskip}
\newcommand{\ssk}{\smallskip}
\newcommand{\BB}{{\bf B}}
\newcommand{\Am}{{\bf A}_{2m}}
\newcommand{\bL}{\BB^*}
\newcommand{\bLs}{\BB}
\renewcommand{\a}{\alpha}
\renewcommand{\b}{\beta}
\newcommand{\g}{\gamma}
\newcommand{\ka}{\kappa}
\newcommand{\G}{\Gamma}
\renewcommand{\d}{\delta}
\newcommand{\D}{\Delta}
\newcommand{\e}{\varepsilon}
\renewcommand{\l}{\lambda}
\renewcommand{\o}{\omega}
\renewcommand{\O}{\Omega}
\newcommand{\s}{\sigma}
\renewcommand{\t}{\tau}
\renewcommand{\th}{\theta}
\newcommand{\z}{\zeta}
\newcommand{\wx}{\widetilde x}
\newcommand{\wt}{\widetilde t}
\newcommand{\noi}{\noindent}
\newcommand{\lb}{\left (}
\newcommand{\rb}{\right )}
\newcommand{\lsb}{\left [}
\newcommand{\rsb}{\right ]}
\newcommand{\lab}{\left \langle}
\newcommand{\rab}{\right \rangle }
\newcommand{\gap}{\vskip .5cm}
\newcommand{\bz}{\bar{z}}
\newcommand{\bg}{\bar{g}}
\newcommand{\Ba}{\bar{a}}
\newcommand{\bt}{\bar{\th}}
\def\com#1{\fbox{\parbox{6in}{\texttt{#1}}}}

\title
{\bf  On continuous
 branches\\ of very singular
 similarity solutions \\ of the
stable
thin film equation} 

\author {J.D.~Evans and V.A.~Galaktionov}

\address{Department of Mathematical Sciences, University of Bath,
 Bath BA2 7AY, UK}
\email{masjde@maths.bath.ac.uk}

\address{Department of Mathematical Sciences, University of Bath,
 Bath BA2 7AY, UK}
\email{vag@maths.bath.ac.uk}


\date{\today}



\begin{abstract}

The fourth-order thin film equation (the TFE--4)
 \[
 u_t = -\nabla \cdot (|u|^n \nabla \Delta u) + \D( |u|^{p-1}u),
\,\, \mbox{where} \,\, n>0, \,\,
 p>1,
 \]
 with the {\em stable} second-order diffusion term is considered.
For the first critical exponent
 $$
 p=p_0= \mbox{$n+1 + \frac 2N$} \quad \mbox{for} \,\,\, n \in (0,\mbox{$\frac
32$}),
 $$
 where $N \ge 1$ is the space dimension,
the standard free-boundary problem (FBP)  with zero height, zero
contact angle, and zero-flux conditions is shown to admit
continuous sets (branches) of source-type {\em very singular}
self-similar similarity solutions (VSSs),
 $$
 u(x,t) = t^{-\frac N{4+n N}} f(y), \quad y =
 x/t^{\frac 1{4+nN}}. 
 $$
  For the Cauchy problem (CP),
  continuous branches
  of oscillatory self-similar patterns of changing sign, which become ``limits"
  of countable sets of  FBP ones,
  are identified.

For  $p \not = p_0$, the set of VSSs
is shown to be finite and to consist of a countable family of
$p$-branches of similarity profiles that originate at a
sequence of critical  exponents $\{p_l, \, l \ge 0\}$. At $p=p_l$,
these $p$-branches appear via a nonlinear bifurcation mechanism
from a countable set of similarity solutions of the {\em second
kind}
 of the pure TFE
 $$
 u_t = -\nabla \cdot (|u|^n \nabla \Delta u)
  \quad \mbox{in}
  \quad \ren \times \re_+.
 $$
 Such solutions are
 detected by a combination of linear and nonlinear ``Hermitian spectral theory" (in both the CP and FBP settings), which
allows us to apply an analytical $n$-branching approach. This means
constructing a continuous path as $n \to 0^+$  to eigenfunctions
of a linear rescaled  operator for $n=0$, i.e., for the
bi-harmonic equation $u_t=-\D^2 u$. Numerics are used to confirm
several analytical conclusions, which do not admit a fully rigorous
study.


\end{abstract}

\maketitle

\section{Introduction: the stable TFE and main results}

\subsection{The model and preliminary discussion}

This paper is devoted to the study of large time  behaviour of
solutions of higher-order quasilinear degenerate parabolic equations
of not-divergent forms.
 More precisely, we
construct
   global in time self-similar {\em very singular solutions} (VSSs) of
the fourth-order quasilinear parabolic {\em thin film equation}
(TFE--4) with the {\em stable} homogeneous second-order diffusion
term,
\be
\label{GPP}
    u_t = -\nabla \cdot (|u|^n \nabla \Delta u)+ \D (|u|^{p-1} u ) ,
     \quad \mbox{where} \,\,\,
      n>0 \,\,\mbox{and} \,\, p>1.
\ee
 The present results
  complete the analysis of the TFEs performed
in \cite{Bl4, Gl4}, where countable sets and continuous branches
of {\em blow-up} and {\em global} similarity solutions were
obtained for the {\em limit unstable TFE} with the {\em backward
diffusion} parabolic term,
\be
\label{GPPun}
    u_t = -\nabla \cdot (|u|^n \nabla \Delta u) - \D (|u|^{p-1} u )
      \quad (n>0, \,\,\, p>1).
\ee
 The main mathematical approaches for (\ref{GPPun}) also exhibit
 certain
similarities with those applied in \cite{GBl6} to the sixth-order
limit unstable TFE--6
 \be
\label{GPP6}
  u_t = \n \cdot (|u|^{n}\n \D^2 u) - \D ( |u|^{p-1} u) \quad 
%
  (n>0, \,\,\,
 p>1).
\ee
Surveys and  extended lists of related references on the physics and
mathematics of such thin film PDEs can be found in these papers \cite{GBl6} and \cite{Bl4}.
As usual, for our future analysis, both pioneering papers of higher-order nonlinear
diffusion theory in the 1990s by Bernis--Friedman \cite{BF1} (mainly, FBP theory for TFEs),
Bernis \cite{Bern88} and Bernis--McLeod  \cite{BMcL91}, where oscillatory similarity solutions
of the Cauchy problem for the fourth-order porous medium-like equations (the PME--4) were constructed,
are key; see the monograph \cite[Ch.~4]{Wu01} for further details in these areas.

\ssk

We consider for (\ref{GPP}) the standard  FBP with {\em
zero-height}, {\em zero contact angle}, and {\em conservation of
mass} ({\em zero flux}) conditions,
\be
 \label{GPP1}
 \mbox{$
 u=\n u=  -  {\bf n} \cdot (|u|^n  \nabla \D u - \nabla (|u|^{p-1} u) ) =0
 $}
\ee
at the singularity surface (interface) $\Gamma_0[u] $, which
is the lateral boundary of $ {\rm supp} \,u \subset \ren \times
\re_+$ with  the unit outward normal ${\bf n}$.
 For the range of
$n\in(0,\frac 32)$ plus the semilinear case $n=0$ also discussed
here, these three
conditions are expected to give a correctly specified problem for
the fourth-order parabolic equation, which is completed with
bounded, smooth, integrable, compactly supported initial data
 \be
 \label{u00}
 u(x,0)=u_0(x) \quad \mbox{in} \,\,\, \Gamma_0[u] \cap \{t=0\}.
 \ee

We will also treat the Cauchy problem  (CP) for (\ref{GPP}) with
compactly supported initial data (\ref{u00})  in $\ren$. The CP
admitting oscillatory solutions of ``maximal" regularity will
need a special setting.

We begin our study in the critical ``conservative" case
 \be
 \label{p011}
p=p_0 = n+1 +\mbox{$ \frac 2N$},
  \ee
which is easier technically and reveals specific properties of
similarity patterns. Eventually, we extend our approach to $p \not
= p_0$ (more precisely, for $p<p_0$). Notice that, for $n=0$,
equation (\ref{GPP}) is the  {\em limit stable Cahn--Hilliard
equation}
 \be
 \label{CH11}
u_t = -\Delta^2 u + \Delta (|u|^{p-1}u),
 \ee
 which occurs in
 various applications; see references in  
 \cite{EGW}.


\subsection{Main results and layout of the paper}

We construct  very singular  self-similar (or source-type for
$p=p_0$) solutions of (\ref{GPP}) in certain ranges of the parameters
$n$, $p$ and $N$. For small enough $n>0$, we will often refer to
analogies with the semilinear Cahn--Hilliard equation
(\ref{CH11}). Typically,  we assume that
 \be
\label{nn22} n \in (0, \mbox{$\frac 32$}) \quad \mbox{and} \quad
p>n+1,
 \ee
 but sometimes we also treat $n> \frac 32$, for which the CP
 continues to admit sign changing solutions that are infinitely
 oscillatory at the interfaces.


In Section \ref{Sect2}
 we formulate the similarity setting of the problem.
 Our conclusions and further
layout of the paper are as follows: we show that the stable TFE
 (\ref{GPP})
admits:

{\bf (i)} In the critical case $p=p_0$, continuous families of
global similarity solutions of  the CP (Section \ref{SectCP1}) and
of the FBP (Section \ref{SectLocR});

{\bf (ii)} As a co-product, we first study
 countable set of
similarity solutions of the pure TFE (these define special  bifurcation
values $\{p_l\}$  for the full model (\ref{GPP}), Section
\ref{CPTFEn})
\be
  \label{tf1}
  u_t = - \n \cdot (|u|^n \n \D u) \quad \mbox{in}
  \quad \ren \times \re_+; \quad \mbox{and}
   \ee

{\bf (iii)} For $p < p_0$ (for $p>p_0$ no such VSSs exist), the
number of similarity solutions (for a given $p$ value) becomes finite, and in
the CP, there exists a countable family of $p$-branches of similarity profiles that
 originate at certain nonlinear bifurcation points $\{p_l>1, \, l
 \ge 0 \}$ (Section \ref{CPpn}). Some of the results are extended to the
 FBP (Section \ref{Sectp0}).

\ssk

The principal issue that occurs is the actual relation between
similarity solutions of the standard FBP and the Cauchy problem
(the latter are infinitely oscillatory near the interfaces for not
that large $n$). More clearly and convincingly than in our
previous research, we show that, in several cases, {\em for each
asymptotic pattern of the CP, there exists a countable set of FBP
patterns, which eventually converge to the CP one}. It turns out
that this is a rather general principal even for the linear
problem for $n=0$,
i.e., for the {\em bi-harmonic equation},
  \be
  \label{tf2}
  u_t= - \D^2 u.
 \ee
 To show this, we will need to develop a type of {\em Hermitian
 spectral theory} for linear rescaled operators for both the CP
 (Section \ref{S4.1}) and for the FBP setting (Section
 \ref{S9.2}). The latter one becomes more difficult and unusual, with a multi-dimensional space of eigenvalues, so we discuss
 just initial aspects of such a theory therein.

\ssk

Similar principles of  self-similar asymptotics apply \cite{GBl6}
to the sixth-order stable TFE
 \be
\label{GPP6st}
  u_t = \n \cdot (|u|^{n}\n \D^2 u) + \D ( |u|^{p-1} u), \quad 
%
\ee
 for $n \in [0, \frac 54)$.
 In this case, the first critical exponent is
  $
  p_0 = n+1 + \mbox{$\frac 4 N$}.
   $
  The semilinear case $n=0$ leads to the
sixth-order unstable limit Cahn--Hilliard equation
 \be
\label{CH6} u_t = \D^3 u + \D (|u|^{p-1}u),
 \ee
whose similarity solutions can be studied as in \cite{EGW}.

\section{Global  similarity solutions: general statement and preliminaries}
 \label{Sect2}


 The similarity solutions of (\ref{GPP}) have the  form
\be
\label{ResVars}
 u_S(x,t) = t^{-\a} f(y),
\quad y = x/t^\b, \quad \mbox{with} \quad
  \a =  \mbox{$\frac 1{2p-(n+2)}>0, \,\, \beta=
  \frac{p-(n+1)}{2[2p-(n+2)]}>0$}.
  \ee
The function $f$ solves a quasilinear elliptic equation, namely,
 \be
 \label{RescProf}
{\bf A}_+(f) \equiv - \nabla \cdot \bigl[|f|^n \nabla \D f -
\nabla
( |f|^{p-1} f ) \bigr] + \b y \cdot \nabla f + \a f= 0. 
\,\,
 \ee
For $n>0$, a natural functional setting for both the FBP and the
CP includes the condition
 \be
 \label{csff}
f(y) \quad \mbox{is non-trivial in a bounded domain in  $\ren$.}
 \ee
In the CP, $f(y)$ can be extended by $f(y) \equiv 0$ outside the
support. For FBPs, posed by definition in a bounded domain, such
an extension is not applicable.
 Note that, for the CP,  the elliptic equation (\ref{RescProf})
 admits non-compactly supported solutions with the asymptotics as
 $y \to \iy$ governed by the leading linear first-order operator:
  \be
  \label{lf1}
  \b y \cdot \n f+ \a f+...=0 \quad \Longrightarrow \quad
   f(y) = C |y|^{- \frac \a \b}(1+o(1)),
    \ee
  where $C=C\big( \frac y{|y|}\big)$ is an arbitrary  smooth function on the unit sphere $S^{N-1}$.
  Actually, in order to satisfy the desired condition (\ref{csff}), one needs to demand that
  \be
  \label{lf2}
  C=0 \quad \mbox{in \,\, (\ref{lf1}).}
   \ee


For the CP with $n=0$, (\ref{csff}) is replaced by (see
\cite{EGW})
 \be
 \label{expff}
f(y) \quad \mbox{has exponential decay at infinity,}
 \ee
 meaning that $f$ belongs to a special weighted $L^2$ space.
 The condition (\ref{lf2}) is also necessarily implied.

Under (\ref{csff}) (or equivalently (\ref{expff})), integrating
 (\ref{ResVars}) over $\ren$ yields the following mass
 time-dependence for $p \not = p_0$:
  \be
  \label{mma1}
   \mbox{$
    \int u_{S}(x,t) \, {\mathrm d}x=
    t^{\frac{N(p-p_0)}{2[2p-(n+1)]}} \, \int f(y) \, {\mathrm d}y
 \quad \Longrightarrow \quad \int f=0 \quad(p \not = p_0).
  $}
  \ee
For $p=p_0$, any mass of $f(y)$ is formally allowed; cf.
\cite{Gl4}.

For  general solutions of the TFE--4 (\ref{GPP}), the self-similar scaling
 \be
 \label{sc55}
 u(x,t) = (1+t)^{-\a} \th(y,\t), \quad y=x/(1+t)^\b, \quad \t =
 \ln(1+t),
 \ee
 yields the evolution equation with the same operator as in
 (\ref{RescProf}),
  \be
  \label{Resc0p}
  \th_\t = {\bf A}_+(\th) \quad \mbox{for \,\, $\t > 0$}.
  \ee
Then a typical asymptotic stabilization problem occurs as $\t \to +\iy$,
which in particular,  requires to know all possible equilibria  of
 ${\bf A}_+$, this being the main current problem of concern.



The critical exponent (\ref{p011}) follows from conservation of
mass; {\em q.v.} the same derivation in \cite[\S~3]{Bl4}.
In addition, a countable sequence of other critical exponents
$\{p_l, \, l=0,1,2,...\}$ is expected to exist. This is confirmed
in  the semilinear case $n=0$ (for the limit Cahn--Hilliard
equation (\ref{CH11})),
 where 
\cite[\S~5]{EGW}
 \be
 \label{Crp}
 p_l = 1+ \mbox{$\frac 2{N+l}$} \quad \mbox{for any} \quad  l =0,1,2,... \, .
 \ee
See further comments in \cite[\S~2]{Bl4}.


 \section{The CP: local oscillatory ODE bundles and profiles for $p=p_0$}
\label{SectCP1}

 We next study the similarity  ODE in the radial setting.
 Let $y$ denote the
single spatial variable $|y|\ge 0$.
  The
operator of the equation (\ref{ResVars})  is  then ordinary
differential,
 \be
 \label{ODEop}
 \mbox{$
 {\bf A}_+ (f) \equiv  - \frac
 1{y^{N-1}}\, \bigl[ \,y^{N-1}|f|^n \big(\frac1{y^{N-1}}(y^{N-1}f')'\big)'
 - y^{N-1}( |f|^{p-1} f )' \, \bigr]' + \b y f' + \a f=0.
  $}
 \ee

We first describe the corresponding oscillatory bundle of
asymptotic profiles close to interface points, which are attributed to the CP. It turns out
that we have to begin with the CP. Namely, we will show that similar  FBP profiles are
naturally associated with the CP ones and, in particular, for both
$n=0$ and $n>0$ demand certain extra ``Hermitian spectral theory" as an extended and more
difficult version of that in the Cauchy setting in $\ren$.

 \subsection{The Cauchy problem: local  oscillatory behaviour close to
interfaces}

These questions have been considered before, so,
 following \cite[\S~7.1]{Gl4}, we briefly
 indicate the oscillatory asymptotic
 bundle  of similarity profiles $f(y)$
exhibiting {\em maximal regularity} at the interface $y=y_0$,
so that, being extended by $f=0$ for $y > y_0$, these will give
solutions of the CP. It is easy to see that,
 for $n \in (0,\frac 32)$, the ODE (\ref{OD-}) does not admit nonnegative
solutions
 of the maximal regularity, which can be considered as a
 counterpart of smooth similarity
 solutions of the CP for $n=0$  \cite[\S~5]{EGW},
 i.e., admitting  a regular limit as $n \to 0^+$. 
Hence, the ODE (\ref{OD-}) implies that sufficiently regular
solutions $f(y)$ must be oscillatory near interfaces; cf.
 proofs  in \cite{BMcL91}. It is important to prescribe the precise structure of
such oscillatory singularities of the ODE and determine the
dimension of the asymptotic bundle.

 For
$n >0$, we take the thin film ODE (\ref{ODEop}), keeping the main
terms for $y \approx y_0^-$ and integrating once, to obtain
 \be
 \label{OD-}
  \mbox{$
 |f|^n \big( f'' + \frac{(N-1)}{y} f'  \big)' - ( |f|^{p-1}f )' = \b y_0 f +
  (\mbox{higher-order terms}) \,.
  $}
  \ee
  For $N=1$,
choosing next just two leading terms close to the interface yields
(in fact, one can see that this approximation holds for any
dimension $N \ge 1$)
 \be
 \label{7.1N}
  |f|^n f'''= \b y f+ ... = \l_0 f+ ... \,,
  \ee
   where $ \l_0 = \b
 y_0$. Thus, we need to consider the unperturbed ODE
 \be
 \label{fff.1}
 |f|^n f'''=  \l_0 f \quad (\l_0 = \b y_0>0),
 \ee
 whose orbits will be exponentially small perturbations
near interfaces of that for (\ref{ODEop}).

 The ODE (\ref{fff.1}) has
  the
following representation of the solutions \cite[\S~7]{Gl4}: 
as $y \to y_0^-$,
 \be
 \label{LC11}
 \mbox{$
 f(y) = (y_0-y)^{\mu} \varphi(\eta) \,, \quad \eta =
 \ln(y_0-y), \quad \mbox{with} \,\,\, \mu=\frac 3n,
  $}
 \ee
 where the {\em oscillatory component}  $\varphi$  
 satisfies the autonomous ODE 
 \be
 \label{eqLC}
 \varphi''' + 3(\mu-1) \varphi'' + (3 \mu^2 - 6 \mu +2) \varphi'
+ \mu(\mu-1)(\mu-2)\varphi + \l_0|\varphi|^{-n} \varphi=0.
 \ee
 One can see that, on orbits like (\ref{LC11}), the term $(|f|^{p-1}f)'$, which has been neglected in (\ref{OD-}), is much small than $f$ for all $p>1+ \frac n3$ that is true since always $p>1+n$.

We are interested in {\em periodic solutions} of (\ref{eqLC}),
which will give, according to (\ref{LC11}), oscillatory profiles
changing sign infinitely many times as $y \to y_0^-$, i.e., as   $\eta
\to -\infty$. Via (\ref{LC11}), periodic functions
$\varphi_*(\eta)$ establish the simplest oscillatory connections
with the interface points keeping the maximal regularity of the {\em envelope}:
 $$
 f(y) \sim (y_0-y)^{\frac 3n} \quad \mbox{for} \quad y \approx
 y_0^-,
 $$
 which represents the true scaling-invariant nature of the ODE
 (\ref{fff.1}).
 According to (\ref{LC11}), the regularity at $y=y_0$ then increases as $n \to
0^+$ forming, at $n=0$, analytic  solutions.  In
\cite[\S~6]{GBl6}, we present a discussion and references related
to the theory of periodic solutions of higher-order ODEs. Unlike
the fifth-order case in \cite[\S~4]{GBl6}, the ODE (\ref{eqLC}) is
of third order and can be reduced to a first-order ODE; see
\cite[\S~7.1]{Gl4}. Therefore the existence of a periodic solution
is not principally
 difficult, while the uniqueness (and the stability)
 are harder.
We expect, and this is confirmed by  numerics \cite{Gl4}, that
this limit cycle is ``almost" globally stable (note that all the
orbits of (\ref{eqLC}) are uniformly bounded, so a stable
attractor should be available, though sometimes 0 may have a
stable manifold, which was not observed) and is unique.

Increasing $n$ more, this periodic solution $\varphi_*(s)$
 is destroyed in a heteroclinic bifurcation at the point
 \cite[\S~7.2]{Gl4}
 \be
 \label{n**1}
  \mbox{$
 n_{\rm h}=1.758665... \,  
\big(\mbox{and}
 \,\,
  n_{\rm h} \in (\frac 32,n_{\rm +}), \,\, \mbox{where}
   \,\, n_{\rm +}= \frac 9{3+\sqrt 3}=1.9019238...\,, \mbox{\cite{PetI}}
   \big),
    $}
 \ee
 with a standard scenario of homoclinic/heteroclinic
bifurcations, \cite[Ch.~4]{Perko}. A rigorous
 justification of such non-local bifurcations  is an
  open problem.

Thus, for
 $n$  larger
than $\frac 32$, not all the solutions are oscillatory near the
interfaces. For $n \in (\frac 32, 3)$, there exists a
one-parametric bundle  of positive solutions with constant
equilibria $\varphi(\eta) \equiv  \pm \varphi_0$ given by
\be
 \label{Var55}
 \mbox{$
  \pm \varphi_0 = \pm  \bigl[ -\frac{\b y_0}{\mu(\mu-1)(\mu-2)}\bigr]^{\frac
1n}.
 $}
 \ee
 For matching
purposes, this is not enough and  the whole 2D asymptotic bundle
(\ref{LC11}) of oscillatory solutions has to be taken into
account, so, for the CP,
  the oscillatory behaviour is expected to remain generic
 for all $n \in (0,n_{\rm h})$ (and similar  to the linear case
$n=0$ with the interface at $y_0=\iy$).

 Thus \cite{Gl4}, it is key that, for $n \in (0, h_{\rm h})$,
 there exists a 2D bundle of asymptotic oscillatory orbits near the interface
 with the behaviour (here $y_0>0$ and $s_0 \in \re$ are
 parameters)
  \be
  \label{as1}
   \mbox{$
   f(y)=(y_0-y)^{\frac 3 n} \varphi_*(\ln(y_0-y)+s_0) + ... \quad
   \mbox{as}
   \quad y \to y_0^-.
   $}
   \ee


Another important question is the passage to the limit $n \to 0^+$
that shows convergence to solutions of the semilinear
Cahn--Hilliard equation. This is explained in detail in
\cite[\S~7.6]{Gl4}.  Various oscillatory sign change issues for
nonlinear degenerate higher-order PDEs of different types are
addressed in \cite[Ch.~3-5]{GSVR}.

\subsection{Continuous branches of similarity profiles for
$p=p_0$}

Again, without loss of generality, we treat the case $N=1$, where
the ODE is simpler and takes the form
 \be
\label{OF44}
 |f|^n f'''-\mbox{$ \frac 1{n+4}$} \, y f - (|f|^{n}
f^3)'= 0.
 \ee
 The origin of the existence of continuous
branches (parameterized by e.g. mass) is the fact that
the ODE (\ref{OF44}) is of third order. Thus the single symmetry
condition
\be
 \label{f12.0}
 f'(0)=0 \quad (f(0) \not = 0).
\ee
    is posed at the origin, while the shooting
bundle from the singularity point $y=y_0$ is 2D according to
(\ref{as1}). So one parameter is free. The existence by a shooting
approach is standard, so we refer to \cite{Gl4, EGW, GW2}
as a guide to such equations.

The numerical results below were mainly obtained using {\tt
Matlab}'s two-point boundary value problem collocation solver {\tt
bvp4c} (with default parameter values RelTol= $10^{-3}$, AbsTol=
$10^{-6}$). The standard regularization
 \be
 |f|^q \mapsto \big(\delta^2 + f^2\big)^{\frac q2}
 \label{eq:reg}
\ee
 was used with $\delta=10^{-6}$ (and sometimes up to $10^{-10}$ in order
 to see complicated and refined zero structure of solutions) for
$q=n$ and $\frac 2N$.

\begin{figure}
\centering
\includegraphics[scale=0.75]{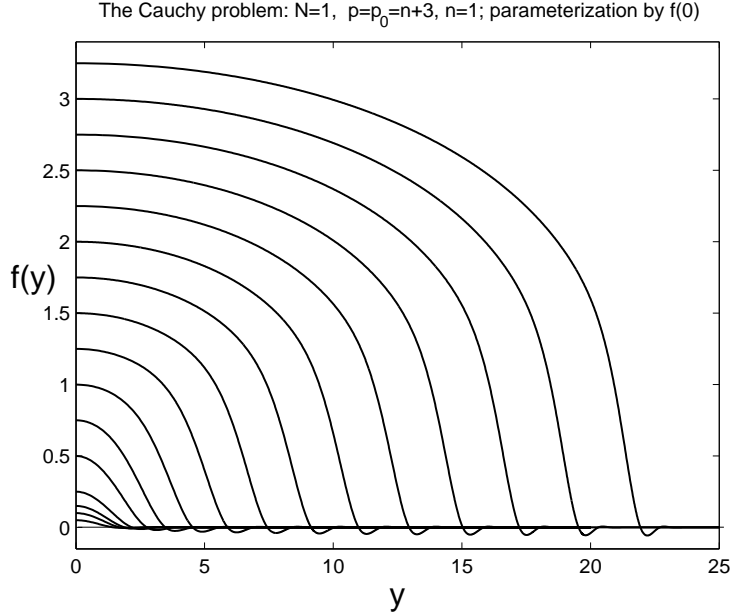}  
\vskip -.2cm \caption{\small Similarity profiles for the CP as
solutions of (\ref{OF44}), (\ref{f12.0}) for $N=n=1$, $p=n+3=4$.}
   \vskip -.1cm
 \label{FCP1}
\end{figure}

In Figure \ref{FCP1}, we present similarity profile for $n=N=1$ and
$p=p_0=n+3=4$, which are parameterized by the values at the origin
$f(0)$.

\begin{figure}
\centering \subfigure[ $f(y)$ for $n=0$ ]{
\includegraphics[scale=0.52]{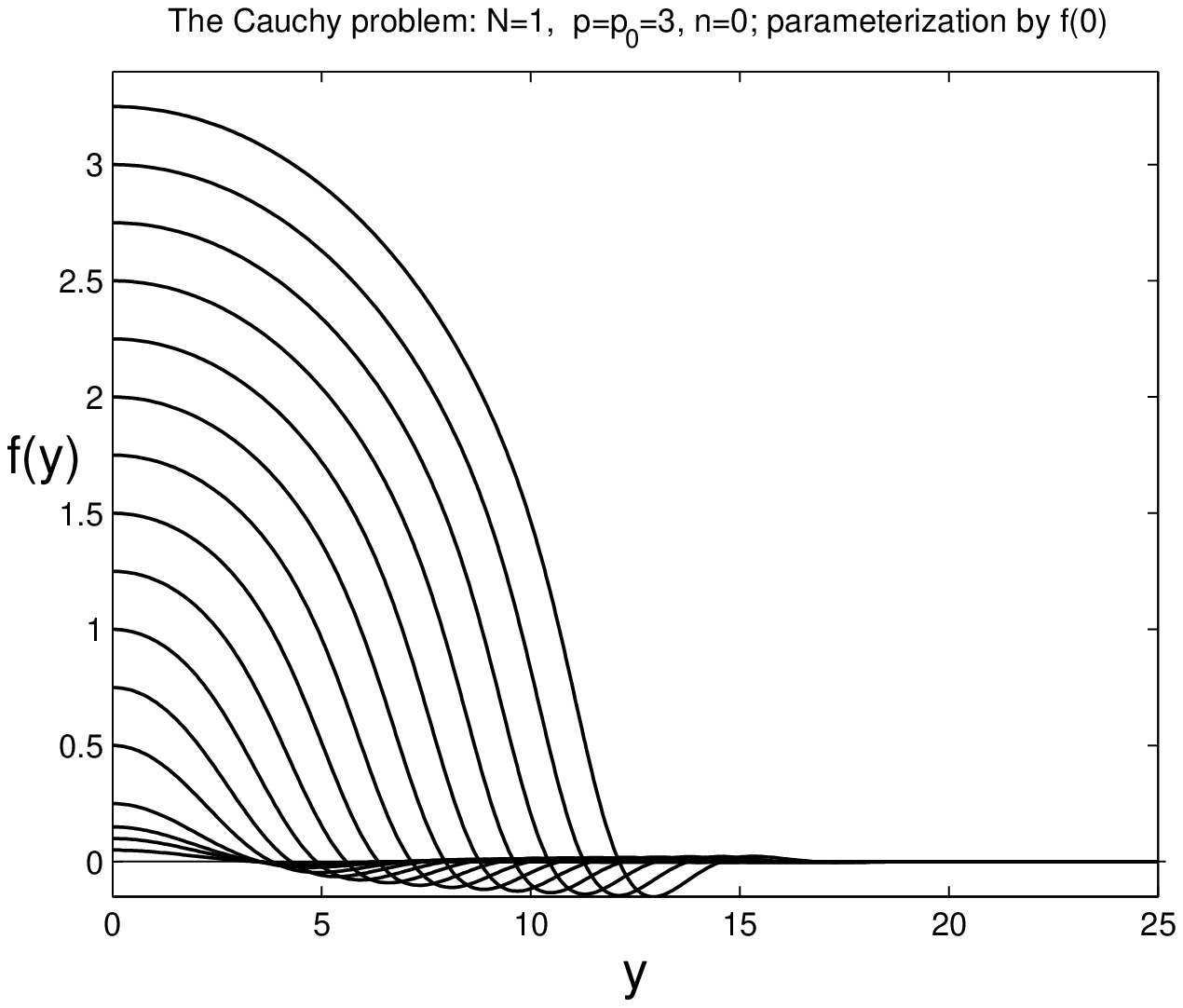}
} \subfigure[zero structure enlarged]{
\includegraphics[scale=0.52]{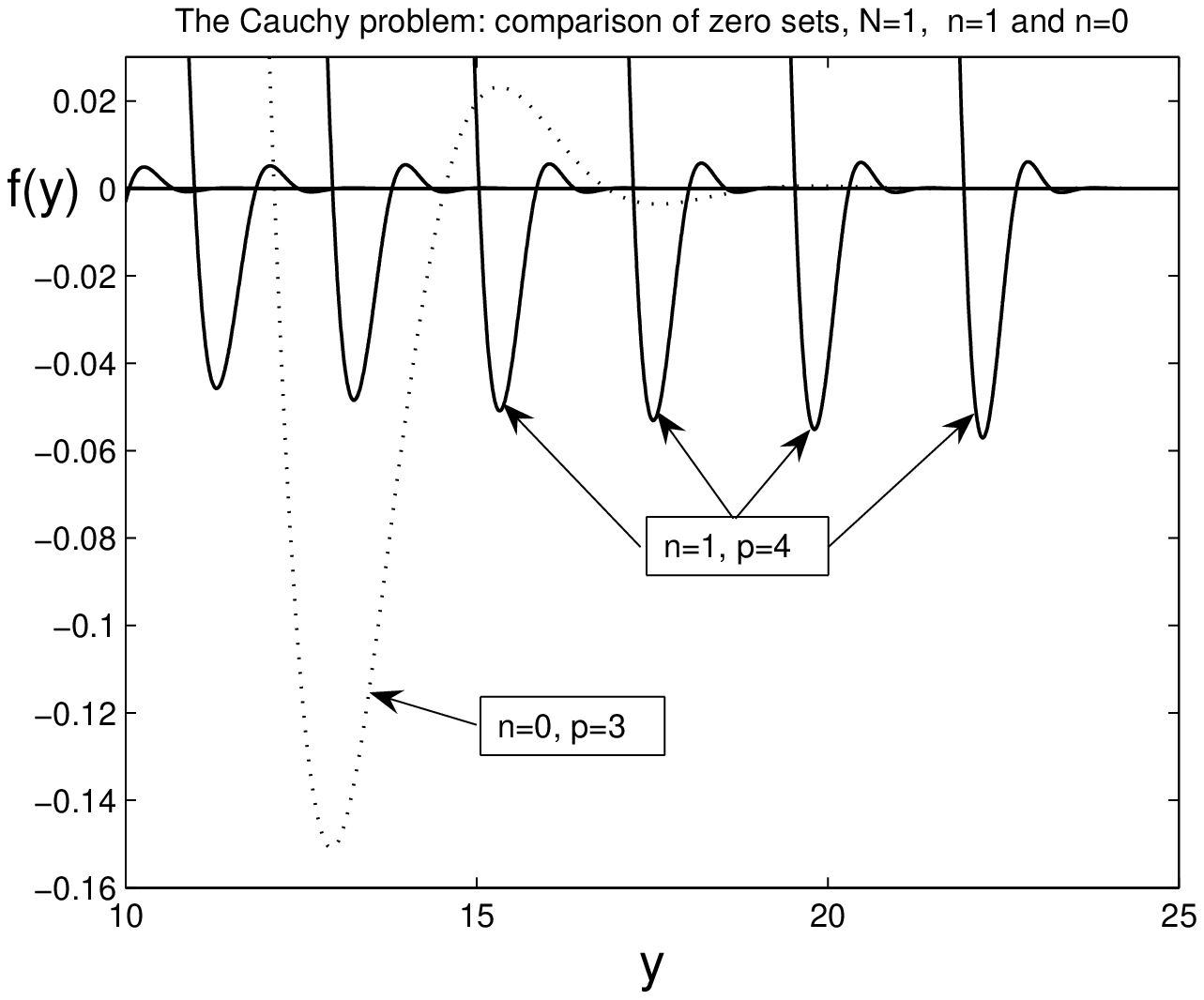}
} 
 \vskip -.2cm
\caption{\rm\small Similarity profiles for the CP satisfying
(\ref{OF44}), (\ref{f12.0}) for $N=1$, $n=0$ $p=3$; profiles (a)
and zero structure (b).}
 \label{CPn0}
 \vskip -0.3cm
\end{figure}

For comparison, in Figure \ref{CPn0}(a), we present similar
profiles
 for the semilinear case $n=0$, i.e., for the limit CH equation (\ref{CH11}),
  where $p=p_0=3$. Figure \ref{CPn0}(b)
 shows a clear difference of the ``tail" behaviour for $n=1$ (nonlinear oscillations (\ref{as1}))
  and $n=0$ (a linearized behaviour).

\begin{figure}
\centering
\includegraphics[scale=0.65]{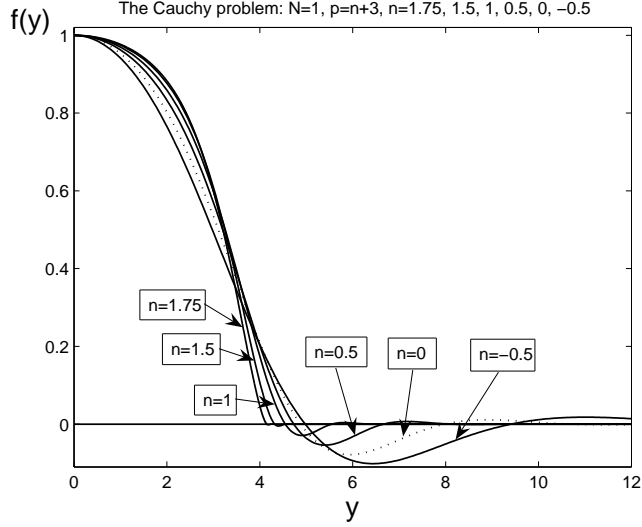}  
\vskip -.2cm \caption{\small Similarity profiles for the CP as
solutions of (\ref{OF44}), (\ref{f12.0}) for $N=1$, $p=n+3=4$ and
various $n \in [-\frac 12, 1\frac 34]$; parameterisation is
$f(0)=1$.}
   \vskip -.3cm
 \label{FCom1}
\end{figure}

In Figure \ref{FCom1}, we explain how the similarity profiles at
$p=p_0=n+3$ are deformed with $n$ starting again  from the
semilinear case $n=0$. It is seen that the profiles get thinner as
$n$ increases and also get less oscillatory near the interface.
 Note that the last case with the maximal $n=1.75$ (the CP profile must still
   change sign, while the FBP one does not \cite{Gl4}) is close to the
 heteroclinic value (\ref{n**1}), after which the similarity
 profiles are assumed to be finitely oscillatory, i.e., can have a
 finite number of sign changes near the interface (or no at all); see further comments
 in \cite[\S~9.4]{Gl4}. This finite oscillation phenomena lead to
 difficult and challenging numerical problems. Our numerics show that,
 for $n=1.75$, the CP profile still changes sign near the
 interface; see the next Figure \ref{FZero}.

In Figure \ref{FCom1}, we also included the case of a single
negative $n=- \frac 12$, which gives a standard source-type
profile but, of course, without finite interface, where $f(y)$ is
oscillatory as $y \to +\iy$. The structure of those oscillations at infinity is
difficult and different from already known ones, and is not studied here.

\begin{figure}
\centering
\includegraphics[scale=0.65]{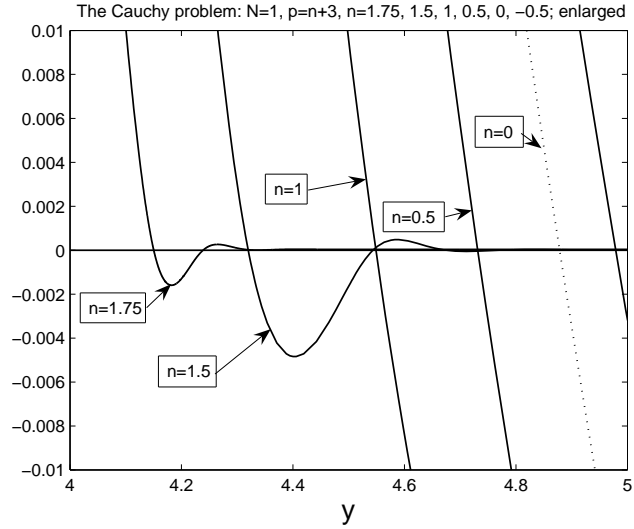}  
\vskip -.2cm \caption{\small Enlarged zero structure of VSS
profiles from Figure \ref{FCom1}; for $p=1.75$, $f(y)$ still
changes sign.}
   \vskip -.3cm
 \label{FZero}
\end{figure}


 \section{The Cauchy problem: countable family of source-type profiles via
 an
 $n$-branching approach}
 \label{CPTFEn}

For future convenience, we now need to postpone our study of the
original PDE (\ref{GPP}) and digress to the pure unperturbed TFE.
In general, construction of various oscillatory source-type
solutions of the CP for {\em pure TFE} (\ref{tf1})
 is a difficult nonlinear problem,
which is harder than that for the FBP.

\subsection{The linear case $n=0$: basics of Hermitian spectral theory}
 \label{S4.1}

On the other hand,
  for $n=0$, i.e., for the {bi-harmonic equation} (\ref{tf2}),
   the first
  such  profile exists, is unique (up to mass scaling), and is just the
   rescaled kernel $F(y)$ of the fundamental
  solution of (\ref{tf2}):
   \be
   \label{tf3}
    \mbox{$
     \begin{matrix}
    b(x,t)= t^{-\frac N4} F(y), \quad  y=x/t^{\frac 14}, \quad \mbox{where}
    \ssk\ssk\ssk\\
     {\bf B}F \equiv - \D^2 F + \frac 14 \, y \cdot \n F +
    \frac N4 \, F=0 \,\,\, \mbox{in} \,\,\, \ren, \quad \int F=1;
     \end{matrix}
     $}
     \ee
   see \cite[\S~4]{Bl4}.
   Moreover, there exists a countable set of eigenfunctions
   $\{\psi_\g, \, l= |\g| \ge 0\}$ of the corresponding rescaled non-self0-adjoint  operator
   ${\bf B}$ with the discrete spectrum \cite{Eg4}
    \be
    \label{tf4}
     \mbox{$
    \s({\bf B})= \big\{- \frac l 4, \,\,\, l =0,1,2,... \big\}.
     $}
     \ee
 The eigenfunctions are derivatives of the rescaled kernel $F$,
  \be
  \label{ps1}
   \mbox{$
    \psi_\g(y) = \frac {(-1)^{|\g|}}{\sqrt{ \g !}} \, D^\g F(y)
     \quad \mbox{for any multiindex} \,\,\, \g.
     $}
    \ee
 In addition, the adjoint operator has the same spectrum:
  \be
  \label{BB*}
   \mbox{$
    \BB^*= - \D^2 - \frac 14 \, y \cdot \n \quad \mbox{with}
    \quad
   \s(\BB^*)= \s({\bf B})= \big\{- \frac l 4, \,\,\, l =0,1,2,...
   \big\},
     $}
     \ee
 and a complete set of eigenfunctions $\{\psi_\g^*(y)\}$, which are
{\em generalized Hermite polynomials}. See more details on such a
Hermitian spectral theory of the pair $\{\BB,\BB^*\}$ in
\cite{Eg4, 2mSturm}.


\subsection{$n$-branching of similarity solutions}
 \label{S7.2}

We apply the $n$-branching approach from the linear case $n=0$ in
order to explain existence of a countable set of similarity
solutions of the TFE (\ref{tf1}). We will follow classic branching
theory in the case of non-analytic nonlinearities of finite
regularity;  see \cite[\S~27]{VainbergTr} and \cite[Ch.~8]{KrasZ}.

We look for solutions of (\ref{tf1}) with small $n>0$ in the
standard form
 \be
 \label{p1}
  \mbox{$
 u_\g(x,t)= t^{-\a} f(y), \quad y= x/t^{\b}, \quad \mbox{where}
  \quad \b= \frac {1-\a n}4,
   $}
   \ee
   where the multiindex $\g$ is used for numbering to be explained
   later on (in fact,  similar to the linear
   eigenfunctions (\ref{ps1})).
 Then $f=f_\g(y)$ solves the elliptic equation
 \be
 \label{TF990}
 \mbox{$
{\bf A}_n(f) \equiv  -\n \cdot (|f|^n \n \D f) + \b y \cdot \n f +
\a f=0 \quad \mbox{in}
 \quad
\ren.
 $}
  \ee
Note that, in general, for $l \ge 1$, we have to assume that
 \be
 \label{p2}
  \mbox{$
  \int f(y) \, {\mathrm d}y=0,
   $}
   \ee
so that the solutions (\ref{p1}) satisfy the mass conservation
condition
 \be
 \label{p3}
  \mbox{$
  \int u_\g(x,t) \, {\mathrm d}x \equiv 0 \quad (|\g| \ge 1).
   $}
   \ee
For $l=0$, where
 $
 \a= \b N$ and  $ \b = \frac 1{4+n N},
  $
 the assumption (\ref{p2}) is not necessary, since the PDE
 (\ref{TF990}) is fully divergent  and admits integration
 over $\ren$.


For small $n>0$ in (\ref{TF990}), we have
 \be
 \label{Ex.11}
  \mbox{$
 \b= \frac 14 -  \frac \a 4 \, n,
  $}
  \ee
  and use the following expansion:
   \be
   \label{Ex.6}
 |f|^n  = 1 + n  \ln |f| + o(n).
  \ee
Here,  (\ref{Ex.6}) should be understood in the weak sense, which
is necessary  for  using  in the equivalent integral equation; see
Proposition \ref{Pr.W} below.


Substituting expansions  (\ref{Ex.11}) and (\ref{Ex.6}) (still
completely formal) into (\ref{TF990}) yields
 \be
 \label{An.1}
  \mbox{$
 {\bf A}_n(f) \equiv {\bf B}f +
 \big(\a- \frac N4 \big)f +
  n {\mathcal L}(f) + o(n)=0,
   $}
  \ee
  with the  perturbation operator
 \be
 \label{LL.11}
 \mbox{$
 {\mathcal L}(f) =- \n \cdot \bigl( \ln|f| \n \D f \big) - \frac \a{4} \, y
\cdot \n f.
 $}
  \ee


We next describe the behaviour of solutions for small $n>0$ and
apply the classical Lyapunov-Schmidt method \cite[Ch.~8]{KrasZ} to
equation (\ref{An.1}).
 Recall that, in this linearized setting,
  we
 naturally arrive at the functional framework that is suitable for
  the linear operator ${\bf B}$, i.e., it is $L^2_\rho(\ren)$, with the
  domain $H^4_\rho(\ren)$, etc., and a similar setting
   for the adjoint operator ${\bf B}^*$; see above and further details in
\cite{Eg4}.

Therefore, for $n=0$, we have to study branching of a nonlinear
eigenfunction from the linear one, where $f$ is a certain
nontrivial finite linear combination of eigenfunctions from a
given eigenspace with fixed $\l_\g= -\frac{|\g|}4 \equiv -\frac l4$, i.e.,
 \be
 \label{a4}
  \mbox{$
 f= \phi_l = \sum_{|\g|=l} C_\g \psi_\g \quad (\not = 0).
   $}
   \ee
Those eigenfunctions are just derivatives (\ref{ps1}) of the
analytic radially symmetric rescaled kernel $F(|y|)$ of the
fundamental solution (\ref{tf3}). Therefore, the nodal (zero) set
of such $f$ in (\ref{a4}) is well understood and consists of a
countable set of isolated sufficiently smooth hypersurfaces which
can concentrate as $y \to \infty$, where
 \be
 \label{FF1}
 \phi_l(y) \to 0 \quad \mbox{as \,\, $y \to \infty$ uniformly and
 exponentially fast.}
  \ee

  Hence, returning to the key limit (\ref{Ex.6}),
 we have:

\begin{proposition}
 \label{Pr.W}
 For a function $f$ given by $(\ref{a4})$,
 in the sense of distributions,
 \be
 \label{Is.1}
 \mbox{$
  \frac 1n \, (|f|^n -1) \rightharpoonup \ln |f| \quad \mbox{as} \quad
  n \to 0^+.
   $}
   \ee
   \end{proposition}


 According to (\ref{Is.1}), analyzing the integral equation for $f$,
  we can use
   the fact that,
 for any function $\phi \in {\mathcal L}$
   (and $\phi \in
 C_0(\ren)$),
 \be
 \label{Is.2}
  \mbox{$
  \int_{\ren} (|f(y)|^n -1) \phi(y) \, {\mathrm d}y =
n \big[\int_{\ren} \ln|f(y)|\phi(y) \, {\mathrm d}y + o(1) \big].
 $}
 \ee

 It follows from (\ref{An.1}) that branching
is possible under the following non-trivial kernel assumption: for
$n=0$,
 \be
 \label{k1}
  \mbox{$
\a- \frac N4= - \l_l=  \frac l4
 \quad \Longrightarrow \quad \a_l(0)= \frac{N+l}4, \quad l \ge 0.
  $}
   \ee
This gives the countable sequence of critical exponents
$\{\a_l(n),\b_l(n)\}$ (to be determined) of the similarity
patterns (\ref{p1}) of the TFE for small $n>0$.

 By   \cite[Lemma~4.1]{Bl4}, the kernel of
the linearized operator
 $$
 E_0={\rm ker\,}({\bf B}- \l_l I)={\rm
Span\,}\{\psi_\b, \, |\b|=l\}
 $$
  is finite-dimensional.
Hence, denoting by $E_1$ the complementary (orthogonal to $E_0$)
invariant subspace, we set
  \be
  \label{kk1}
   \mbox{$
 f= \phi_l+V_1, \quad \mbox{where} \,\,\, \phi_l \in E_0
  \,\,\, \mbox{and} \,\,\,V_1= \sum_{|\g| >
l} c_\g \psi_\g \in E_1.
 $}
   \ee
  According to the known spectral
properties of operator ${\bf B}$, we define $P_0$ and $P_1$,
$P_0+P_1=I$, to be projections onto $E_0$ and $E_1$ respectively.
We also introduce a perturbation of the parameter $\a$ by setting
 \be
 \label{a1}
 \a_l(n)= \a_l(0)+ \d, \quad \mbox{with}
 \quad \d=\d(n).
  \ee
This perturbation $\d$ is obtained from the orthogonality
condition by substituting into (\ref{An.1}) and multiplying by
$\psi_\g^*$. This gives
 \be
 \label{a2}
 \d(n) = c_l n + o(n),
  \ee
  where $c_l$ is obtained from the system
  \be
  \label{a3}
 \langle {\mathcal L}(\phi_l), \psi_\g^* \rangle = c_l, \quad
 |\g|=l.
  \ee
  Since according to (\ref{kk1}),  $\phi_l$ is given by (\ref{a4}),
 (\ref{a3}) is  an algebraic system for unknowns $\{C_\g\}$ and $c_l$.
 It can be  solved, for instance, in the radial geometry
 and in some other cases (including those where the dimension
 of the kernel is odd; even dimensions are known to need additional treatment);
 see more details in \cite[App.~A]{GHCo}. However, the total number of solutions of the non-variational system (\ref{a3}) remains unclear.


Finally, 
 setting
 \be
 \label{a5}
 V_1 = n Y + o(n),
  \ee
 we obtain,  passing to the limit $n \to 0^+$,
 the following equation for $Y$:
 \be
 \label{YY.11}
{\bf B}Y =- c_l \phi_l +{\mathcal L}(\phi_l).
 \ee
 By  Fredholm's theory, in view of the orthogonality, it admits a unique solution
  $Y \in E_1$.

  \ssk

In general, the above analysis shows that, up to solvability of
the nonlinear algebraic systems, the TFE admits a countable set of
different source-type similarity solutions (\ref{p1}) at least for
small $n>0$, where the parameters $\a_l(n)$ are given by
 \be
 \label{par1}
  \mbox{$
  \a_l(n)=  \frac {N+l}4 + c_l n + o(n) \quad \mbox{as}
   \quad n \to 0^+; \quad l=0,1,2,... \, .
    $}
    \ee
 At $n=0$, these solutions are originated from
suitable eigenfunctions of the linear operator in (\ref{tf3}). The
global extensions of these $n$-branches of similarity solutions
for larger $n>0$ represent a difficult open problem, to be treated
numerically later on.

\subsection{Nonlinear eigenfunctions of the TFE in one dimension}

We consider the Cauchy problem for the 1D TFE with
continuous compactly supported initial data,
 \be
 \label{eq1}
  u_t=- (|u|^n u_{xxx})_x \quad \mbox{in}
  \quad \re \times \re_+, \quad u(x,0)=u_0(x) \in C_0(\re).
  \ee
Then, for $N=1$, the {\em nonlinear eigenvalue problem} for the
elliptic equation (\ref{TF990}) is formulated as follows:
 \be
 \label{nl1}
  \fbox{$
 -(|f|^n f''')' + \frac {1-\a n}4 \, y f' + \a f=0 \quad \mbox{in}
 \quad \re, \quad f(y) \not \equiv 0, \quad f \in C_0(\re).
  $}
  \ee

Actually, (\ref{nl1}) is also about self-similarity of
{\em second kind}, where the desired  set of parameters
(nonlinear eigenvalues) $\{\a_l(n), \, l \ge 0\}$
 is obtained not by a pure dimensional analysis, but via
 solvability of a nonlinear ODE in a given functional class
 $C_0(\re)$ of compactly supported functions satisfying the
 condition of maximal regularity. The term {\em similarity of the
 second type} was  introduced by
 Ya.B.~Zel'dovich in 1956,
 \cite{Zel56}.

  Note that,
for $n=0$, (\ref{nl1}) in $L^2_\rho(\re)$, where we replace the
last condition by
 $$
  f \in L^2_\rho(\re),
  $$
is a standard {\em linear eigenvalue problem} for a non
self-adjoint operator with the point spectrum (\ref{tf4}) and
complete-closed set of eigenfunctions $\{\psi_\b\}$ given in
(\ref{ps1}), \cite{Eg4}.

\begin{figure}
\centering \subfigure[eigenfunctions]{
\includegraphics[scale=0.52]{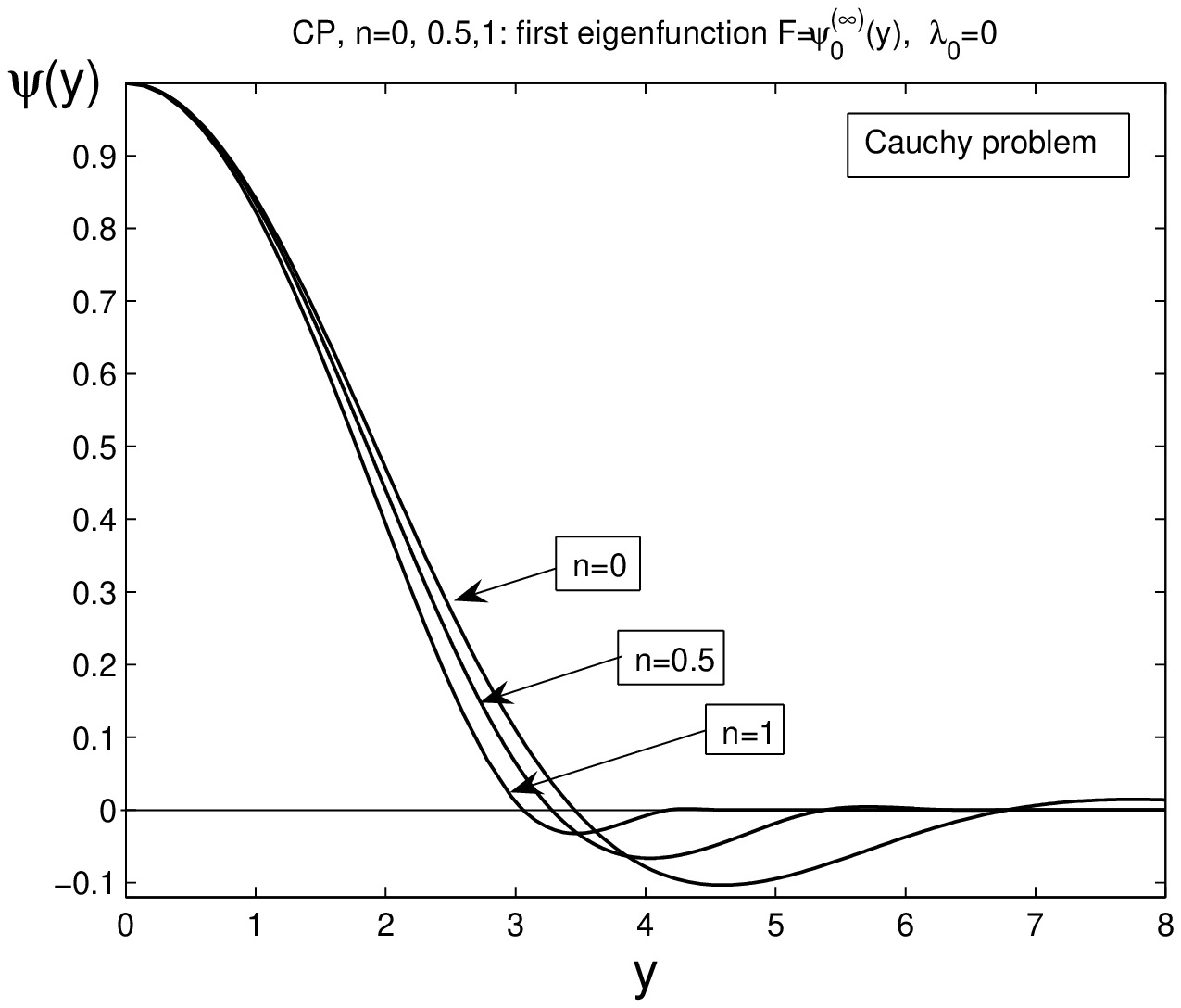}
} \subfigure[zeros in log-scale]{
\includegraphics[scale=0.52]{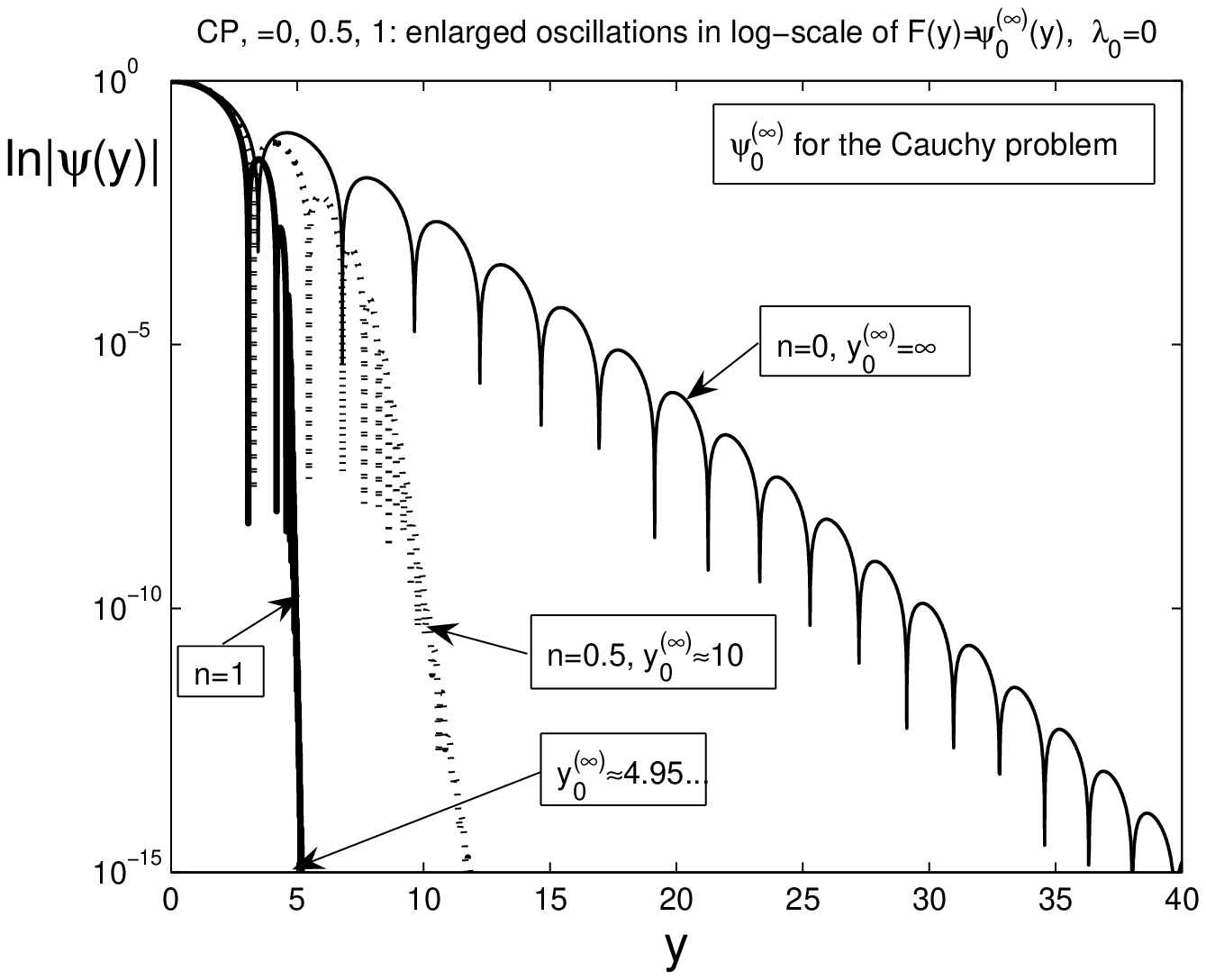}
}
 \vskip -.4cm
\caption{\rm\small The first eigenfunction of (\ref{nl1}),
(\ref{nn1}) for $n=0$, $0.5$, and $1$.}
 \label{FFF01}
\end{figure}

The first nonlinear eigenvalue-eigenfunction pair $\{F_0, \,
\a_0\}$ of (\ref{nl1}) was proved to exist for $n \in (0,1]$; see
\cite[\S~9]{Gl4}. In this case, the first eigenvalue is
 \be
 \label{nn1}
  \mbox{$
  \a_0(n) =\frac {N}{4+n N}\big|_{N=1}= \frac 1{4+n}.
   $}
  \ee
It turns out that, in view of the highly oscillatory nature of those
profiles near interfaces, even identifying the position of
interfaces numerically, is not an easy problem.
 Therefore, we begin with Figure \ref{FFF01}, where the first even
 nonlinear eigenfunction is presented for $n=0$, $\frac 12$, and
 $1$. A careful study of their zero structure in the log-scale in
 (b) allows us to find an approximate and rather rough interface
 location, according to the expansion (\ref{LC11}), which yields
 \be
 \label{log1}
  \mbox{$
   \ln |f(y)| = \frac 3n \, \ln(y_0-y) + \ln|\varphi(\ln(y_0-y)|+... \, .
    $}
    \ee
 For $n=0$, the expansion is exponential (cf. (\ref{.7}) below) and is entirely different, which is seen in
  (b); recall the regularization  such as in (\ref{eq:reg}) eventually
  entering the expansion for $|f|$ very small.

Thus, in what follows, we use the  parameterisation:
 \be
 \label{hh1}
 F_l(0)=1, \,\,\, l=0,2,4,... \, ; \quad F_l'(0)=1, \,\,\,
 l=1,3,... \, .
  \ee
We then obtain that,
respectively,
 \be
 \label{hh2}
 \a_1=0.2534... \, , \quad \a_2=0.320... \, \quad(n=1).
  \ee
In addition, numerics show that $\a_3 \approx 0.38$ for $n=1$.

\begin{figure}[htp]

\includegraphics[width=7.5cm]{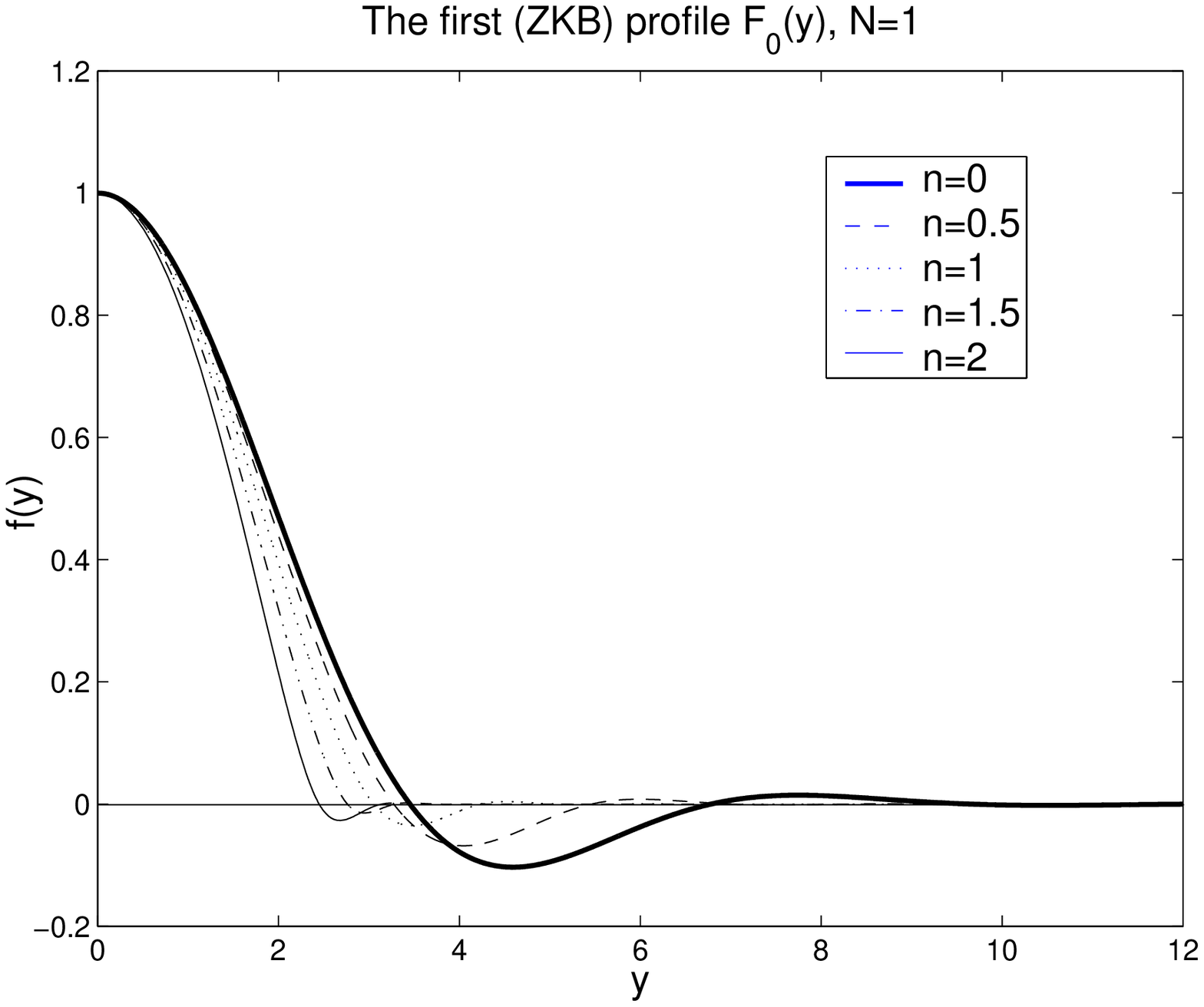}
\includegraphics[width=7.3cm]{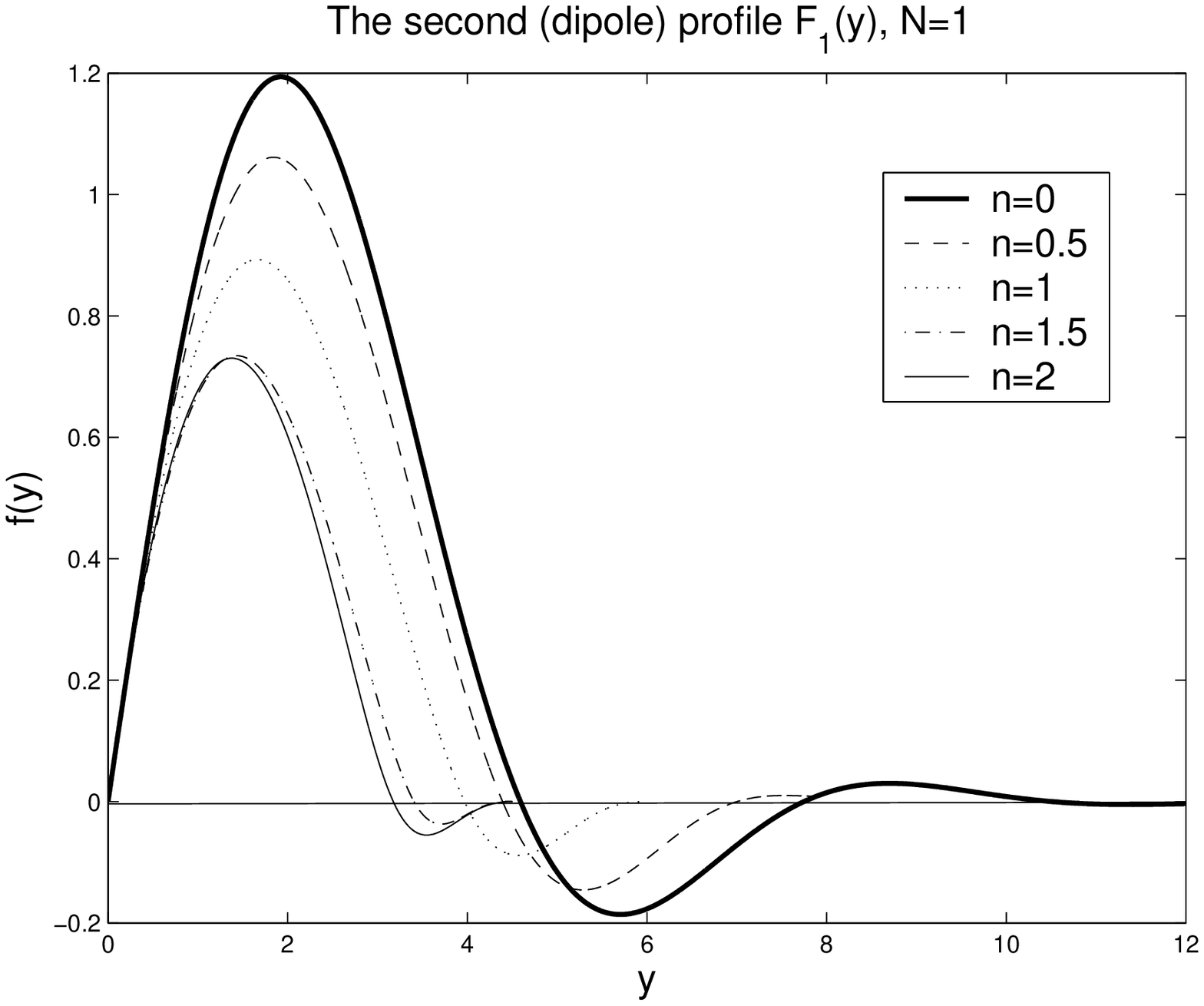}
\includegraphics[width=7.3cm]{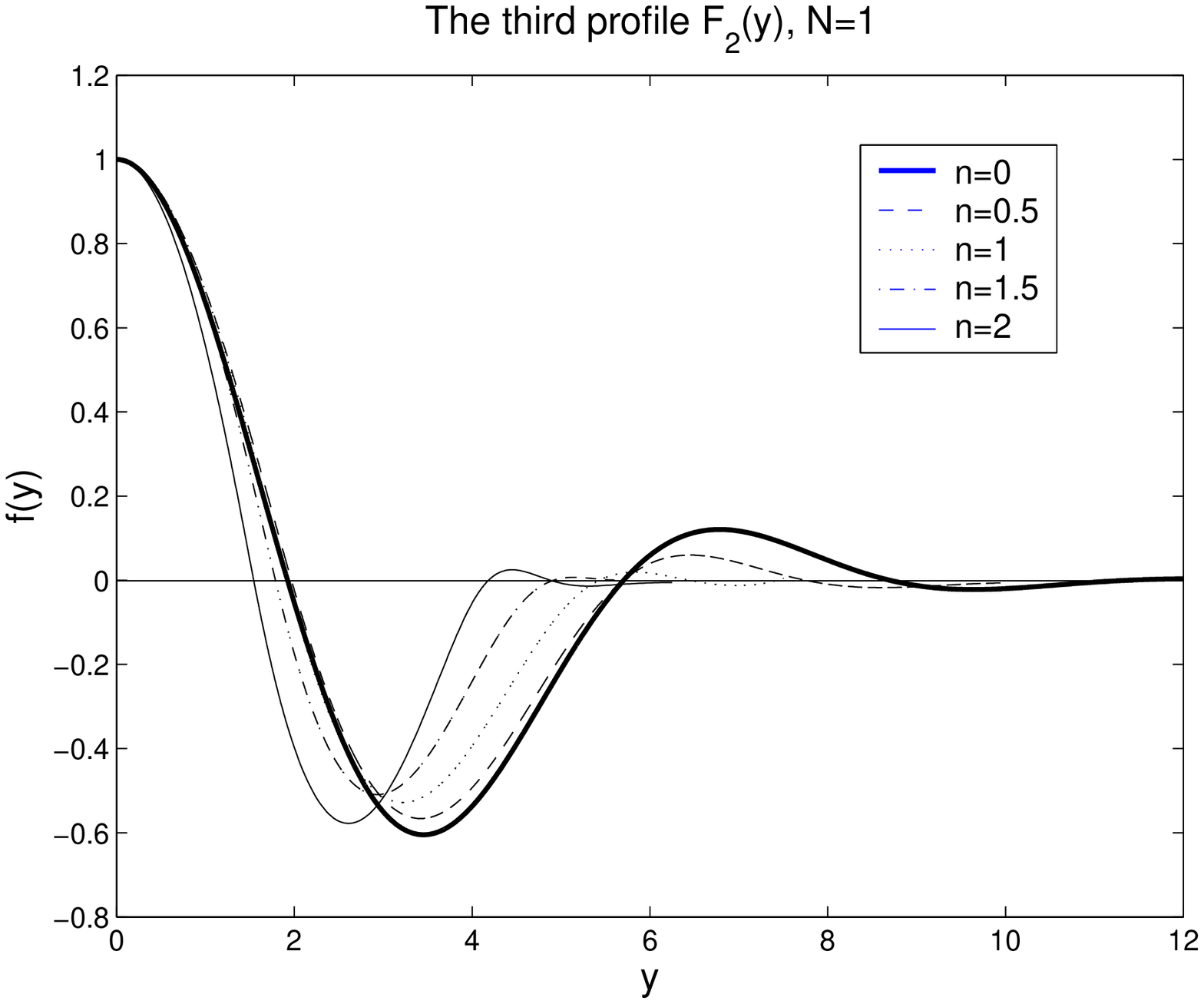}
\includegraphics[width=7.3cm]{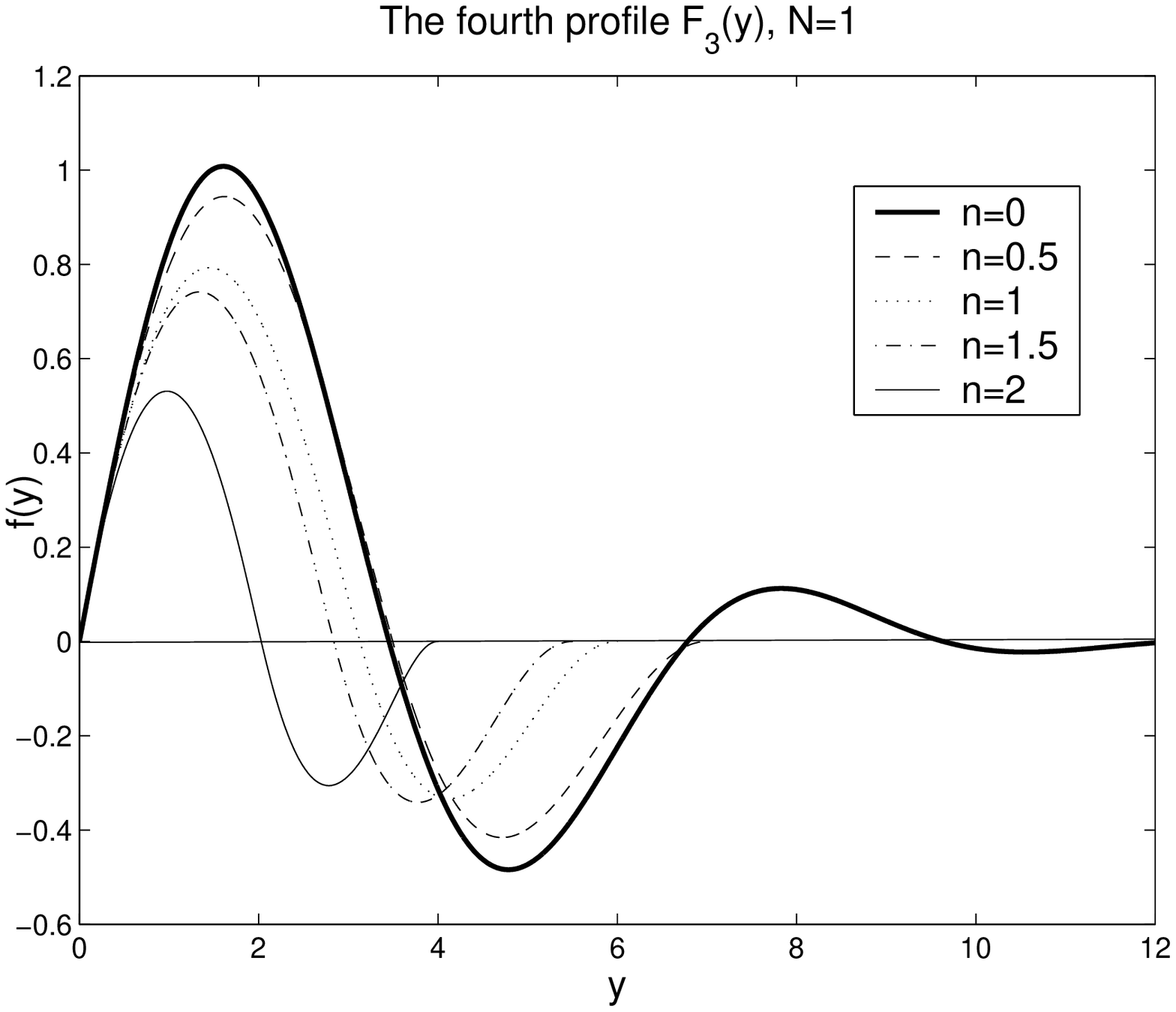}

\vskip 0cm
\caption{ \small Illustrative numerical solutions for the first
 four nonlinear eigenfunctions $F_{l}(y)$, $l=0,1,2,3,$ shown for selected $n$ in one
  dimension $N=1$.
The corresponding behaviour of the eigenvalues are shown in Figure \ref{Fn1a}.
 The regularisation (\ref{eq:reg}) with $\delta=10^{-2}$ was used
in the numerical shooting scheme, where the even numbered profiles
 satisfy $F_{l}(0)=1,F_{l}'(0)=F_{l}'''(0)=0,$ whilst the odd numbered
  profiles have $F_{l}'(0)=1,F_{l}(0)=F_{l}''(0)=0.$ }
 \label{num1N}
\end{figure}

The results of an accurate numerical study of the first four
nonlinear eigenfunctions for $N=1$ and various $n \in [0,2]$ are
presented in Figure \ref{num1N}. Further eigenfunctions are very
difficult to obtain numerically, to say nothing about an
analytical proof of their existence.

\begin{figure}
\centering
\includegraphics[scale=0.65]{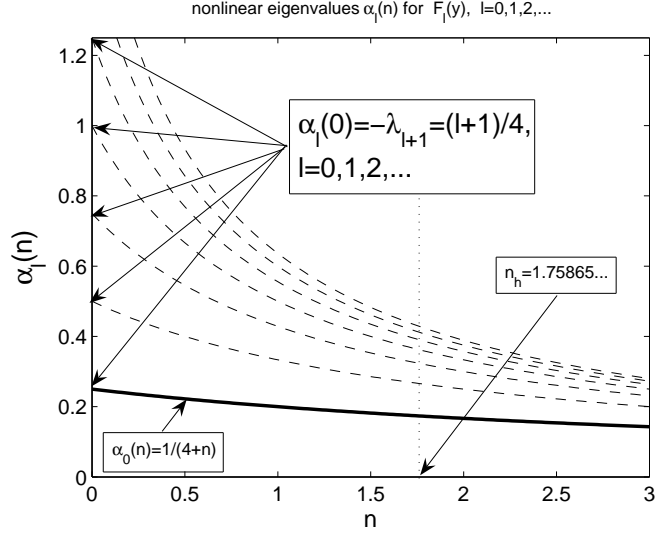}  
\vskip -.5cm \caption{\small Formal $n$-branches of nonlinear eigenvalues
$\a_l(n)$ of (\ref{nl1}).}
   \vskip -.3cm
 \label{Fn1}
\end{figure}

In Figure \ref{Fn1}, we show a formal schematic behaviour of the
nonlinear eigenvalues $\a_l(n)$ of (\ref{nl1}). The first
$n$-branch, according to (\ref{nn1}), has the explicit form
 $$
  \mbox{$
 \a_0(n)= \frac {1}{4+n}, \quad n \in [0,3).
 $}
  $$
Other $n$-branches in Figure \ref{Fn1} are not explicit and are
hypothetical. According to the $n$-branching approach, all these
branches originate at the eigenvalues of the linear problem
(\ref{tf4}), i.e.,
 \be
 \label{ss1}
  \mbox{$
  \a_l(0)= - \l_{l+1} = \frac {l+1}4 \quad \mbox{for}
  \quad l=0,1,2,... \, ,
   $}
    \ee
and moreover, after scaling, we may assume that, at $n=0$, the
similarity profiles $F_l(y)$ coincide with the eigenfunctions
(\ref{ps1}), and hence by continuity  mimic their geometric shapes
for $n>0$.

\begin{figure}
\centering
\includegraphics[scale=0.65]{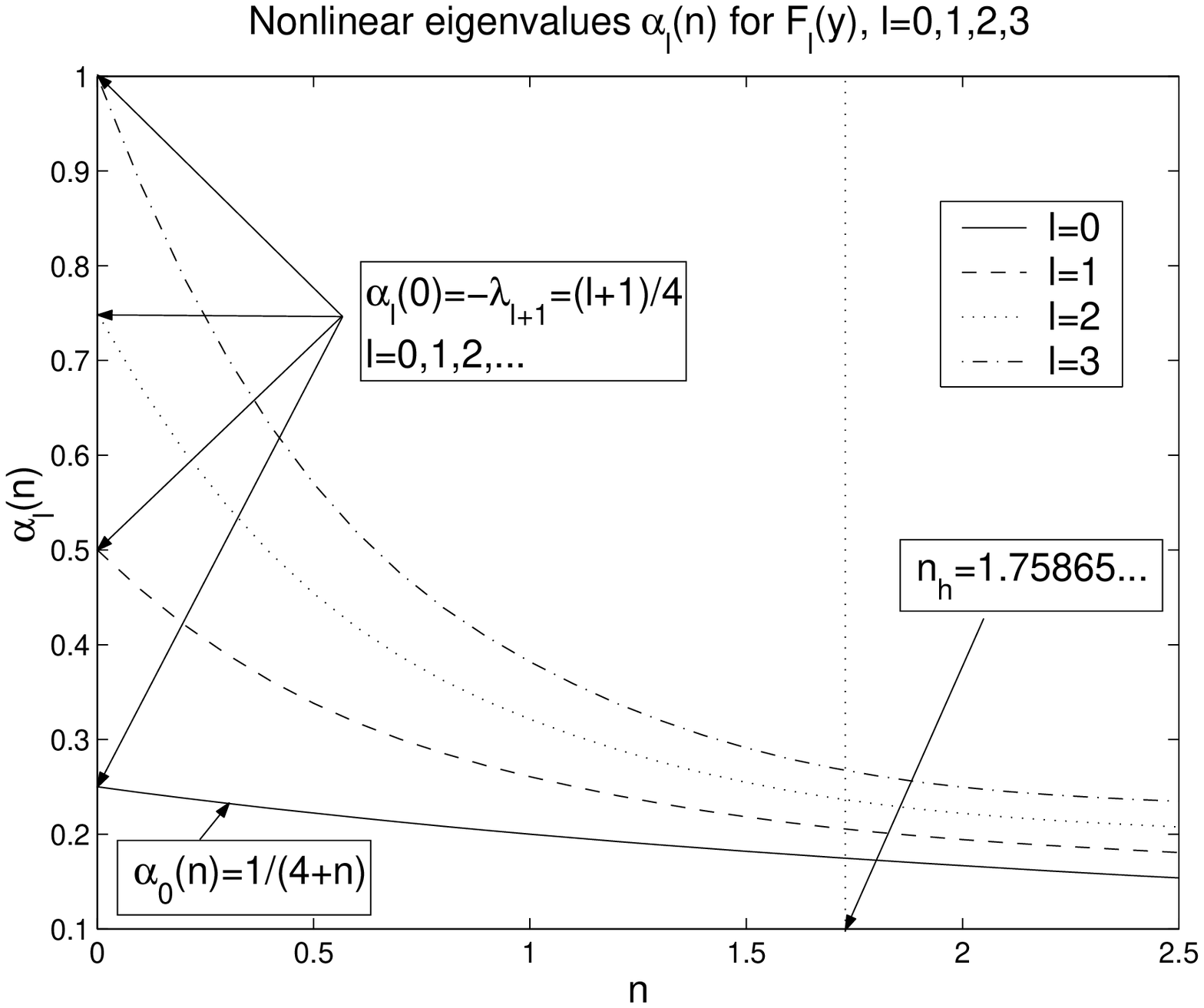}  
\vskip 0cm \caption{\small The actual $n$-branches of the first
four  nonlinear eigenvalues $\alpha_{l}(n)$, $l=0,1,2,3$,
constructed numerically. }
   \vskip -.5cm
 \label{Fn1a}
\end{figure}

Figure \ref{Fn1a} shows the actual numerical construction of first four $n$-branches,
and even these involve technical difficulties. Note that, at the critical heteroclinic bifurcation
value (\ref{n**1}), the similarity profiles $F_l(y)$ are supposed
to loose their oscillatory behaviour at the interface and become
finite oscillatory for $n > n_{\rm h}$ (or even non-oscillatory at
all); see \cite[\S~7.2]{Gl4}.

Analytical difficulties for the eigenvalue problem (\ref{nl1})
begin already with $l=1$, i.e., with the dipole profile $F_1(y)$.
This study has a well-developed history (see \cite{BerHK00,
BHK, BW06} and references therein), but still there are no
definite results of existence and uniqueness of $F_1$ in both the FBP
and Cauchy problem settings.


We end this discussion with the following:

\ssk

\noi {\bf Conjecture \ref{CPTFEn}.1.} (i) {\em For any $n \in
(0,n_{\rm h})$,
 the nonlinear eigenvalue problem $(\ref{nl1})$ admits
 a countable set of sufficiently smooth solutions of \underline{maximal
regularity}\footnote{More details on this are given in \cite{Gl4}.}
 \be
 \label{nl2}
 \Phi=\{F_l(y), \,\, l=0,1,2,...\},
  \ee
  where the nonlinear eigenvalues $\{\a_l\}$ form a strictly
  increasing sequence and
   \be
   \label{nl3}
    \mbox{$
 \a_l \to \frac 1n \quad \mbox{as} \quad l \to \infty.
  $}
   \ee
 }
 \noi (ii) {\em The eigenfunction subset $(\ref{nl2})$ is
 evolutionary complete in $C_0(\re)$ for the TFE $(\ref{eq1})$,
 i.e.,
 for any $u_0 \not =0$, there exists a finite $l
\ge 0$ and a constant $b=b(u_0) \not = 0$ such that
 \be
 \label{uu1}
 u(x,t) = t^{-\a_l} \big[ b F_l(x/t^{\b_l}|b|^{n/4})+ o(1)\big]
 \quad \mbox{as} \quad t \to \infty.
 \ee
  }

We expect that an analogous countable set of radially symmetric
similarity solutions exists for the TFE (\ref{tf1}) in any
dimension $N \ge 2$, though numerical calculations become much
more difficult than for $N=1$. Moreover, the branching approach
 in Section \ref{S7.2} shows that there many other non-radial
 similarity solutions that have a more complicated geometry but which for
 small $n>0$ mimic the eigenfunctions (\ref{ps1}).

 \ssk

 The evolution completeness of nonlinear eigenfunctions is known
 rigorously for the PME
  \be
  \label{pp1NN}
   u_t= (|u|^n u)_{xx} \quad (n>0)
    \ee
    in a bounded interval
 \cite{CompG}; see also \cite{GHCo} for results
in $\ren$ for initial data $u_0 \in C_0(\ren)$ and on an
$n$-branching technique. Existence of a  countable set of radial
similarity solutions of (\ref{pp1NN}) in $\ren \times \re_+$ was
proved by Hulshof \cite{H}.

\section{The Cauchy problem: on  nonlinear $p$-bifurcations}
 \label{CPpn}

\subsection{Semilinear Cahn--Hilliard equation: countable set of
critical exponents}

In order to explain the essence of the nonlinear
bifurcation analysis, we first digress to the stable CH equation
(\ref{CH11}), for which the analysis is much simpler; cf.
\cite{EGW, GW2}.

Namely, studying the behaviour as $t \to +\iy$, we perform the
standard scaling
 \be
 \label{M1}
 u(x,t)= (1+t)^{-\frac 1{2(p-1)}} v(y,\t), \quad y=x/(1+t)^{\frac
 14}, \quad \t= \ln(1+t),
  \ee
 where $v(y,\t)$ solves the following rescaled equation:
  \be
   \label{M2}
    \mbox{$
   v_\t = - \D^2 v + \frac 14 \, y \cdot \n v+ \frac 1{2(p-1)} \,
   v + \D (|v|^{p-1}v) \equiv {\bf B}v + c_0 v +\D (|v|^{p-1}v).
 $}
 \ee
 It then follows from (\ref{tf4})  that a centre manifold
 behaviour is formally possible in the critical cases (\ref{Crp}) only:
  \be
  \label{M3}
   \mbox{$
c_0= \frac 1{2(p-1)}- \frac N4= \frac l4 \quad \Longrightarrow
\quad p=p_l=  1 + \frac 2{N+l} \quad (l \ge 0).
 $}
  \ee
Checking the necessary condition of such a centre manifold
behaviour
 and looking, say, for a solution moving
along the centre
eigenspace, 
 \be
 \label{M4}
  \mbox{$
  v(y,\t) = a_\g(\t) \psi_\g(y) + w, \quad w \bot \psi_\g, \quad
  \sup_y |w(y,\t)|=o(a_\g(\t))\,\,\,\mbox{as} \,\, \t \to \iy,
   $}
   \ee
   and substituting into (\ref{M2}) yields on multiplication by
   $\psi_\g^*$ (see \cite{Eg4} for details)
   \be
   \label{M5}
    \mbox{$
     \dot \a_\g(\t)= \mu_\g |a_\g|^{p-1} a_\g + ..., \quad \mbox{where}
     \quad \mu_\g= \langle  \D|\psi_\g|^{p-1}\psi_\g,\psi_\g^*
     \rangle \equiv \langle  |\psi_\g|^{p-1}\psi_\g,\D\psi_\g^*
     \rangle.
      $}
      \ee
 Since $\psi_\g^*(y)$ is a $\g$-degree polynomial \cite{Eg4}, we
 then conclude that the necessary condition of existence of such a
 centre subspace behaviour is as follows:
  \be
  \label{M6}
    \mu_\g \not = 0 \quad \mbox{at least, for $|\g| \ge 2$.}
   \ee
 Note that, in the limit $p \to 1$, the following holds:
  \be
  \label{M7}
  \mu_\g =1 \quad \mbox{by bi-orthonormality of eigenfunctions,}
   \ee
 so that (\ref{M6}) is true for $l \gg 1$ by continuity of the integral relative to the parameter $p$.

Eventually, the centre subspace behaviour (\ref{M5}) generates the
following asymptotic patterns for the CH (\ref{CH11}):
 \be
 \label{M8}
  \mbox{$
  u_\g(x,t) \sim C_\g (t \ln^2 t)^{-\frac{N+l}4} \psi_\g\big( \frac
  y{t^{1/4}}\big)+... \quad \mbox{as} \quad t \to \iy \quad
  (l=|\g| \ge 2),
   $}
   \ee
 where constants $C_\g$ are independent of initial data $u_0$.

\subsection{Local bifurcations from $p_l$}
 \label{S5.2}

We now return to the VSSs of the TFE with the stable PME term
(\ref{GPP}) and perform a formal {\em nonlinear} version of a
$p$-bifurcation (branching) analysis for $n>0$. As usual,
according to classic branching theory \cite{KrasZ, VainbergTr}, a
justification (if any) is performed for the equivalent quasilinear
integral equation with compact operators. For simplicity, we
present computations for the differential setting.

Thus, we consider the elliptic PDE (\ref{RescProf}). The critical
exponents $\{p_l\}$ are then determined from the equality ({\em
q.v.} (\ref{M3}))
 \be
 \label{pp1}
  \mbox{$
 \a \equiv  \frac 1{2p_l(n)-(n+2)}= \a_l(n)
 \quad \Longrightarrow \quad p_l(n)= \frac {n+2}2 + \frac 1{2
 \a_l(n)} \quad (l \ge 0).
  $}
  \ee
In particular, for the semilinear case $n=0$, we have $\a_l(0)=
\frac {N+l}4$ from (\ref{k1}), so that (\ref{pp1}) leads to
(\ref{M3}), i.e., to
 the known sequence of critical exponents (\ref{Crp}).

We next use an expansion relative to the small parameter $\e= p_0-p$, i.e., as $\e \to 0$,
 $$
 \mbox{$
 \begin{matrix}
 \a= \frac 1{2p_l-(n+2)-2 \e}= \a_l + 2 \a_l^2 \e +... \, ,
 \ssk\ssk\\
 \b= \frac{1-n \a_l}4 + c_l \e + ... \, , \,\,\, c_l=
 \frac{1-n \a_l}4\big[ n+2 + \frac 1{\a_l} - \frac 2{1-n
 \a_l}\big].
  \end{matrix}
   $}
    $$
Substituting these expansions and the last one in (\ref{Ex.11})
into (\ref{RescProf}) and performing the same standard
linearization yields
  \be
 \label{eq4}
  \mbox{$
 {\bf A}_n(f) + \D(|f|^{p_l} f) + \e
 \big[ - \D(|f|^{p_l}f \ln |f|) + {\mathcal L}_1 F \big] + O(\e^2)=0,
 $}
   \ee
 $$
  \mbox{$
  \mbox{where} \quad {\mathcal L}_1=
 c_l y \cdot \nabla + 2 \a_l^2 I
  $}
  $$
  is a linear operator, and ${\bf A}_n$ is the rescaled operator (\ref{TF990}) of
the pure TFE with the parameter $\a= \a_l(n)$ (an eigenvalue), for which there exists
the corresponding similarity profile $F_l(y)$ (the nonlinear eigenfunction).
 The fact that the operator
 ${\bf A}_n$ with $\a=\a_l$ in (\ref{eq4}) occurs in the
rescaled
 pure TFE correctly describes the essence of a
 {\em ``nonlinear bifurcation phenomenon"} to be revealed.

To this end, we use the additional invariant scaling of the operator
${\bf A}_n$ by setting
  \be
 \label{ep1}
 f(y)=b  F(y/b^{\frac n4}) \quad (b>0),
  \ee
 where $b=b(\e)>0$ is a small parameter satisfying
  \be
  \label{bb1}
   b(\e) \to 0 \quad \mbox{as} \quad  \e \to 0,
   \ee
    to be determined.
 Substituting (\ref{ep1}) into (\ref{eq4}) and omitting
 all higher-order terms (including the one with the logarithmic multiplier $\ln
|b(\e)|$)
  yields
  \be
 \label{ep2}
 {\bf A}_n(F) + b^{p_l- \frac n2} \D(|F|^{p_l} F)+
  \e {\mathcal L}_1  F =0.
  \ee

 Finally, we perform linearization about the nonlinear
 eigenfunction $F_l(y)$ by setting
  $$
  F=F_l+Y.
   $$
    This yields the following linear
 non-homogeneous problem:
   \be
  \label{ep3}
  {\bf A}_n '(F_l)Y +  b^{p_l- \frac n2} \D(|F_l|^{p_l} F_l) +
   \e {\mathcal L}_1 F_l=0.
  \ee
 Here the derivative is given by
  $$
   \mbox{$
  {\bf A}_n '(F)Y= - \n \cdot [|F|^n (\frac n F \, (\n \D F)Y + \n \D
  Y)] + \b_l y \cdot \n Y + \a_l Y.
   $}
   $$

 The rest of the analysis depends on assumed good  spectral properties of
 the linearised operator ${\bf A}_n'(F_l)$. We follow the lines of a
 similar analysis performed for the FBP case in
 \cite[\S~2]{PetI}, where
  the operator ${\bf A}_n'(F_l)$ for $n=1$ turns out to possess   a
 (Friedrichs') self-adjoint  extension with compact resolvent and discrete
 spectrum. Such a self-adjoint extension does not exist for the oscillatory
$F(y)$.
 Here we use general theory of non-self-adjoint operators;
  see e.g.,
 \cite{GGK}.
  A proper functional setting of this
 operator is more straightforward for $N=1$
 (and in the radial setting), where, using the behaviour of
 $F(y) \to 0$ as $ y \to 1$, it is possible to check whether the
 resolvent is compact in a suitable weighted $L^2$ space. In
 general, this is a difficult problem; see below.

We assume that such a proper functional setting is available
 for ${\bf A}_n$, so
we deal with operators having solutions with ``minimal"
singularities at the boundary of the support $S_l$, where the
operator is degenerate and singular. 
  Namely,  we assume that
   ${\bf A}_n'(F_l)$ has discrete spectrum and a complete and closed set of
eigenfunctions
   denoted again by $\{\psi_\g\}$. We also assume that the kernel is finite
dimensional
   and we are able to
  determine the spectrum, eigenfunctions $\{\psi_\g^*\}$, and  the kernel of the
 adjoint operator $({\bf A}_n'(F_l))^*$ defined in  a natural
  way using the topology of the dual space $L^2$ and having
  the same point spectrum  (the latter  is true for compact operators in a
suitable space
  \cite[Ch.~4]{KolF}).


Further, we assume that there exists the orthogonal subspace ${\rm
   Span}\{\psi_\g, |\g| > l\}$ of eigenfunctions of
   ${\bf A}_n'(F_l)$, and we look for solutions of (\ref{ep3}) in the form
    $$
    Y=  \phi_l+ w,
     $$
     where $\phi_l$ belongs to the kernel and hence is analogously given by
     (\ref{a4}) and  $w$ belongs to the orthogonal complement of the kernel.
      In doing so,
    we need to transform (\ref{ep3}) into an equivalent integral equation
    with compact operators, but for convenience, we continue our
    computations using the differential version; see additional details
   in \cite[\S~3]{GW2}.

    Thus,
multiplying (\ref{ep3}) by $\psi_\g^*$ with any $|\g|=l$ in
$L^2$
 and, if necessary, integrating by parts in the differential term $y \cdot
\n F_l$ in ${\mathcal L}_1 F_l$, we obtain the following
orthogonality condition
 of solvability 
(Lyapunov-Schmidt's branching equation \cite[\S~27]{VainbergTr}):
  \be
 \label{ep41}
  \mbox{$
 b^{p_l- \frac n2}\langle \D(|F_l|^{p_l-1} F_l), \psi_\g^* \rangle =
 -\e \langle {\mathcal L}_2 F_l, \psi_\g^* \rangle \quad \mbox{for all} \quad |\g|=l.
   $}
    \ee
    These are algebraic equations for the expansion coefficients
    $\{C_\g\}$
    in (\ref{a4}) and the function $b=b(\e)$.
    Similar to (\ref{M5}), one needs to check whether
    the constants are non zero,
 \be
 \label{CC1}
 \mbox{$
\langle \D(|F_l|^{p_l-1} F_l), \psi_\g^* \rangle \not =0 \quad
\mbox{and} \quad \langle {\mathcal L}_2 F_l, \psi_\g^* \rangle
\not = 0, $}
 \ee
which is not a simple problem and can lead to restrictions
for such a behaviour. The analysis is much simpler if the kernel
is 1D, which always happens in the radial geometry where we deal
with ordinary differential operators. Then  (\ref{ep41}) is a
single and easily solved algebraic equation, for which the
``transversality" problem (\ref{CC1}) also occurs.

Under the conditions (\ref{CC1}), the parameter $b(\e)$ in
(\ref{ep1}) for $p \approx p_0$ is given by
  \be
 \label{ep5}
 \mbox{$
 b(\e) \sim [\g_l(p_l-p)]^{\frac{2 \a_l}{1+2 \a_l}}.
 $}
     \ee
The direction of each $p_l$-branch and, whether the bifurcation is
sub- or supercritical, depends on the sign on the coefficient
$\g_l$ that follows from (\ref{ep41}). This can be checked
numerically only, but, in general, we expect that the most of
these nonlinear bifurcations are subcritical so the $p_l$-branches
exist for $p < p_l$.

\smallskip

    For $n=0$, a rigorous justification of this bifurcation
    analysis can be found in \cite[\S~6]{GW2}, where a
    countable number of $p$-branches was shown to originate at bifurcation
    points (\ref{Crp}) and were
 detected on the basis of known spectral properties of the
 corresponding linear operator in (\ref{tf3}); see details in
  \cite{Eg4}.
    For $n>0$, as we have seen, the justification needs spectral properties of
the
    linearised operator ${\bf A}_n'(F_l)$ and the corresponding
    adjoint one $({\bf A}_n'(F_l))^*$, which
is very difficult for non-radial nonlinear eigenfunctions $F_l$
and is an open problem.
 In particular, it would be important to
  know that the bi-orthonormal eigenfunction
  subset $\{\psi_\g\}$ of the operator ${\bf A}_n'(F_l)$
  is complete and  closed in a weighted $L^2$-space  or in some specially defined closed subspace
   (for $n=0$,
   such results are available \cite{Eg4}).   We expect that for $n \approx 0$,
there
  exist critical exponents for the TFE with absorption that are
 close to those in (\ref{Crp}) at $n=0$. This can be checked by
 standard
 branching-type calculus; see  \cite[App.~A]{GHCo},
 where nonlinear eigenfunctions of the rescaled PME in $\ren$ were
 studied by a branching approach.


\section{The Cauchy problem: towards global extensions of $p$-branches}
 \label{SectGlob}


\subsection{Examples of various profiles}

 Recall that, for $p \not =p_0$ in the ODE (\ref{ODEop}), we still have a 2D bundle
 at the singular interface point (\ref{as1}), but now, for even
 profiles, we also need to satisfy two symmetry boundary conditions at
 the origin:
 \be
 \label{2d1}
  f'(0)=f'''(0)=0.
   \ee
 Therefore, unlike the third-order problem (\ref{OF44}),
 (\ref{f12.0}) in the critical case $p=p_0$, we cannot expect continuous sets of
 solutions.
 Actually, as was shown in \cite{Gl4, GBl6, EGW, GW2}, in these
 non-critical cases, there occurs a countable set of $p$-branches
 of similarity profiles, which originate at the standard (for
 $n=0$) or nonlinear bifurcation points $\{p_l\}$ as explained in Section
 \ref{CPpn}. The global behaviour of such
$p$-branches can be complicated and we do not intend to
study these delicate open questions in any detail, restricting
ourselves to examples only.

\begin{figure}
\centering
\subfigure[$n= 1$, $p=3$]{
\includegraphics[scale=0.52]{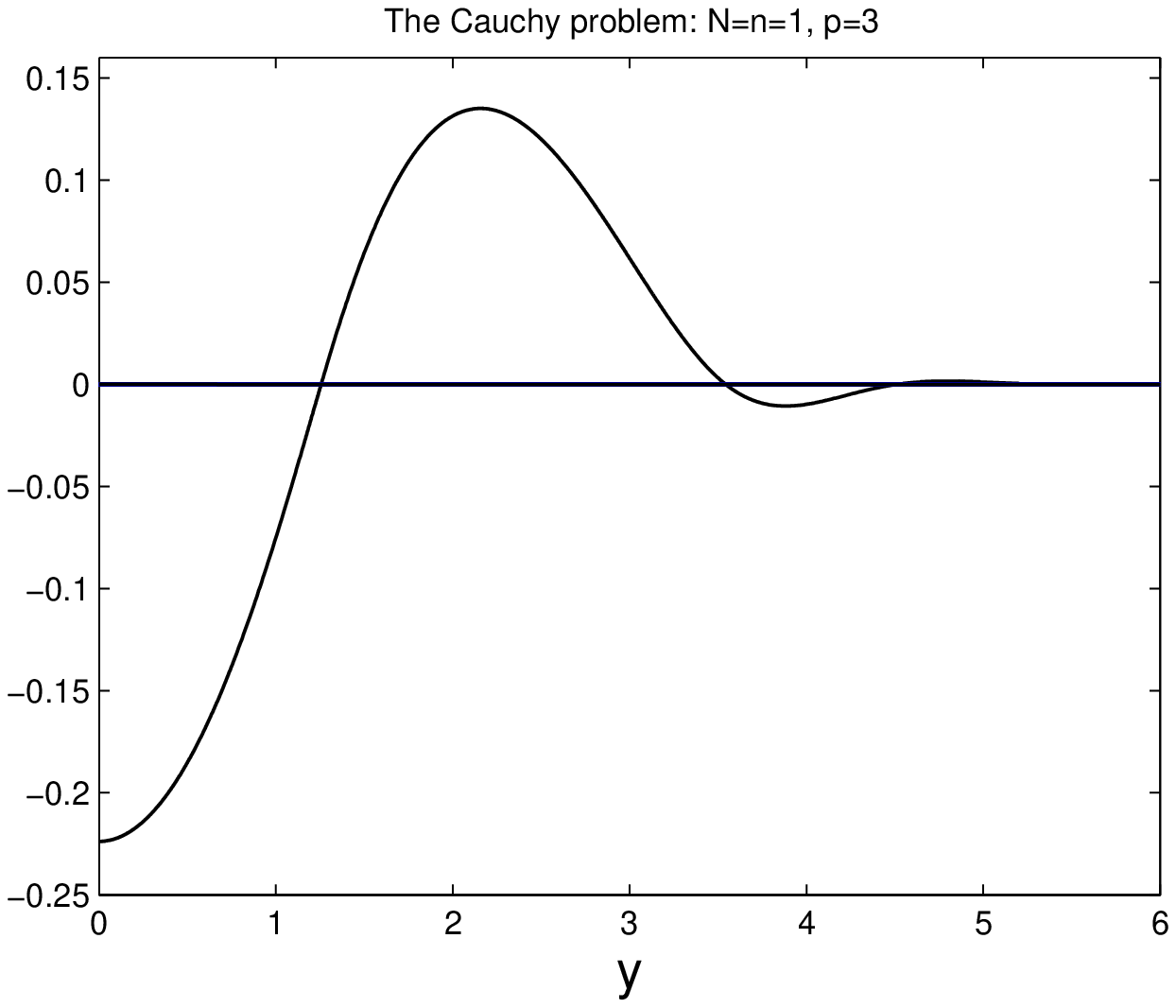} 
}
\subfigure[$n= \frac 12$, $p=2$]{
\includegraphics[scale=0.52]{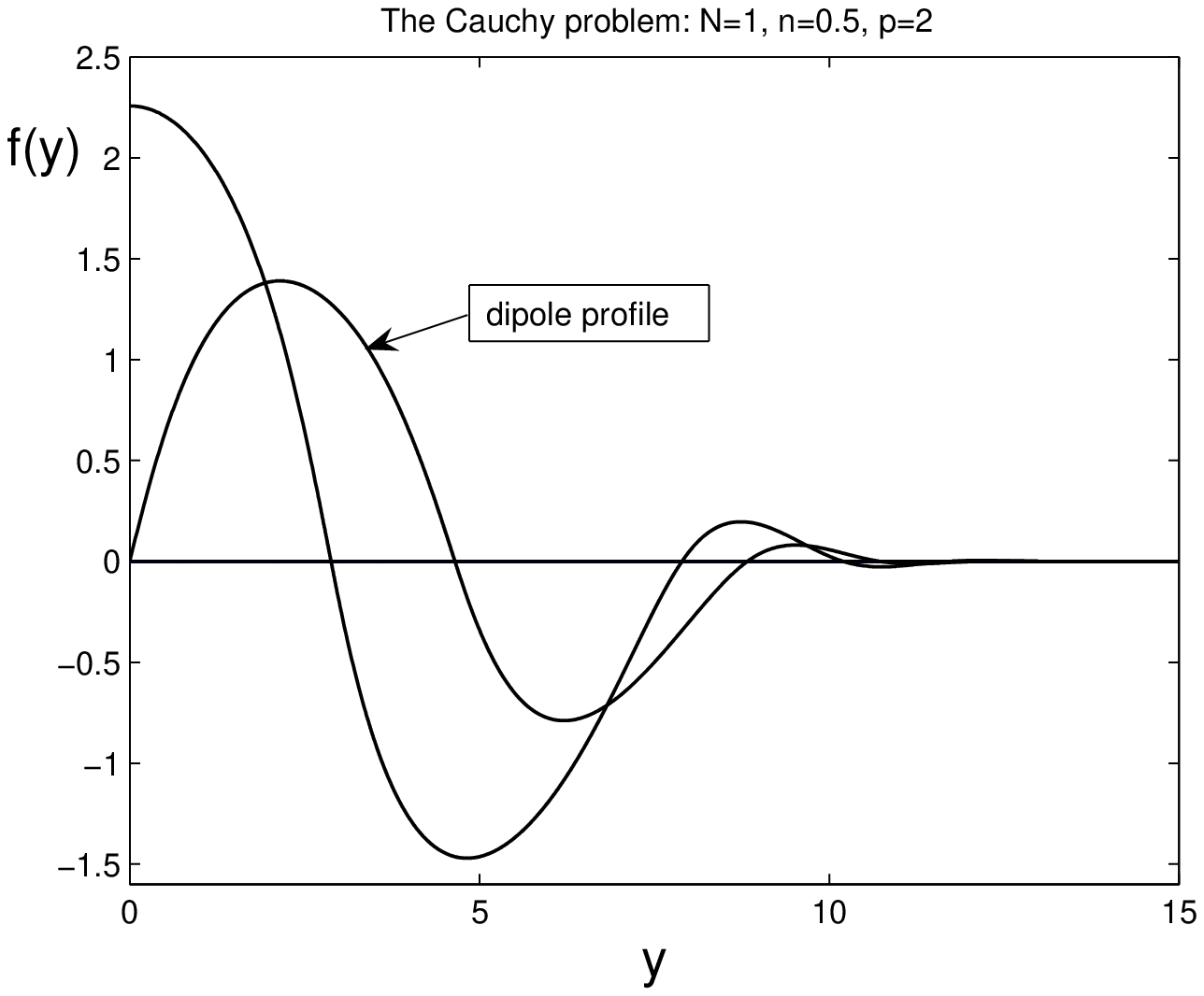} 
}
   \vskip -.3cm
    \caption{\small  Examples of the VSS  profile for the CP
    satisfying (\ref{ODEop}), (\ref{2d1}) for $N=1$;  $n=1$, $p=3$ (a) and
    $n = \frac 12$, $p=2$ (b).}
   \vskip -.3cm
 \label{F13}
\end{figure}

In Figure \ref{F13}, we present some VSS profile $f(y)$ for $N=1$
in two cases: $n=1$ and $p=3< p_0=4$ in (a) and $n= \frac 12$,
$p=2$ in (b).
In (b), we also show the first dipole
profile $f_1(y)$ that,
 instead of (\ref{2d1}), satisfies the
anti-symmetry conditions at the origin,
 \be
 \label{ans1}
 f(0)=f''(0)=0 \quad \Longrightarrow \quad f(-y) \equiv -f(y).
 \ee

 Note an important feature of such compactly supported profiles that is  seen in the figures:
 by (\ref{mma1}),
 {\em their mass must be zero}.
 This necessary condition essentially ``deforms" the VSS similarity
 profiles,
 so that it gets difficult to distinguish in Figure \ref{F13} their
 Sturmian-like properties on the numbers of dominant extrema and transversal
 zeros (if these apply at all). Note that the orthogonality property in (\ref{mma1}) is
 perfectly valid for the eigenfunctions (\ref{ps1}) for $n=0$ (see (\ref{p3})), which
 made it possible to develop the above branching theory.




\subsection{$p$-bifurcation branches: numerics}

\begin{figure}[htp]

\vskip 0cm
 \hspace{-1.5cm}
\includegraphics[width=12cm]{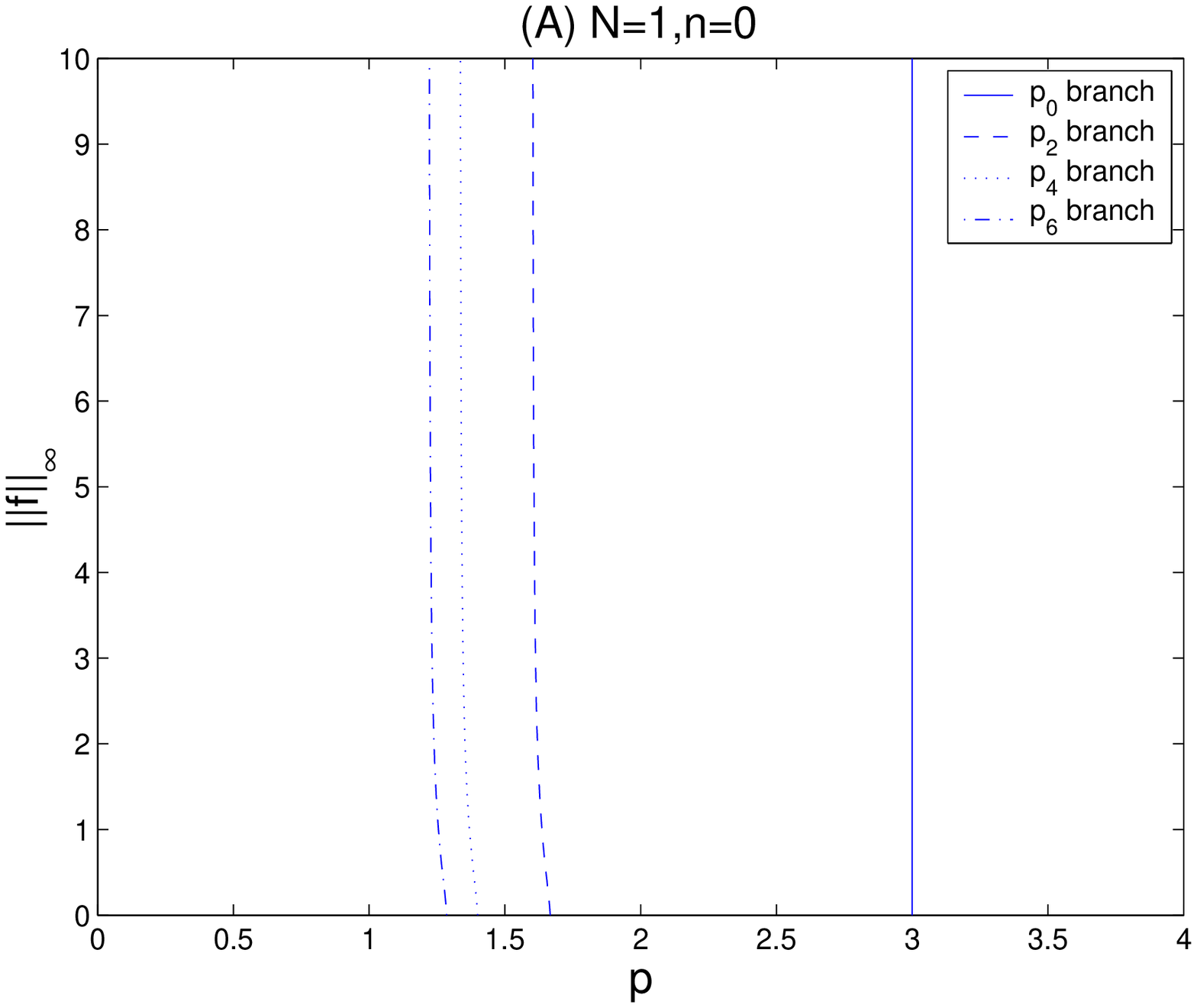}
\hskip 7cm
\includegraphics[width=7.3cm]{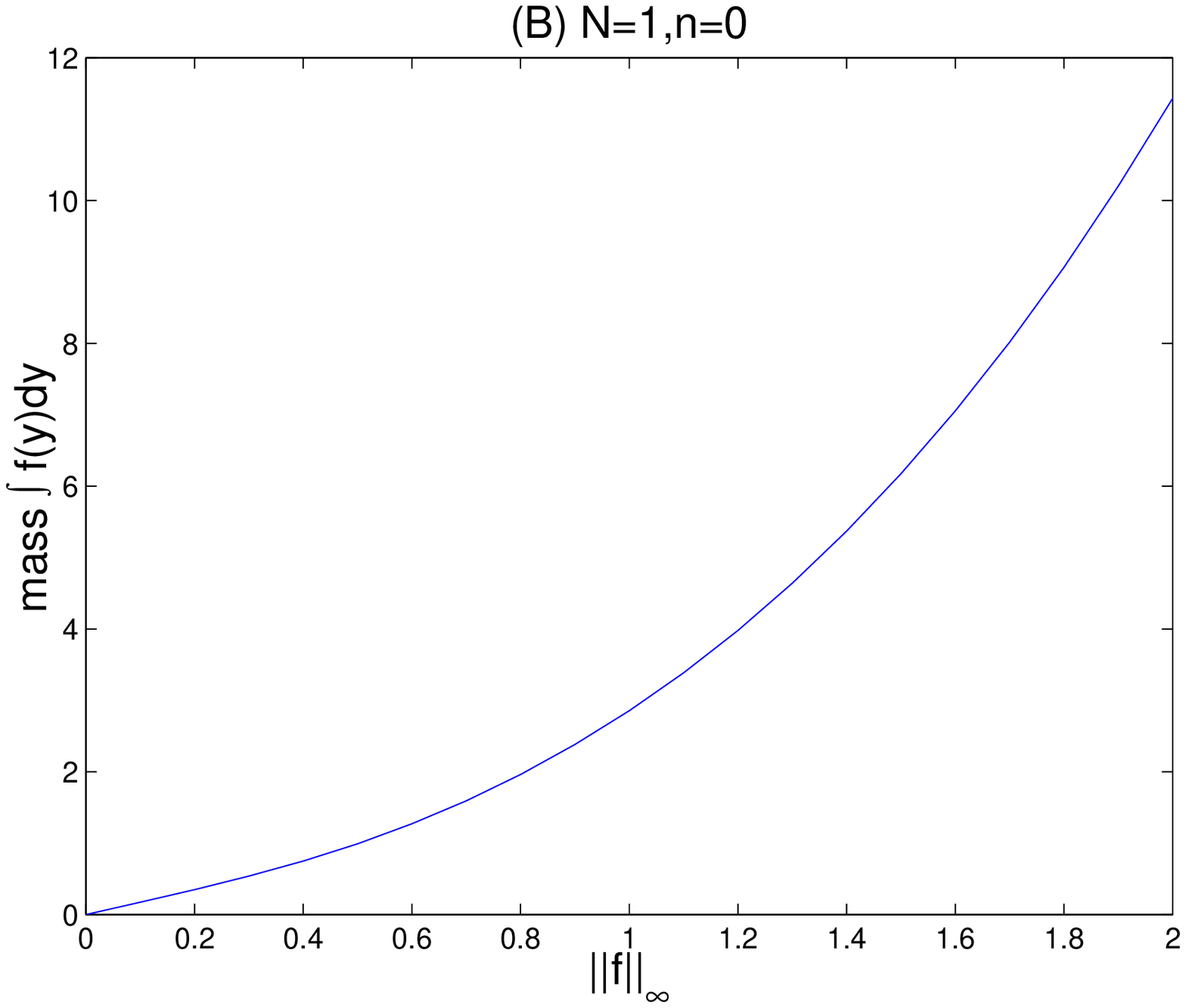}
\hskip 0.5cm
\includegraphics[width=7.3cm]{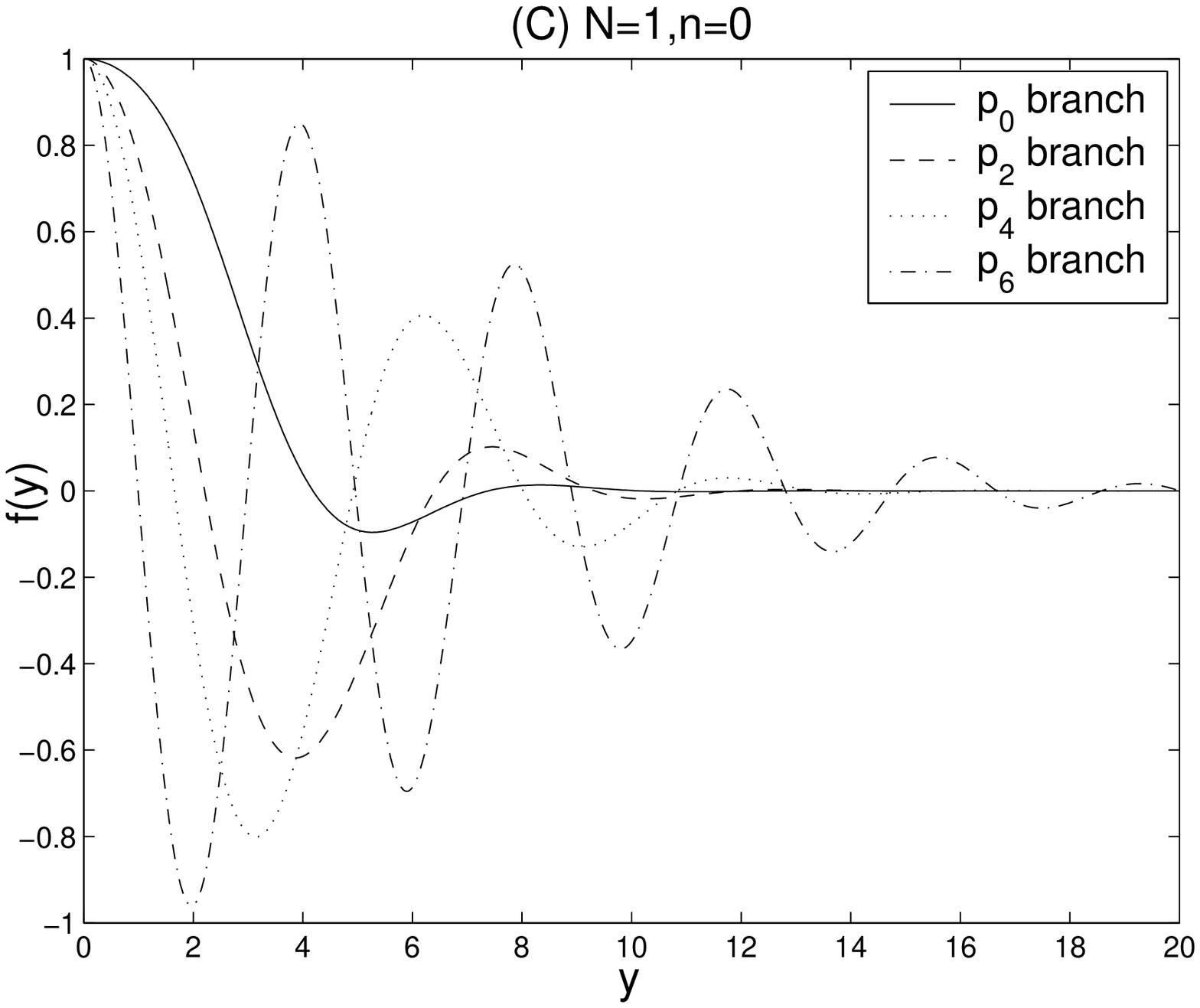}

\vskip 0cm \caption{ \small The bifurcation $p$-diagram and
associated plots for the Cauchy case when $N=1,n=0$. (A) shows the
$p$-bifurcation branches emanating from the critical exponents
$p=p_{l} = 1 + \frac 2{1+l}$ on the $p$-axis. The first four
(even) branches $l=0,2,4,6$ are plotted. (B) illustrates the
monotonicity of the mass of solutions in the critical case
$p=p_0=3$, whilst (C) shows selected profiles on the four branches
in (A) that have $||f||_{\infty}=1$.
 }
 \label{pfig}
\end{figure}

Consider the semilinear case $n=0$, which is known to be simpler,
but correctly describes the
expected general behaviour of the $p$-branches, at least, for
sufficiently small $n>0$ (bearing in mind, that a continuous
``homotopic" deformation as $n \to 0$ is observed in a number of
papers mentioned above). Thus,
 Figure \ref{pfig} illustrates this subcritical case for even symmetry
conditions (\ref{2d1}) at the origin. The branches are seen to
remain distinct, which contrasts markedly with the supercritical
case \cite{Gl4, EGW}, where the branches are increasing with $p$,
so that  the $p_2$ and higher branches ``intersect" the vertical
$p_0$ branch, $\{p=3\}$ at points (profiles) $f$ with the zero
mass as in (\ref{mma1}).

In Figure \ref{pfig}(A), we observe a strong, almost vertical,
growth of these $p$-branches, which bifurcate, respectively, at
 \be
  \label{p111}
  \mbox{$
  p_2= 1+ \frac{2}{1+2}= \frac 53, \quad p_4=1+ \frac 2{1+4}=\frac
  75, \quad p_6=1+ \frac 2{1+6}= \frac 97.
  $}
   \ee
   This is not surprising, since the ODE (\ref{ODEop}) for $N=1$
   and $n=0$ assumes, as $p \to 1^-$, balancing the terms
    \be
    \label{p1112}
    \mbox{$
   ...\, + (|f|^{p-1}f)''+ \,...\, + \frac 1{2(p-1)}\, f=0\,\,
   \Longrightarrow \,\, f=C \, \hat f, \,\,\, \mbox{where}\,\,\, C(p) \sim (p-1)^{-\frac 1{p-1}}
 $}
  \ee
  (the scaled function $\hat f(y)$ is then supposed to be ``almost" uniformly bounded, probably
   up to slower factors).
  Therefore, by (\ref{p1112}), $f(y)$ has a super-exponential growth as $p \to 1^+$.
  Since the bifurcation values in (\ref{p111}) are already
  sufficiently close to 1 and the bifurcations are subcritical,
  these explain  such a strong growth of all the $p$-bifurcation
  branches in Figure \ref{pfig}(A).

In the nonlinear case $n>0$, such an convincing justification of the
general $p$-diagram is not available. Indeed, as we have shown in
the previous section, in the simplest case $l=0$, i.e., $p=p_0$,
the bifurcation of this vertical (in $p$) branch occurs from a
nonlinear eigenfunction, which is the scaled source-type profile
of the nonlinear thin film operator. Other nonlinear
eigenfunctions of the thin film operator are still unknown
possibly excluding the second dipole-like eigenfunction.  We
expect that the discrete nature of the $p$-bifurcation branches
discovered in \cite{EGW} for $n=0$ remains  valid for small $n>0$,
where the nonlinear branching  points cannot be calculated
explicitly as in (\ref{Crp}) and follows from a complicated
nonlinear eigenvalue problem for the thin film operator.

\section{FBP: local behaviour of radial self-similar profiles}
   \label{SectLocR}

\subsection{Interface conditions in the radial setting}

 As in \cite{Gl4}, for the FBP, we need to look for   profiles $f(y)$, which
vanish at  finite $y=y_0>0$ and describe asymptotics of the
general solution satisfying the zero contact angle and zero-flux
conditions: as $y \to y_0^-$,
 \be
 \label{y0pm}
  \mbox{$
 f(y) \to 0, \,\,\,\, f'(y) \to 0,
\,\,\,\, - |f|^n \big(\frac1{y^{N-1}}(y^{N-1}f')'\big)'  +
(|f|^{p-1} f)' \to 0. $}
 \ee
The classification below also applies to
 the corresponding bundles occurring for
the Cahn--Hilliard equation with  $n=0$, \cite[\S~2]{EGW}. We
assume that the conditions (\ref{nn22}) hold (as in the blow-up
case \cite{Bl4}).

\subsection{Two-dimensional asymptotic bundle of similarity profiles}

The derivation of this bundle is standard and
coincides with that in \cite[\S~3.2]{Gl4}. Namely, as above, on
integration, we obtain the ODE (\ref{OD-}).
Therefore, for $n \in [0, \frac 32)$, $p > \frac 32$, the {\em
two-parametric bundle} of solutions is (cf. \cite{BPW,FB0})
 \be
\label{B112} \mbox{$
   f(y) = C_0(y_0-y)^2 - A (y_0 -y)^q (1+o(1)),
   $}
 \ee
where $y_0>0$ and $C_0>0$ are arbitrary parameters and the
correction term depends upon the value of $N$ and $n$, namely:
\[
 \mbox{$
\mbox{(a) for $N=1$}, \quad  q=5-2n, \hspace{0.5cm} A = -
\frac{\beta y_0 C_0^{1-n}}{(5-2n)(4-2n)(3-2n)}
  ;
   $}
\]
\[
\mbox{(b) for $N \ge 2$}, \quad  \left \{ \begin{array}{ll} q=3,
\hspace{0.5cm} A = \frac{C_0(N-1)}{3 y_0}  \hspace{1cm} & \mbox{if
$n<1$},\ssk \\
    q=3,    \hspace{0.5cm} A = \frac{C_0(N-1)}{3 y_0} - \frac{\beta y_0}{6}  \hspace{1cm} & \mbox{if $n=1$} ,\\
    q=5-2n, \hspace{0.5cm} A = - \frac{\beta y_0 C_0^{(1-n)}}{(5-2n)(4-2n)(3-2n)}  \hspace{1cm} & \mbox{if $n>1$} .
   \end{array}
\right.
\]

  This expansion exists also for $n=0$ and has nothing to do with the CP
exhibiting infinite propagation. It can be also used in the
FBP posed for the Cahn--Hilliard equation with zero contact angle
and zero-flux conditions. Proving such expansions demands a rather
involved application of Banach's contraction principle;  see also
\cite{BPW, FB0}.

\section{FBP: source-type similarity 
patterns in the critical case $p=p_0$}
 \label{Sect5}

In the critical case $p= p_0$, the ODE (\ref{ODEop}) can be
integrated once reducing the radial ODE to a third-order equation
of the form
 \be
\label{RescBU}
 \mbox{$ |f|^n \bigl(f''+ \frac{N-1}y \, f'\bigr)' - (|f|^{p_0-1}f)' -\b y f=0,
\quad \mbox{where}
  \,\,\, \b=   \frac{1}{4+nN}. $}
 \ee
     At the interface, we
take  conditions in (\ref{y0pm}) written as 
\be
   f(y_0)=f'(y_0)=0 ,
\label{BUfz0}
 \ee
and 
we complete the problem statement by taking the  symmetry
condition at
the origin (\ref{f12.0}).  
 The ODE (\ref{RescBU}) itself  then  implies the second symmetry
 condition
$f'''(0)=0$.

The mass $M$ of $f$ is a parameter, which is useful in
distinguishing the solutions of (\ref{RescBU}) subject to
(\ref{BUfz0}) and (\ref{f12.0}). Here, we set 
\be
 \mbox{$
        M = \int\limits_{0}^{y_0}  y^{N-1}   f(y)\, {\mathrm d}y.
         $}
\label{BUpar1}
 \ee
As another parameter, we can take  $f''(y_0)$ or
$C_0=-\frac{1}{2}f''(y_0)$, corresponding to the bundle (\ref{B112}).
There are also other choices of parameters, such as $\{f(0),f''(0)\}$.

For the FBP,
 the local parameters
$\{C_0,y_0\}$ determine  a two-parameter shooting problem in
order to attain the single symmetry condition at the origin
(\ref{f12.0}).
The global mass-parameter $M$ is useful in classifying our
solutions as bifurcations from critical mass values associated
with non self-similar steady states.

Thus the statement of  the FBP comprises  the ODE (\ref{RescBU})
with the symmetry condition (\ref{f12.0}) within  the bundle
(\ref{B112}) with two parameters. Therefore, we expect
  a countable set of continuous families of solutions to exist, which can
  be  parameterised relative to $y_0$ or with respect to the
mass $M$ of the profiles. The asymptotic structure of the first
(stable) branch of such similarity profiles is easier to detect.


\subsection{Continuous mass-branches of solutions of the FBP}

For simplicity, we again consider the 1D case. The global
structure of similarity profiles does not essentially depend on
$N$. In the critical case (\ref{p011}), i.e.,
$p_0=n+3$ for $N=1$, the ODE (\ref{RescBU})
 takes the form (\ref{OF44}),
where at the interface point $f(y)$ is assumed to belong to the
bundle (\ref{B112}).
 Note that this equation is obviously non-variational.
 The existence of a continuous set of
solutions  with any sufficiently small mass is proved by a
shooting argument exactly as in \cite[\S~5]{EGW}. The only
difference is that all the positive large solutions as $y \to
+\infty$ are concentrated in a three-dimensional bundle around the
profile
 $$
f_*(y) =  y^{\frac 2{n+2}}
 \varphi_*(\ln y),
 $$
 where the oscillatory component $\varphi_*(s)$ is a periodic
 function of changing sign of a certain autonomous ODE that is
 easy to derive.



\subsection{Numerical construction of the global similarity patterns}

\begin{figure}[htp]

\vskip -2.2cm
 \hspace{-1.5cm}
\includegraphics[width=7.3cm]{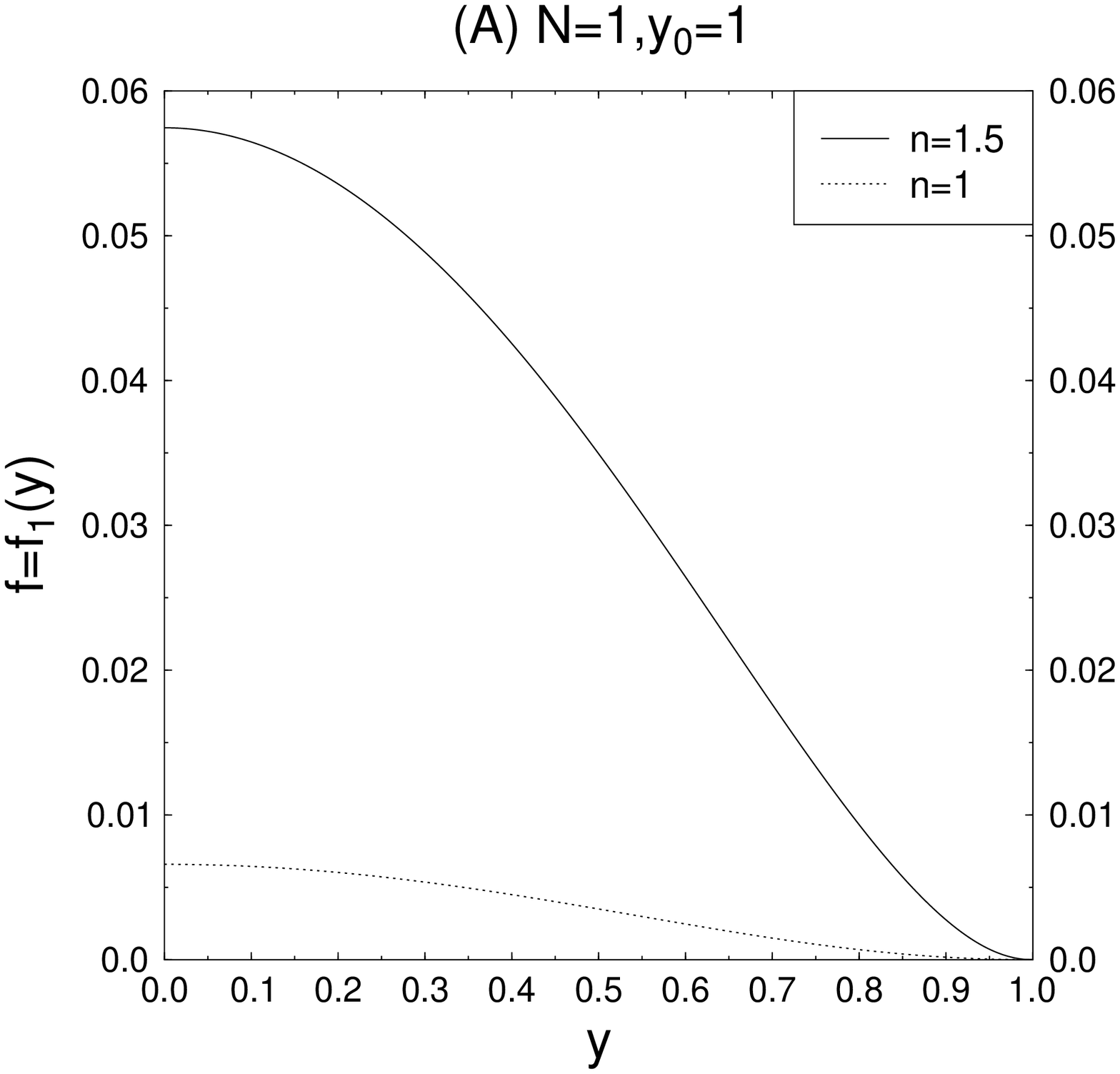}
\includegraphics[width=7.3cm]{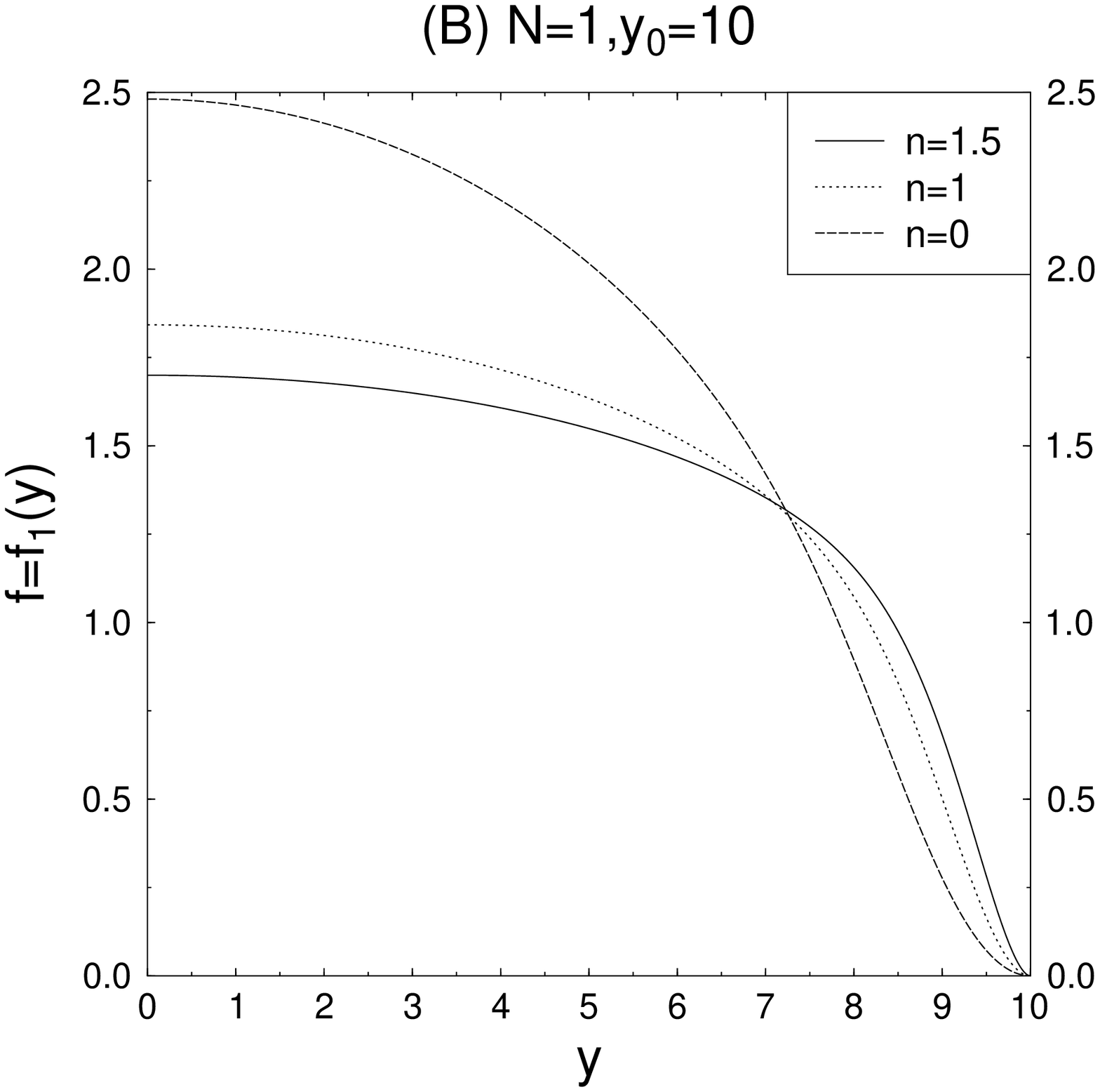}
\vskip -3.2cm \hspace{-1.5cm}
\includegraphics[width=7.3cm]{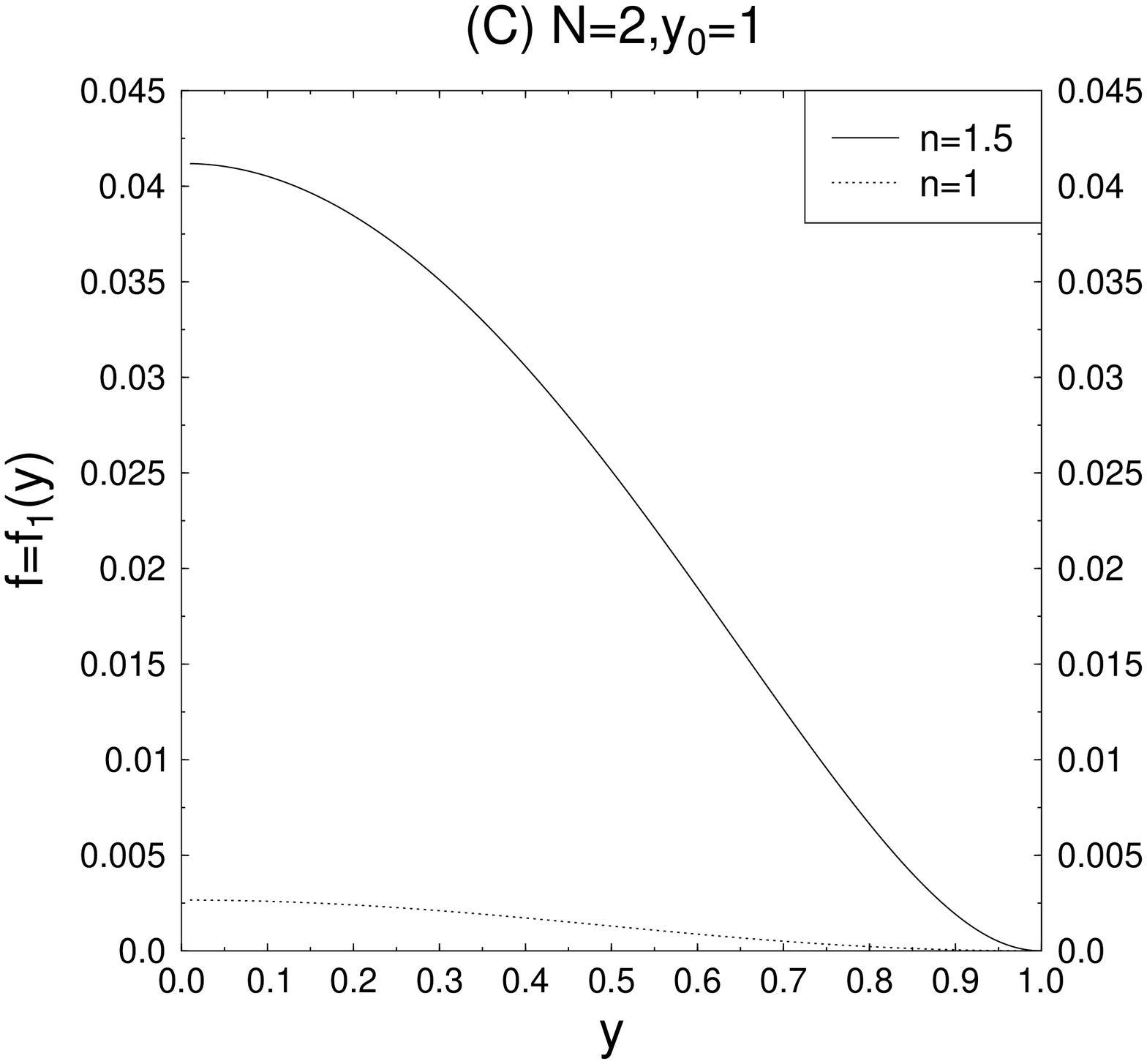}
\includegraphics[width=7.3cm]{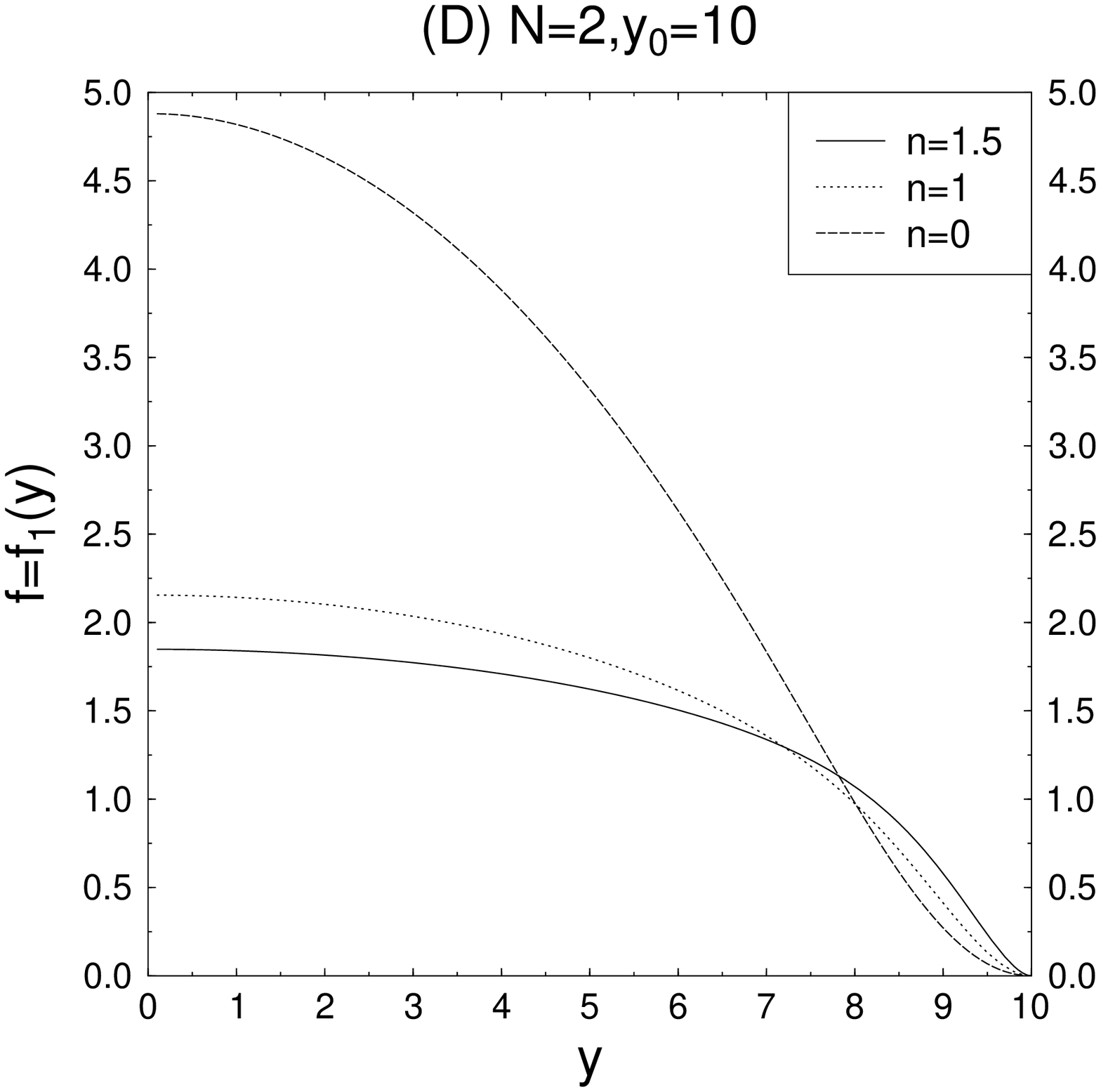}

\vskip -1.2cm \caption{ \small Illustrative numerical solutions of
the base $f=f_1(y)$ global similarity profiles in the cases
$N=1,2$. Parameter values of $y_0=1$ and $10$ were used for the
support length together with the selected values $n=0,1,1.5$ for
the index $n$. (A) and (B) show the one-dimensional case $N=1$,
(C) and (D) the two-dimensional case $N=2$. In both dimensions,
profiles when $y_0=1$ were not satisfactorily obtainable for
smaller $n$ values less than $1$ (due to their values being below
scheme tolerances). }
 \label{num1}
\end{figure}

\begin{figure}[htp]

\vskip -2.2cm
 \hspace{-1.5cm}
\includegraphics[width=7.3cm]{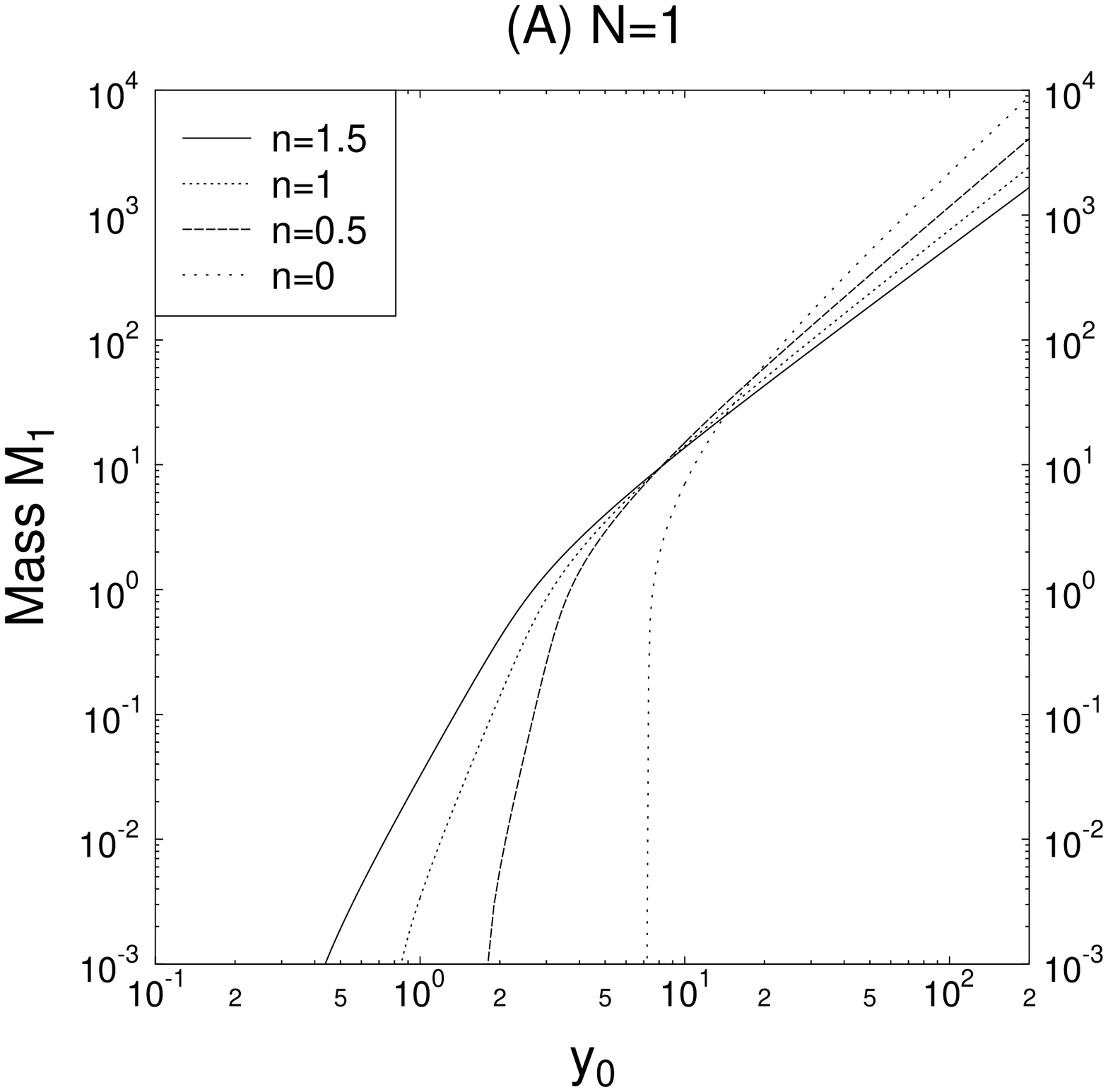}
\includegraphics[width=7.3cm]{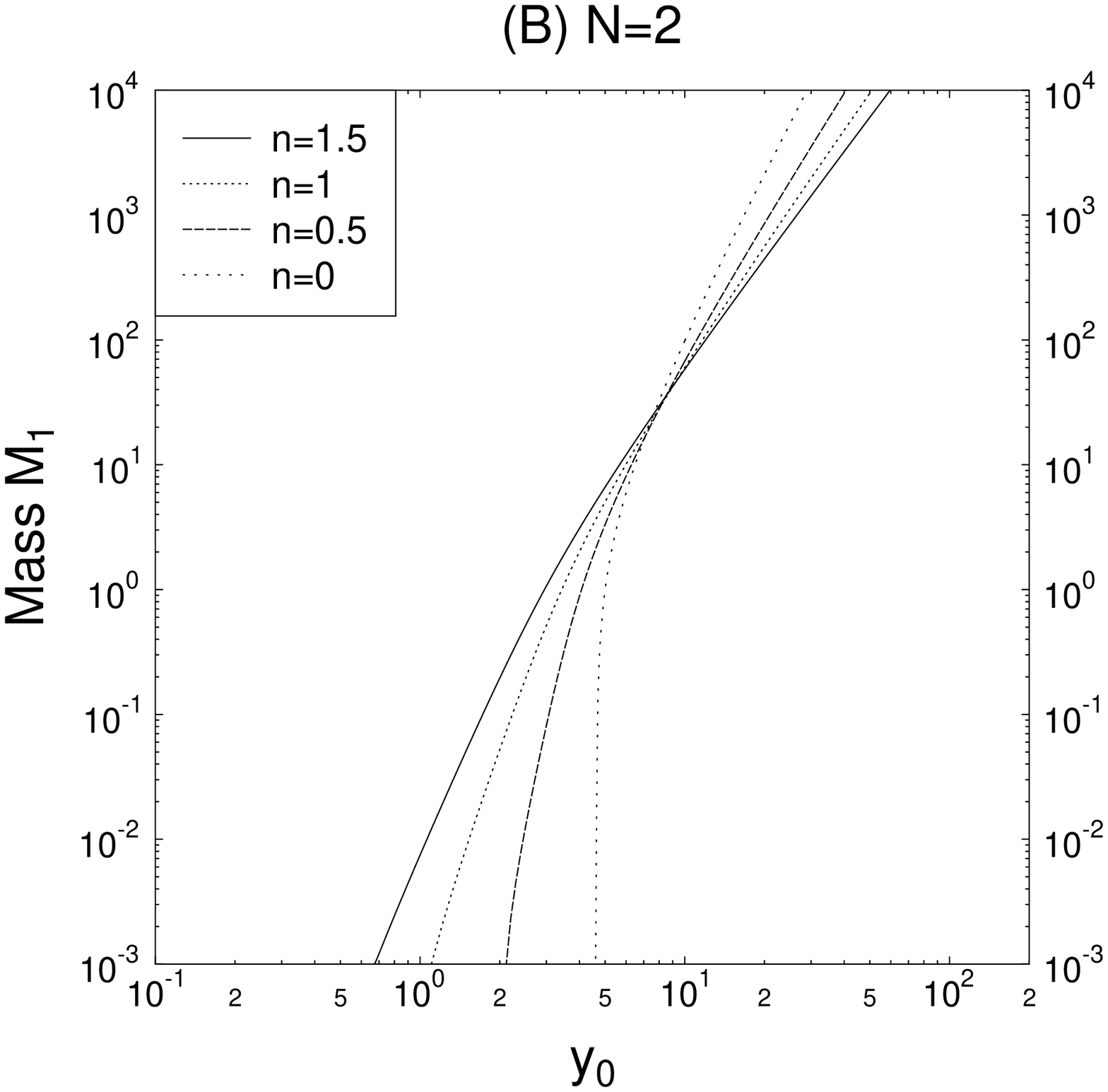}

\vskip -1.2cm \caption{ \small Parameter plots of the mass $M=M_1$
of the base profiles solutions $f=f_1(y)$ with $y_0$ for selected
$n$ values. (A) gives the one-dimensional case $N=1$ and (B) the
two dimensional case $N=2$. }
  \label{num2}
 \end{figure}

\begin{figure}[htp]

\vskip -2.2cm
 \hspace{-1.5cm}
\includegraphics[width=7.3cm]{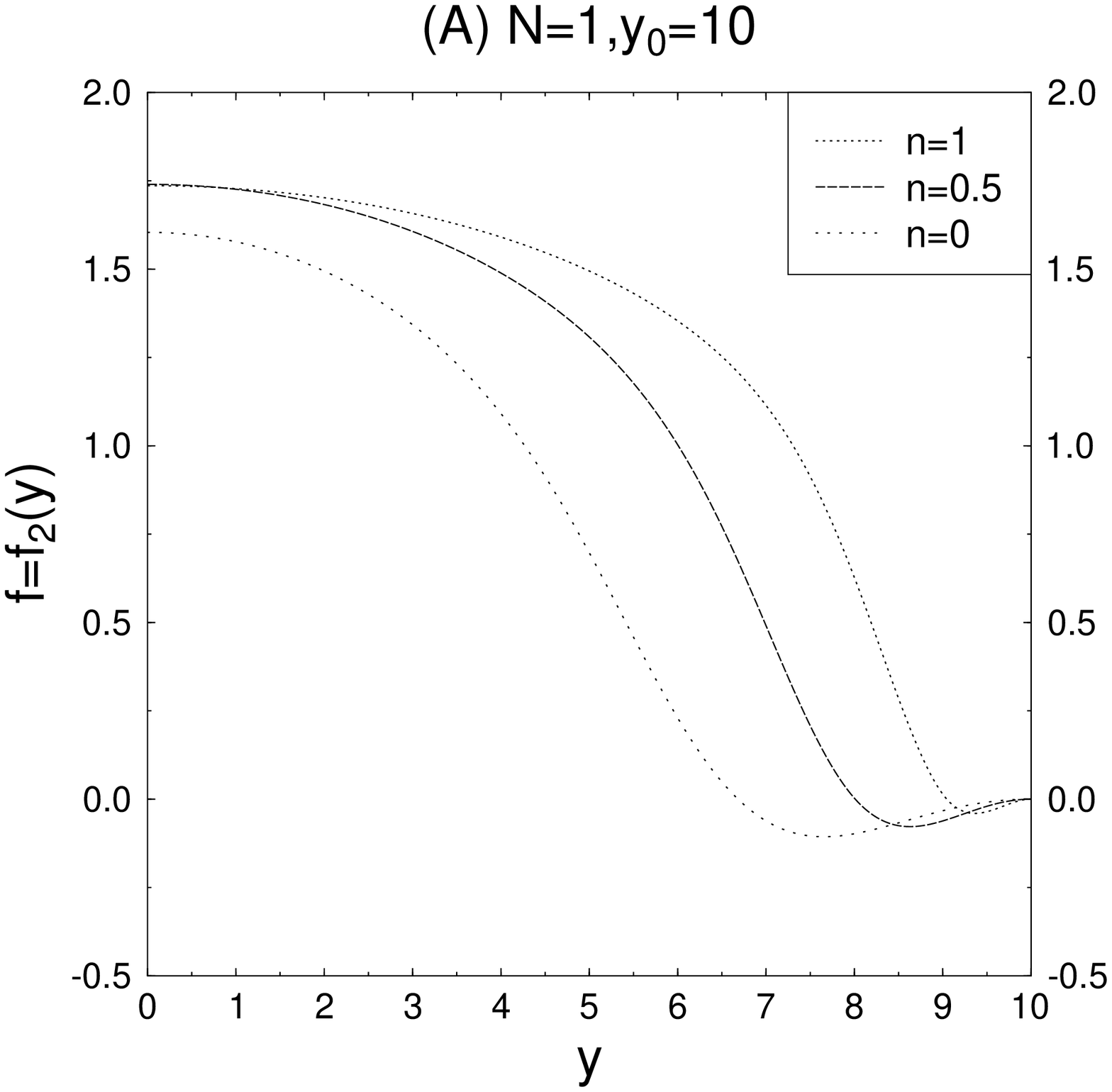}
\includegraphics[width=7.3cm]{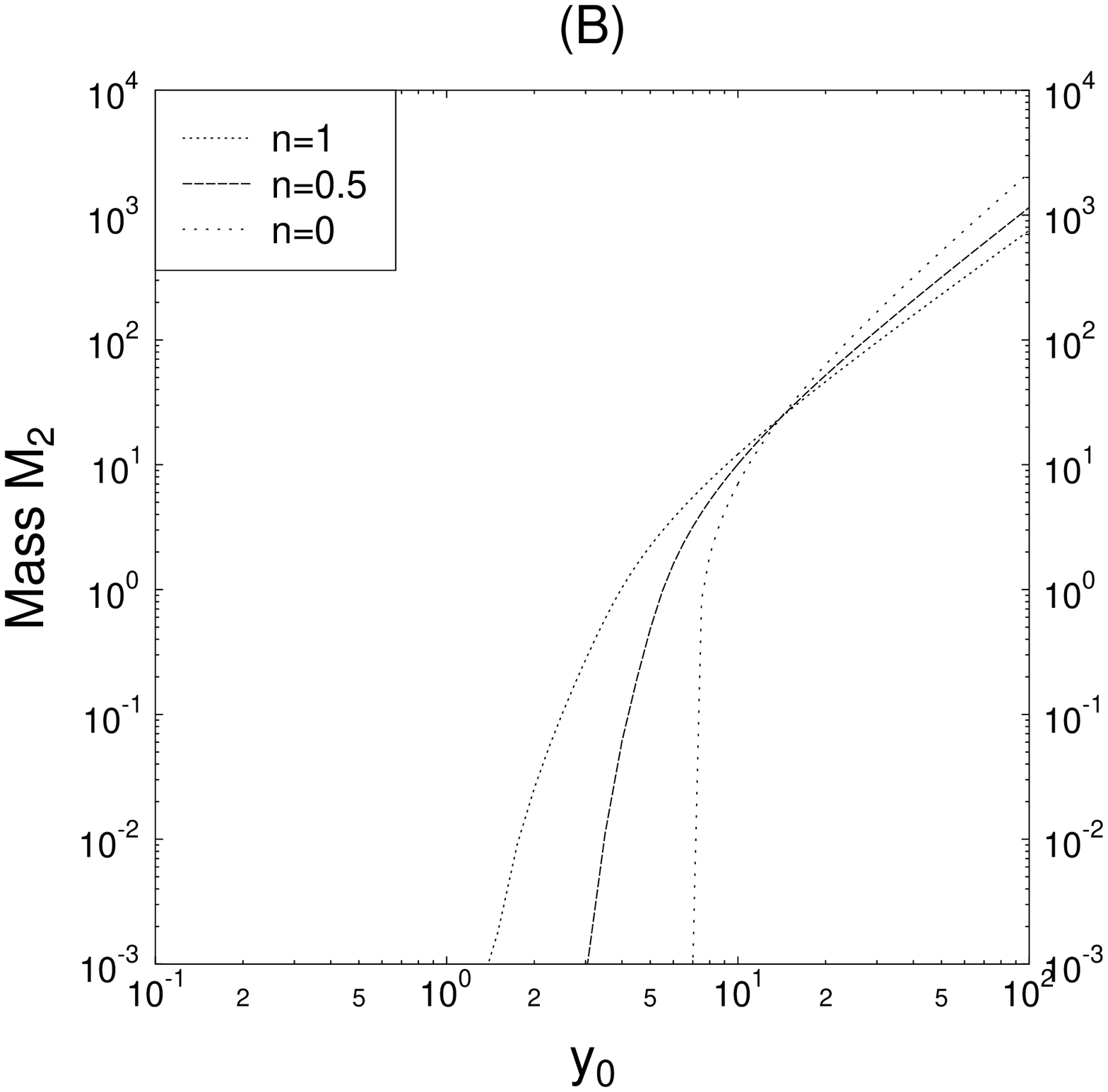}
\vskip -1.2cm \caption{ \small Parameter plots of the second
profile solutions $f=f_2(y)$ for selected $n$ in one-dimension
$N=1$. (A) gives the profiles for $y_0=10$, whilst (B) gives the
mass as the support $y_0$ varies. }
  \label{num3}
 \end{figure}

The boundary-value problem is now (\ref{RescBU}) written as
\be
\label{RescBUb}
 \mbox{$ f''' + \frac{N-1}{y} \bigl(  f'' - \frac{1}{y}\, f' \bigr)
   - p_0 |f|^{\frac{2}{N}} f' -\b y f |f|^{-n}=0 , $}
 \ee
together with (\ref{BUfz0}), (\ref{f12.0}), and (\ref{BUpar1}).
For fixed $N$ and $n$, the numerical results presented below
suggest that there is a countable number of solutions for a given
$y_0$. We denote the profiles with positive mass as $f_k$ where
the index $k=1,2,\,...$ represents the number of sign changes
($k-1$) of the profile over the interval $[0,y_0]$ (it also being
related to the number of maxima and minima).  The results for the
base profiles $f_1(y)$ are presented in Figure \ref{num1}
(similarity profiles) and Figure \ref{num2} ($y_0$-bifurcation
diagram), whilst Figure \ref{num3} gives the corresponding results
for the next profiles $f_2(y)$ in the set. We remark that there
are the corresponding reflected profiles $-f_k(y)$ which have negative
masses.

\subsection{Asymptotic expansions}

\subsubsection{The small mass limit  $M \to 0$}

As in \cite[\S~4.2]{Gl4}, we present an approach for the asymptotic
expansion in the limit of small mass.
 This limit is easier for $n=0$
\cite[Prop.~6.2]{EGW}, and, after rescaling, $f$ converges to
the fundamental similarity rescaled profile $F$ defined in
(\ref{tf3}) as the first (with the eigenvalue $\l_0=0$) normalized
eigenfunction of the linear non-self-adjoint operator ${\bf B}$
there. For $n>0$, we identify the corresponding zero-mass limit as
follows. Let us
introduce the small parameter 
\be
       \e = M^{\frac{2}{N}}
\label{eq:epsilon}
 \ee
   and perform the  scaling
 \be
\label{sc++}
 f(y) = \e^{\frac{2N}{4+nN}} g(z),\quad y=
\e^{\frac{nN}{2(4+nN)}} z,
 \quad y_0= \e^{\frac{nN}{2(4+nN)}} z_0 .
\ee
 Then the problem (\ref{RescBU})--(\ref{BUpar1}) becomes
\be
\label{eq:gbc} 
  \left\{
  \begin{matrix}
 |g|^n \bigl(g''+ \frac{N-1}z \, g'\bigr)' - \b z g=  \e (|g|^{p_0-1}g)',
\qquad \qquad\quad \,\,\,\ssk\ssk \\
   g(z_0)=g'(z_0)=0 , \,\,\, g'(0)=0, \,\,\,
         1 = \int_{0}^{z_0}  z^{N-1}   g(z)\, {\mathrm d}z ,
 \end{matrix}
  \right.
 \ee
  where $'$ now denotes $\frac{\mathrm d}{{\mathrm
d}z}$. We pose the regular expansions
\[
    g = g_0 - \e g_1 + ..., \quad z_0 = z_0^0 - \e z_0^1 + ...\quad
\mbox{as\,\, $\e \to 0$},
\]
to obtain the leading-order problem
\be
\label{eq:g0}
 \left\{
 \begin{matrix}
  |g_0|^n\bigl(g_0''+ \frac{N-1}z \, g_0'\bigr)' - \b z g_0 = 0, \,\,\, z>0,
  \qquad\qquad\qquad\,\,\,
 \ssk\ssk \\
   g_0(z_0^0)=g_0'(z_0^0)=0 , \,\, g_0'(0)=0, \,\,
         1 = \int_{0}^{z_0^0}  z^{N-1}   g_0(z)\, {\mathrm d}z .
  \end{matrix}
  \right.
\ee
 The first-order problem, with $ z_{0}^1 = -
 \frac{g_1'(z_0^0)}{g_0''(z_0^0)}$, is
\be
\label{eq:g1} \left\{
 \begin{matrix}
  |g_0|^n\bigl(g_1''+ \frac{N-1}z \, g_1'\bigr)' + (n-1) \b z g_1  = - p_0 g_0'
|g_0|^{p_0-1}, \quad \quad
 \ssk\ssk \\
   g_1(z_0^0)=g_1'(z_0^0)=0 , \,\, g_1'(0)=0, \,\,
         0 = \int_{0}^{z_0^0}  z^{N-1}   g_1(z)\, {\mathrm d}z .
         \end{matrix}
          \right.
 \ee

\ssk

These leading and first-order problems coincide with those in the
unstable case described in \cite[\S~4.2.1]{Gl4}. As such the
details will not be repeated, other than to correct typographical
errors in (4.18) therein,  which should read
\[
 \mbox{$
  z_{01}^1 = - \frac{p_0}{4 z_{01}^0 [8(N+2)(N+4)]^{\frac{2}{N}}} \,
  p_N(z_{01}^0), \,\,\, \mbox{where} \,\,
   $}
  \mbox{$
p_1(z_{01}^0) = -\frac{1024}{45045} \bigl(z_{01}^0\bigr)^{12},
\,\,\, p_2\bigl(z_{01}^0\bigr) = -\frac{1}{60}(z_{01}^0)^{8}.
 $}
 \]

We remark that the zero mass limit requires the support to vanish
(as shown through the scaling (\ref{sc++})). The corresponding
limit of vanishing support but with the  mass remaining non-zero does
not possess non-trivial solutions in contrast to the unstable
case. This is consistent with the physical interpretation of the
terms on the RHS in (\ref{GPP}), since the effect of the second
order term (the gravity) is no longer opposed by the fourth-order term
(the surface tension) as in the unstable case (\ref{GPPun}).

\subsection{The large mass limit $M\to \infty$}

This limit occurs simultaneously with increase in support length
$y_0$. We thus introduce the scalings
 \be
 y=y_0 z, \quad f(y) = y_0^{\frac{2}{p_0-1}} g(z),\quad  M = y_0^{N+\frac{2}{p_0-1}}
 m,
\label{lmscale} \ee so that (\ref{RescBU})--(\ref{BUpar1}) becomes
\be
\label{eq:lmprob} 
  \left\{
  \begin{matrix}
  y_0^{-2( 1+ \frac{2}{(p_0-1)N})}\bigl(g''+ \frac{N-1}z \, g'\bigr)' =  \b z g |g|^{-n} + p_0 |g|^{\frac{2}{N}}g',
 \ssk\ssk \\
    g(1)=g'(1)=0 , \,\,\, g'(0)=0, \,\,\,
         m = \int_{0}^{1}  z^{N-1}   g(z)\, {\mathrm d}z ,
 \end{matrix}
  \right.
 \ee
  where $'$ again denotes $\frac{\mathrm d}{{\mathrm d}z}$. This gives a
  singular perturbation problem in the limit $y_0 \to \infty$, comprising an outer
  region $0\leq z <1$ together with an inner region near $z=1$.

\ssk

{\bf (I) Outer problem.} We pose the regular expansions
\[
    g = g_0 + o(1), \quad m = m_0 + o(1) \quad
\mbox{as\,\, $y_0 \to \infty$},
\]
to obtain the leading-order outer problem
\be
\label{eq:outerg0}
 \left\{
 \begin{matrix}
   \b z g_0 |g_0|^{-n} + p_0 |g_0|^{\frac{2}{N}} g_0'  = 0, \,\,\,
  \qquad\qquad\qquad\,\,\,\qquad
 \ssk\ssk \\
   g_0(1)=0 =0 , \,\, g_0'(0)=0, \,\,
         m_0 = \int_{0}^{1}  z^{N-1}   g_0(z)\, {\mathrm d}z .
  \end{matrix}
  \right.
\ee We thus obtain the explicit solution
\be
 \mbox{$
     g_0 = \big[ \frac{\beta (p_0-1)}{2p_0} (1-z^2) \big]^{\frac{1}{p_0-1}} ,
   \quad m_0 = \frac{1}{2} \big[
   \frac{\beta(p_0-1)}{2p_0}\big]^{\frac{1}{p_0-1}}\,\,
        \frac{\Gamma(\frac{p_0}{p_0-1} ) \Gamma(\frac{N}{2})}{\Gamma(\frac{p_0}{p_0-1}+\frac{N}{2})}.
         $}
\label{eq:g0sol} \ee This solution for $g_0$ does not satisfy the
condition $g_0'(1)=0$ and thus we require an inner region near
$z=1$. We note that this leading order outer solution is common to
all profiles $f_k$ in this large mass limit.

\ssk

\begin{figure}[htp]

\vskip -2.2cm
 \hspace{-1.5cm}
\includegraphics[width=7.3cm]{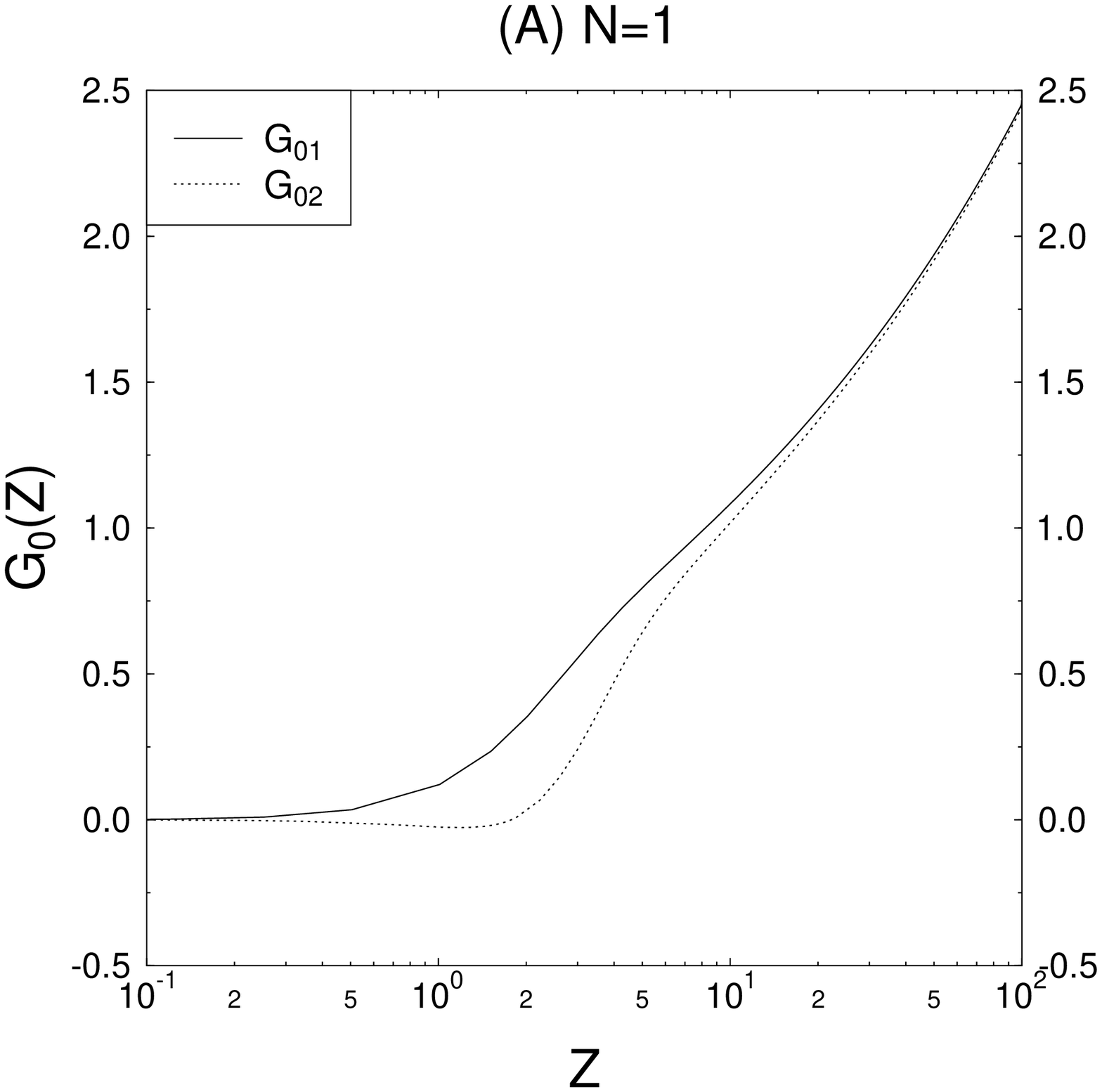}
\includegraphics[width=7.3cm]{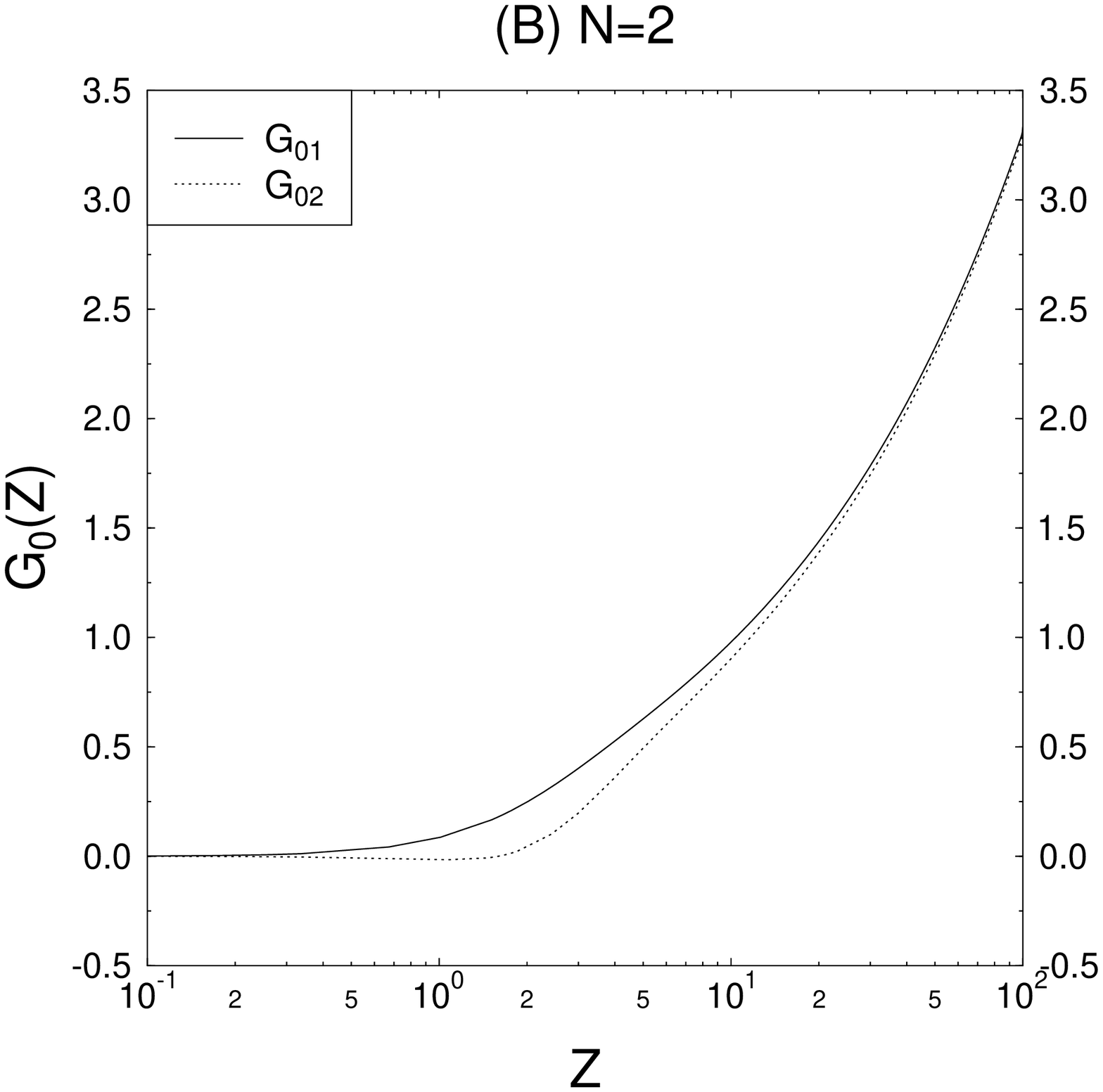}

\vskip -1.2cm \caption{ \small Numerical solution for the leading
order inner problem in the large mass limit. Illustration of the
first two profiles in a sequence of increasing changing sign
profiles. The one-dimensional case $N=1,n=1$ is shown in (A),
whilst (B) shows the two-dimensional case $N=2,n=1$. The curvature
values at the origin may be used to distinguish the profiles and
for comparison, we record the values
$G_{01}''(0)=0.2917,G_{02}''(0)=-0.0960$ for (A), whilst we have
$G_{01}''(0)=0.2130,G_{02}''(0)=-0.0667$ in (B).
 }
 \label{mlarge}
\end{figure}

{\bf (I) Inner problem.} Near $z=1$ we introduce the scalings
\be
    z = 1 - \delta Z , \quad g = \delta^{\frac{1}{p_0-1}} G,
\label{eq:lminscales} \ee where dominant balance in
(\ref{eq:lmprob}) gives the small parameter
$\delta=\delta(y_0)=y_0^{-\frac{nN+4}{nN+3}}
$, and for $Z=O(1)$ we pose
\[
    g = G_0 + o(1) \quad \mbox{as\,\, $y_0 \to \infty$},
\]
to obtain the leading-order inner problem
\be
\label{eq:innerG0}
 \left\{
 \begin{matrix}
    G_0'''  =  -\b G_0 |G_0|^{-n} + p_0 |G_0|^{\frac{2}{N}} G_0',  \,\, \,\,\,
  \qquad\qquad\qquad\,\,\,\qquad\quad
 \ssk \\
   G_0(0)= G_0'(0)= 0 , \quad  G_0 \sim \big[ \frac{(p_0-1) \beta Z}{p_0} \big]^{\frac{1}{p_0-1}} \;\;\; \mbox{as $Z \to +\infty$} ,
 \end{matrix}
  \right.
\ee
   where $'$ is now $\frac{\mathrm d}{{\mathrm d}Z}$. The last
condition arises from the  matching with the outer solution
(\ref{eq:g0sol}).  Consistent with the set $f_{k}$, we anticipate
a countable set of solutions denoted by $G_{0k}$, $k=1,2,...\,$,
and  distinguished by the number of sign changes ($G_{0k}$ having
$k-1$ sign changes). Numerical solutions for $G_0$ are shown in
Figure \ref{mlarge}, where the first two profiles are shown for
the parameter value $n=1$ in the $N=1$ and $N=2$ cases.

\section{FBP: on countable sets of $p$-branches of similarity patterns
for $p \not = p_0$}
 \label{Sectp0}


Similar to Sections \ref{CPTFEn} and \ref{CPpn} for the Cauchy
problem, we now intend to develop an analogous analytical approach
to show existence of $p$-branches of similarity profiles in the FBP
setting. This will demand rather unusual a special ``spectral theory" for the corresponding
linear problem for $n=0$.

\subsection{The origin of countable $p$-branches}

 In view of essential differences and additional difficulties arising
 for the FBP, we restrict ourselves to the simplest case $N=1$, so that the FBP setting
includes standard three free-boundary conditions
 \be
 \label{fbc1}
 f(y_0)=f'(y_0)=(|f|^n f''')(y_0)=0 \quad \mbox{at an unknown
 boundary} \,\,\, y=y_0>0,
 \ee
 accomplished with two symmetry (\ref{2d1}) or anti-symmetry ones
 (\ref{ans1}) at the origin.
 Overall, for the 1D fourth-order ODE
 \be
  \label{.1}
   \mbox{$
  {\bf A}_+(f) \equiv - (|f|^n f''')' + \b\, f'y + \a f +
  (|f|^{p-1}f)''=0, \quad y \in (0,y_0),
   $}
    \ee
    where $ \a= \frac
  1{2p-(n+2)}$ and $ \b= \frac{p-(n+1)}{2[2p-(n+2)]}$,
these give five conditions plus an extra free parameter $y_0$ (a
nonlinear eigenvalue). This looks like a correctly posed problem,
which, in the standard analytic setting, could not have more than
a countable set of solutions, or a finite number of uniformly
bounded ones.

As usual, we next need to consider the corresponding 1D pure TFE
(\ref{eq1}), for which, in the FBP setting there appears the
following {\em nonlinear eigenvalue problem} (cf. (\ref{nl1})):
 \be
 \label{nl1FBP}
  \left\{
   \begin{matrix}
  -(|f|^nf''')'+ \frac{1- \a n}4\, y f' + \a f=0\,\,\,
  \mbox{on}\,\,\, (y_-,y_+),
  \,\,\, \ssk\\
   f=f'=|f|^nf'''=0 \,\, \mbox{at} \,\, y=y_\pm.
   \qquad\qquad\quad\,\,\,\,
    \end{matrix}
    \right.
  \ee
 Note that, unlike the CP one (\ref{nl1}), the space of eigenvalues
 is three parametric, and, besides the usual nonlinear eigenvalue $\a \in \re$ includes
 two free boundary positions $y_\pm$.
  In the two basic simpler cases with the symmetry (\ref{2d1}) or antisymmetry  (\ref{ans1}),
 the eigenvalue space is 2D:
  \be
  \label{.2}
  \mbox{Eigenvalues:} \quad \mu=(\a,y_0)^T \in \re^2,
 \ee
 to which we concentrate upon in what follows.

\subsection{$n=0$: first aspects of linear ``Hermitian spectral theory"}
 \label{S9.2}

 For $n=0$, the nonlinear eigenvalue problem (\ref{nl1FBP})
 becomes a linear one for the operator $\BB$ in (\ref{tf3}):
  \be
  \label{.3}
   \mbox{$
    \BB \psi \equiv -\psi^{(4)} + \frac 14\, \psi'y + \frac 14\,
    \psi= \l \psi \,\,\, (\l=\frac 14 - \a), \quad
    \psi=\psi'=\psi'''=0 \,\,\, \mbox{at} \,\,\, y=y_\pm.
    $}
    \ee
   Even in the simpler case  (\ref{.2}), this is not a standard spectral problem,
     and, moreover, it is not clear whether this can be
    attributed to such classes. Let us comment that each
    ``eigenfunction" $\psi_k(y)$ is supposed to be defined on its
    own interval $(y_{k-},y_{k+})$, but using
    eigenfunction subsets together with notions of completeness,
    closure, etc. may not make sense.

   We now discuss some particular aspects of the problem
   (\ref{.3}) and restrict our analysis to the {\em first} eigenvalue
   $\l_0=0$ and the corresponding {\em even} eigenfunctions.
 On integration, the ODE becomes of the third order,
 \be
  \label{.4}
   \mbox{$
     -\psi''' + \frac 14\,y \psi=0 \,\,\,\mbox{on}
     \,\,\,(0,y_0),\,\,\,
     \psi(0)=1, \,\,\, \psi'(0)=0; \quad \psi(y_0)=\psi'(y_0)=0,
      $}
      \ee
 where for convenience we fix the normalization $\psi(0)=1$ and take symmetry conditions at the origin (the subscript + then being dropped on the domain end point). We can prove the following:

  \begin{proposition}
   \label{Pr.01}
   The $y_0$-eigenvalue problem $(\ref{.4})$, corresponding to
   $\l_0=0$,
    admits a countable set of eigenfunctions
   $\{\psi_0^{(k)}, \, k \ge 1\}$ defined on the intervals
    $\{(0,y_0^{(k)})\}$, and,
     \be
     \label{.5}
  \mbox{ as $k \to \iy$}, \quad   y_0^{(k)}  \sim k^{\frac 34} \to +\iy \quad \mbox{and}
     \quad \psi_0^{(k)}(y) \to F(y),
      \ee
 where $F(y) \equiv \psi_0^{(\iy)}(y)$ is the first eigenfunction $(\ref{ps1})$
 of the rescaled operator $\BB$ in $(\ref{tf3})$ defined for the CP
 $($i.e., in the whole $\re)$.
  \end{proposition}

 \noi{\em Proof.} A proper solvability of the problem (\ref{.4})
 follows from the known WKBJ-type asymptotics of solutions of this ODE for
 large $y$:
  \be
  \label{.6}
 \mbox{$
  \psi(y) \sim {\mathrm e}^{a y^{4/3}} \quad \Longrightarrow
  \quad a^3= \frac 14\big(\frac 34\big)^3.
   $}
    \ee
 This gives three roots: the  real positive one $a_0= 3 \cdot 4^{-\frac 43}$
 and two complex:
  $$
 \mbox{$
  a_\pm= \big(-\frac 12 \pm {\rm i}\, \frac {\sqrt
 3}2\big)a_0.
 $}
 $$
   A full WKBJ expansion includes also a slow
 growing algebraic multiplying factor $y^{-1/3}$ (not important for the final estimates):
 as $y \to \iy$,
  \be
  \label{.7}
   \mbox{$
    \psi(y) \sim y^{-\frac 13} \Big\{ C_1{\mathrm e}^{a_0 y^{4/3}}
    + C_2 {\mathrm e}^{- \frac {a_0}2\, y^{4/3}}
    \cos \big( \frac{a_0 \sqrt 3}2 y^{\frac 43} +C_3\big)
    \Big\} +
    ...\, ,
     $}
     \ee
     where $C_{1,2,3}$ are real constants.
      Solving the problem at the free boundary $y_0 \gg 1$, we see
      that acceptable roots are concentrated about the roots of
      the $\cos$ or $\sin$ functions, so that $y_0^{(k)}$ satisfy
       \be
       \label{.8}
        \mbox{$
\cos \big( \frac{a_0 \sqrt 3}2 (y_0^{(k)})^{\frac 43}+C_3\big)
\approx 0
 \quad \Longrightarrow
 \quad
\frac{a_0 \sqrt 3}2 (y_0^{(k)})^{\frac 43} \sim \pi\big(\frac 12
+k \big)-C_3,
 $}
  \ee
  whence the estimate in (\ref{.5}). Since $y_0^{(k)} \to \iy$,
  the convergence to the rescaled kernel $F$ in the CP is obvious
  (see Figures below for a further justification). $\qed$

\ssk

\begin{figure}
\centering
\includegraphics[scale=0.65]{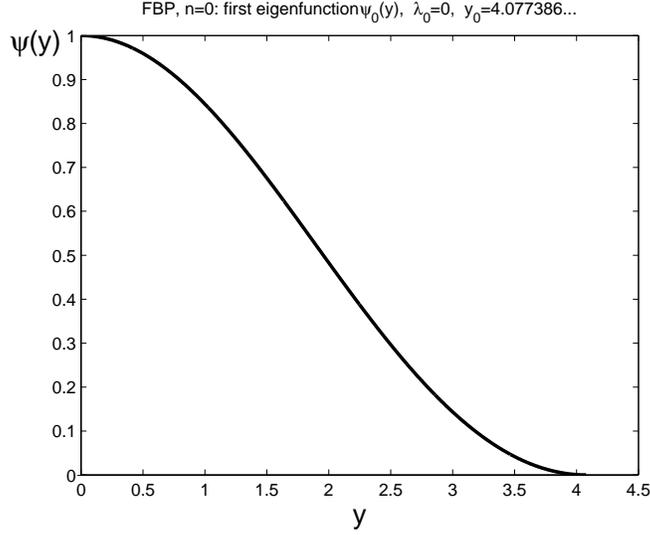}
\vskip -.5cm
 \caption{\small The first positive eigenfunction of (\ref{.4}), with the interface at
  $y_0^{(1)}=4.077386...\,$.}
   \vskip -.3cm
 \label{F01}
\end{figure}

\begin{figure}
\centering \subfigure[four eigenfunctions]{
\includegraphics[scale=0.52]{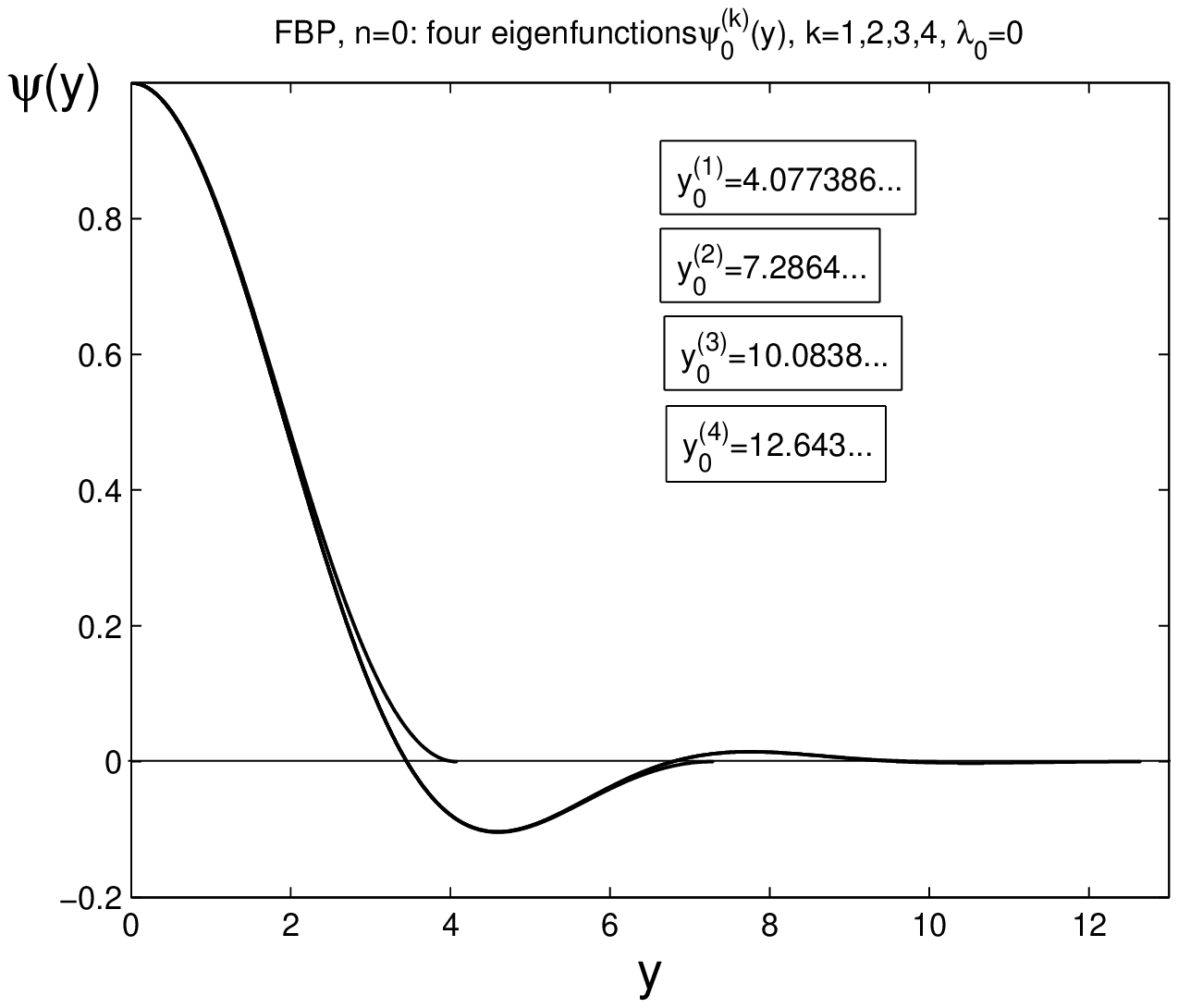}
} \subfigure[enlarged zeros]{
\includegraphics[scale=0.52]{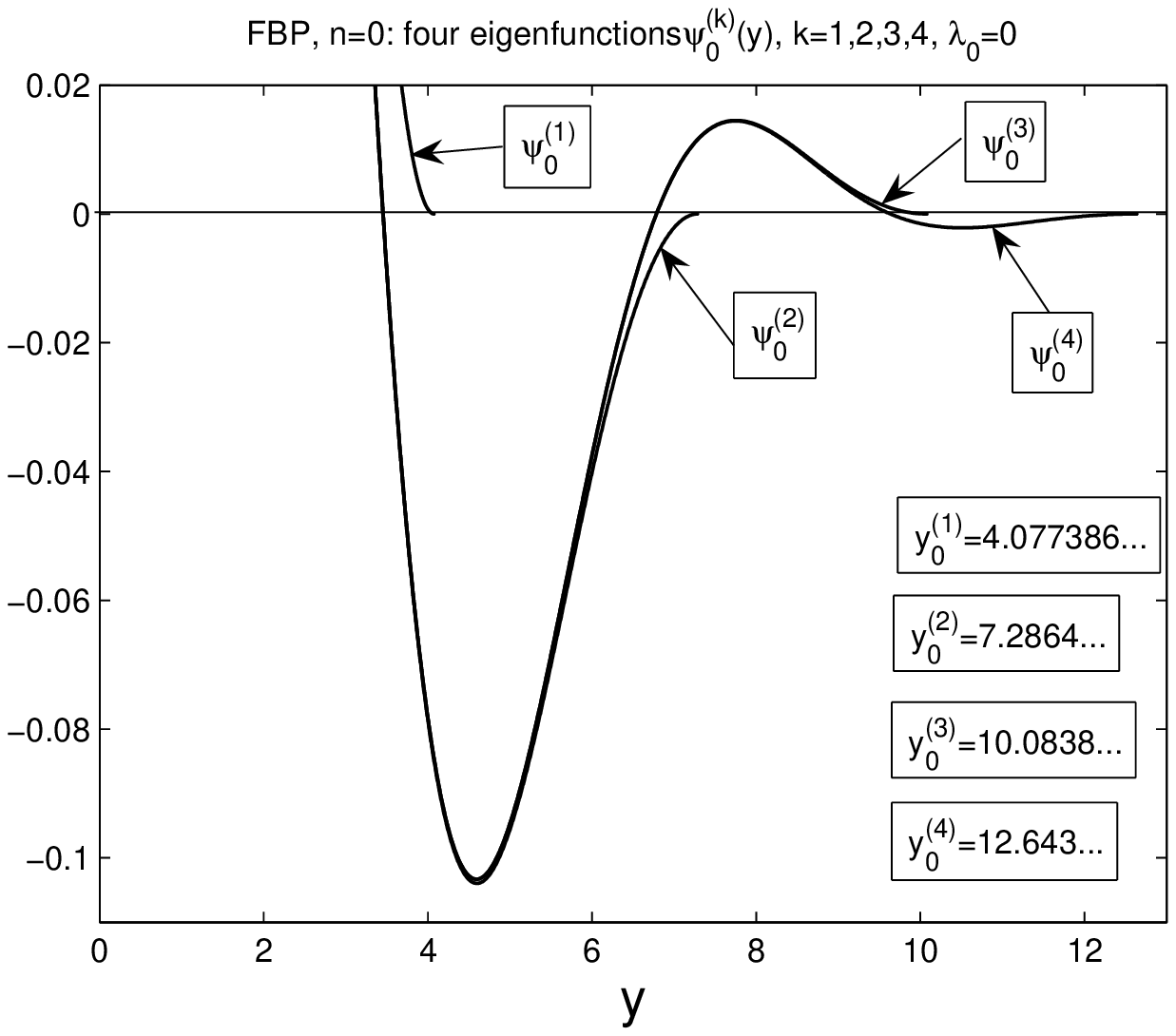}
}
 \vskip -.4cm
\caption{\rm\small First four eigenfunction of (\ref{.4}).}
 \label{F02}
\end{figure}

In Figure \ref{F01}, we show the first positive eigenfunction
$\psi_0^{(1)}$ of (\ref{.4}). Figure \ref{F02} shows first four
eigenfunctions of (\ref{.4}). This illustrates the convergence to
$F$ in (\ref{.5}), where $\psi_0^{(4)}(y)$ is already very close
to $F$, so that the next $\psi_0^{(5)}$ is difficult to detect
numerically.

\begin{figure}
\centering \subfigure[four eigenfunctions]{
\includegraphics[scale=0.52]{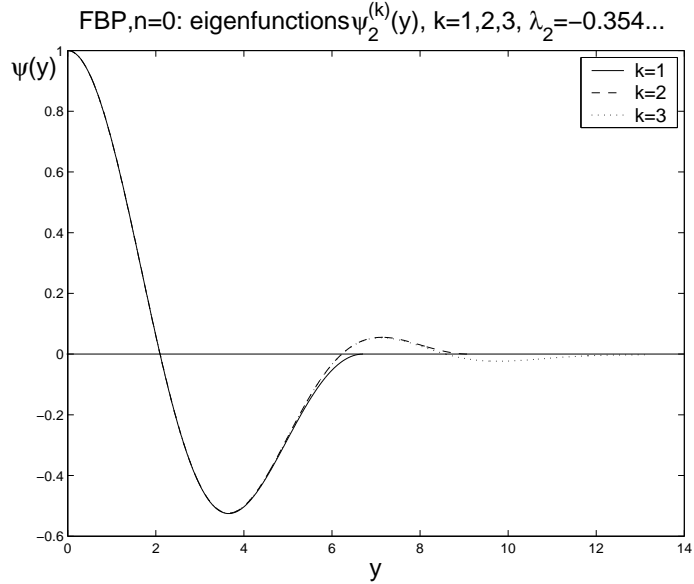}
} \subfigure[enlarged zeros]{
\includegraphics[scale=0.52]{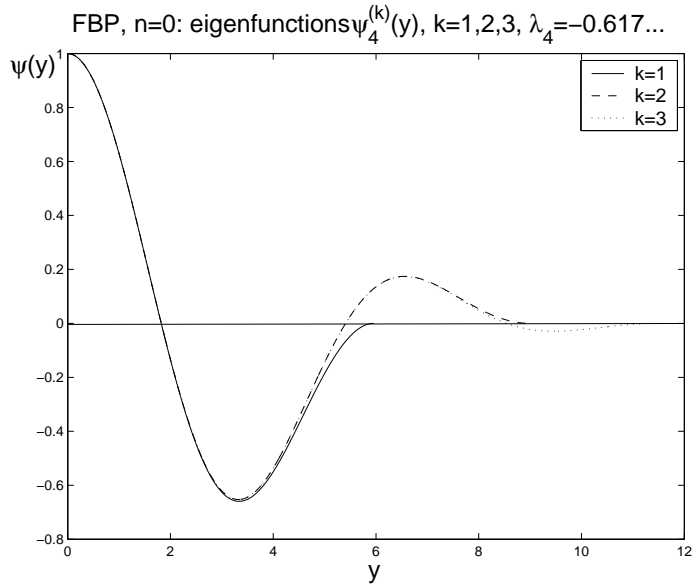}
}
 \vskip -.4cm
\caption{\rm\small First three profiles for the second and third even eigenfunctions $m=2,4$ of (\ref{lnotzero}).}
 \label{F2F4n0}
\end{figure}

For $\lambda \neq 0$, the corresponding statement to (\ref{.4}) is
\be
  \label{lnotzero}
    \left\{
     \begin{matrix}
     -\psi''' + \frac 14\,y \psi = \lambda \int\limits_{0}^{y} \psi (s)\, {\mathrm d}s
      \,\,\,\mbox{on}
     \,\,\,(0,y_+), \qquad\qquad\qquad\,\, \ssk\\
    \,\,\, \psi(0)=1, \,\,\, \psi'(0)=0; \quad \psi(y_+)=\psi'(y_+)=\psi'''(y_+)=0.
     \end{matrix}
      \right.
      \ee
We distinguish the sets of even eigenfunctions $\{ \psi=\psi_m, \lambda=\lambda_m, y_+=y_m \}$, using the subscript $m=0,2,4,\ldots$. Within each set with $m$ fixed, we have a countable set of eigenfunctions $\{ \psi_m=\psi_m^{(k)},y_m=y_m^{(k)} \}, k=1,2,3,...\,$.
The first three eigenfunctions for the second $m=2$ and the fourth $m=4$ sets are shown in Figure \ref{F2F4n0}.


\subsection{Nonlinear ``Hermitian spectral theory"}

Consider now the nonlinear eigenvalue problem (\ref{nl1FBP}) for
$\l_0=0$ in the setting as in (\ref{.4}),
 \be
 \label{.10}
  \mbox{$
 \a= \frac 1{4+n}
 \quad \Longrightarrow \quad
 -|f|^n f'''+ \frac 1{4+n} \, y f =0, \quad f=f'=0 \,\,\, \mbox{at
 the interfaces}.
 $}
 \ee
 Obviously, the first eigenfunction $\psi_0^{(1)}$ is positive and
  is the classic one first obtained in \cite{BPW} (see also
  \cite{FB0} for $N>1$). However, we claim that,  similarly to
  Proposition \ref{Pr.01}, the following holds:

  \ssk

  \noi{\bf Conjecture \ref{Sectp0}.2.} {\em Besides the positive
  solution \cite{BPW, FB0},  in the
  oscillatory range $n < n_{\rm h}$ given in  $(\ref{n**1})$, the problem
   $(\ref{.10})$ admits a
  countable set of sign changing patterns $f^{(k)}=\psi_0^{(k)}(y)$, $k=1,2,...\,$, where
  each $k$-th one has precisely $k$ zeros for $y \in (0,y_0^{(k)}]$.}

  \ssk

   A rigorous proof is
  still incomplete since it requires detailed knowledge of
  oscillatory structures near interfaces such as (\ref{LC11}).
  This should play a similar role to the linear expansion as in (\ref{.7}). Such a deep understanding
  of the nonlinear expansion is still not achieved.

\begin{figure}
\centering
\includegraphics[scale=0.65]{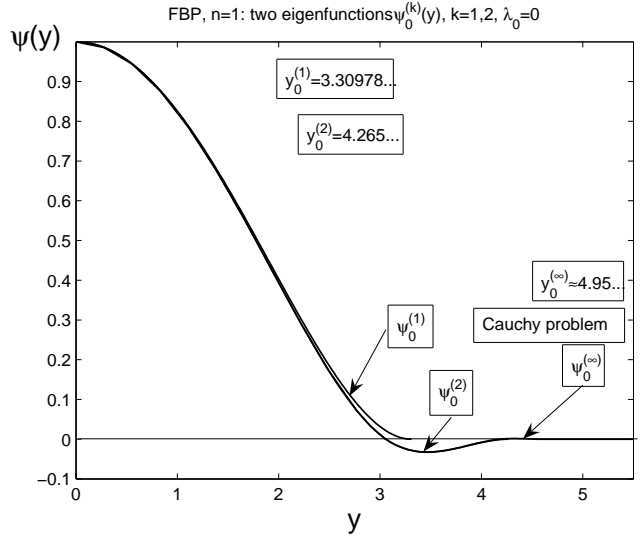}
\vskip -.5cm
 \caption{\small Two nonlinear eigenfunctions of (\ref{.10}) for $n=1$.}
   \vskip -.3cm
 \label{F03}
\end{figure}

\begin{figure}
\centering \subfigure[$\sim 10^{-3}$]{
\includegraphics[scale=0.52]{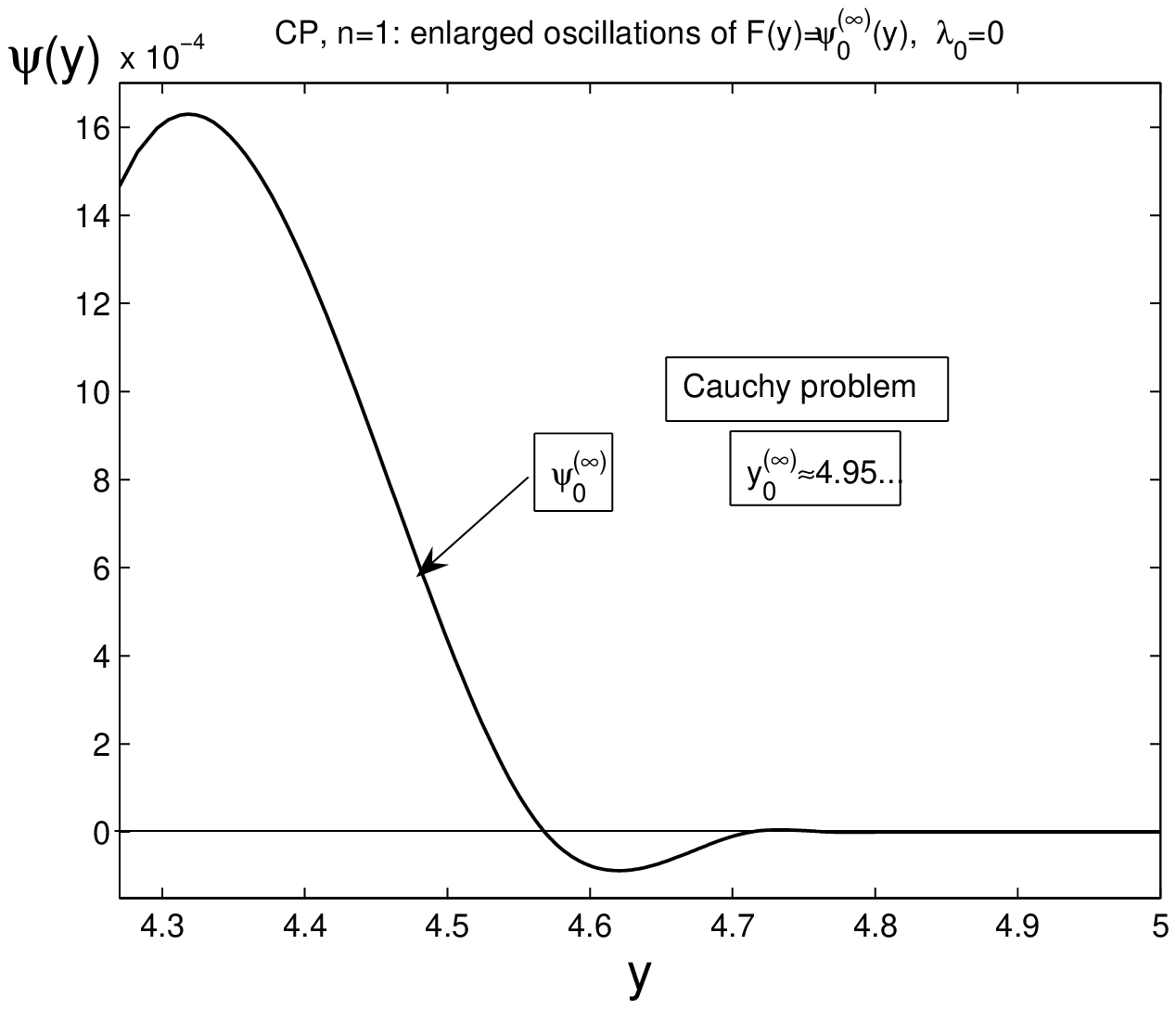}
} \subfigure[$\sim 10^{-6}$]{
\includegraphics[scale=0.52]{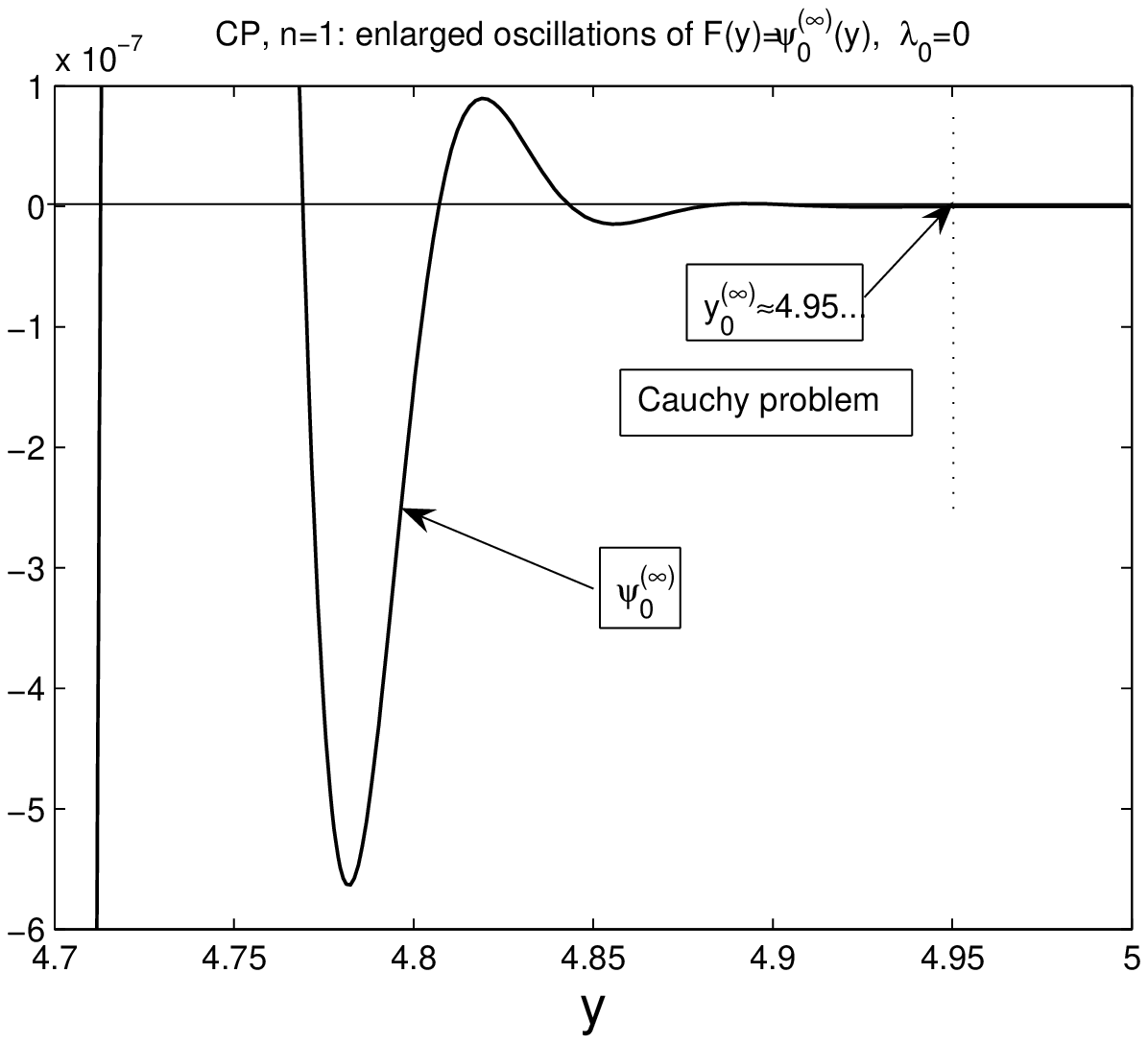}
}
 \vskip -.4cm
\caption{\rm\small The CP-eigenfunction of (\ref{.10}), with
$n=1$, for $k=\iy$, with the interface at $y_0^{(\iy)}
=4.95...\,$.}
 \label{F04}
\end{figure}

  In Figure \ref{F03}, we show the first two eigenfunctions of the
  problem (\ref{.10}) for $n=1$. It is difficult to demonstrate
  more functions, since the next ones are already very close to
  the Cauchy profile denoted by $\psi_0^{(\iy)}$. The oscillations of
  this CP-profile $F(y)$ are presented in Figure \ref{F04}, in between of the humps of which
  the interfaces of further nonlinear eigenfunction $\psi_0^{(k)}(y)$ are assumed to be situated.


\ssk

\begin{figure}
\centering \subfigure[four eigenfunctions]{
\includegraphics[scale=0.52]{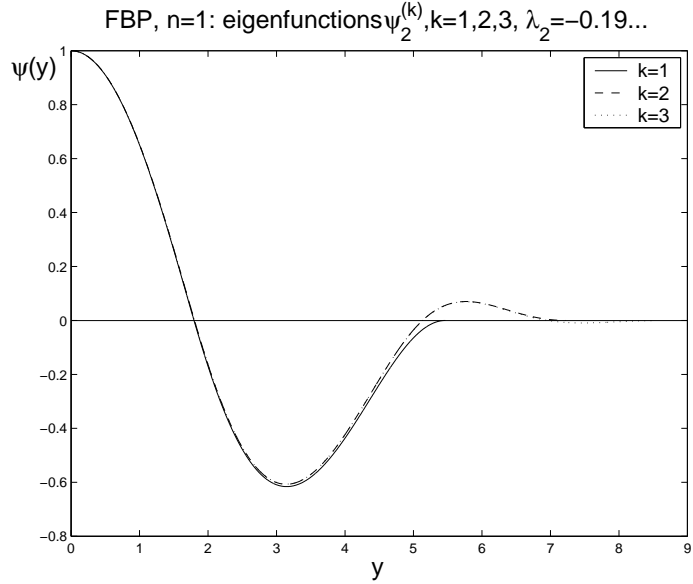}
} \subfigure[enlarged zeros]{
\includegraphics[scale=0.52]{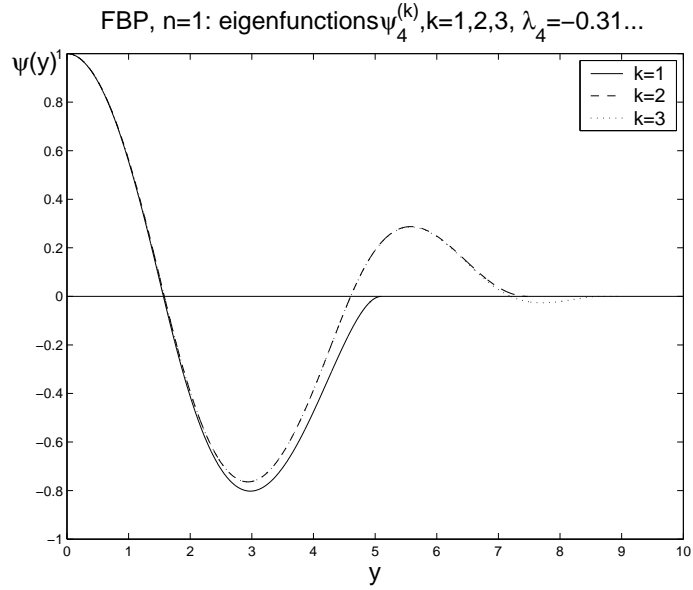}
}
 \vskip -.4cm
\caption{\rm\small The first three profiles for the second and third even eigenfunctions $m=2,4$ of (\ref{lnotzeron1}) for $n=1$. These numerical solutions of (\ref{lnotzeron1}) used the regularisation (\ref{eq:reg}) with $\delta=0.1$. }
 \label{F2F4n1}
\end{figure}

For $\lambda \neq 0$, the corresponding statement to (\ref{.10}) is
\be
  \label{lnotzeron1}
   \left\{
   \begin{matrix}
     -|f|^n f''' + \frac{1+n\lambda}{n+4} \,y f = \lambda \int\limits_{0}^{y} f(s)\, {\mathrm d}s \,\,\,\mbox{on}
     \,\,\,(0,y_+),\qquad\qquad\,\,\,\,\, \ssk \\
     f(0)=1, \,\,\, f'(0)=0; \quad f(y_+)=f'(y_+)=|f|^n f'''(y_+)=0,
     \end{matrix}
     \right.
      \ee
with $\lambda = \frac{1-\alpha(n+4)}4$.
 Again the sets of even eigenfunctions are denoted by $\{ f=\psi_m,
\lambda=\lambda_m, y_+=y_m \}$, using the subscript
$m=0,2,4,\ldots$. Within each set (i.e. $m$ fixed), we have a
countable set of eigenfunctions $\left\{
\psi_m=\psi_m^{(k)},y_m=y_m^{(k)} \right\}, k=1,2,3,...\,$. The
first three eigenfunctions for the second $m=2$ and the fourth
$m=4$ even sets are shown in Figure \ref{F2F4n1} for the case
$n=1$. In the nonlinear case $n>0$, unlike the linear one $n=0$
above, convergence of numerical methods become much more slower.


\subsection{Comments on $n$-branching and $p$-bifurcations}

These properties are assumed to be similar to those
developed earlier for the CP-setting. However, there are
essential difficulties even in doing some formal computations.

Overall, we expect that the nonlinear eigenvalue problem
(\ref{.10}) (see Conjecture \ref{Sectp0}.2) admits $n$-branching as $n \to 0^+$ from linear
eigenfunctions from Proposition \ref{Pr.01}. As usual ``variations" of the free boundary then must be taken into accound, which does not lead to any extra difficulty.

 Furthermore, the general VSS problem (\ref{.1}) possesses $p$-bifurcation branches
obtained via nonlinear bifurcations, whose theory is developed along the same lines used in
Section \ref{S5.2}. In particular, as in (\ref{pp1}), the first
bifurcation exponent will be
 $$
 \mbox{$
 p_0(n)= \frac {n+2}2 + \frac 1{2 \a_0(n)} \equiv  \frac
 {n+2}2+\frac{n+4}2=n+3.
 $}
 $$
 Similar to the analysis in Section \ref{S5.2}, further study is
 necessary to check whether a ``nonlinear bifurcation" occurs at
 $p=p_0$ (possibly not as for $n=0$), so further critical
 exponents $\a=\a_l(n)$ should be revealed.



We also expect that the above set of nonlinear
eigenfunction-eigenvalue $\{f_l(y),\,\a_l(n)\}$ pairs for the CP
for the TFE (\ref{tf1}) also play a role for the FBP.
 More precisely, it is expected that, for any given  CP profile $f_l(y)$,
 with $\a=\a_l(n)$, there exists a countable set of FBP profiles
 $\{f_l^j(y),\, \a_l(n)\}$, which are defined on the expanding
 intervals $[0,y_l^j)$, where $y_0^j \to y_0^{(\iy)}$ (the interface of the CP-profile)
  as $j \to \iy$.
 Eventually, there holds:
  \be
  \label{tt1}
   \mbox{$
   f_l^j(y) \to f_l(y) \quad \mbox{as} \quad j \to \iy
    $}
    \ee
    uniformly on compact subsets. In other words, each CP profile
 $f_l(y)$ can be arbitrarily closely approximated by FBP ones with a finite number of sign changes.
 The convergence (\ref{tt1}) is a difficult open problem, both
 analytically and numerically even for the first values of $l$.
 Of course, (\ref{tt1}) is associated with the discovered
 oscillatory properties of all the VSS CP profiles.








Again, we mention that the FBP problems have a different nature and are
more difficult than the CP, since a new free parameter
the interface location $y_0>0$ (an extra ``nonlinear eigenvalue")
appears in the mathematical setting.





\section{Discussion}
 \label{SectDi}


 We began in Section \ref{Sect2} with the similarity analysis of global, source-type (for
 $p=p_0$) and very singular  (for $p \not = p_0$)
 solutions of the TFE with a stable parabolic term (\ref{GPP}).
 Treating the Cauchy problem   (respectively, the FBP), we first described
 in Section \ref{SectCP1} for the CP
  (respectively,  \ref{SectLocR} for the FBP) the local
asymptotic
 properties of solutions near interfaces. While the FBP asymptotics
 turned out to be standard and reasonably well known since the 1990s,
   it is important that, for the CP, we
 detected the nonlinearity range $n \in (0,n_{\rm h})$ ($n_{\rm h}$ is the point
 of a  {\em heteroclinic bifurcation} for a nonlinear related ODE), in which the rescaled ODE
 exhibits a unique stable periodic motion describing, as expected, generic
 changing sign properties of more general solutions.

In Section \ref{SectCP1}, we studied the CP in the critical case
 $$
  \mbox{$
 p=p_0=n+1 + \frac 2N,
  $}
 $$
 where we  detected continuous branches of similarity profiles.
In Section \ref{CPTFEn}, we  developed $n$-branching theory of
similarity solutions of the 1D pure TFE
 $$
  u_t= - (|u|^n u_{xxx})_x,
   $$
 where we described branching of nonlinear similarity profiles from
 eigenfunctions of a linear rescaled operator at $n=0$.  This
 allowed us in Section \ref{CPpn} to reveal a countable sequence of critical exponents
 $\{p_l\}$ of the original stable TFE (\ref{GPP}) and to describe
 similarity solutions for $p \not = p_0$.

After a detailed study of the CP, we returned to the FBP setting.
We  studied in Section \ref{Sect5}
various branches of similarity patterns for the FBP in the
critical case $p=p_0$ and extended some of the
results to $p \not = p_0$ in  Section \ref{Sectp0}.

By comparing the similarity patterns of the CP and the FBP,
a striking ``limit" property emerges: the infinitely
oscillatory patterns of the CP are the limits of FBP-patterns with
a finite number of sign changes. Naturally, this is required  to take
into account sign changing patterns of the FBP, which have not been previously
studied in any detail.
 In a sense, the above {\em limit property}  can be considered
 as a certain definition of solutions of the CP (besides the
 already existing ones via {\em maximal regularity} at the
 interfaces or via a smooth analytic ``homotopy" deformation to
  the bi-harmonic equation
  $u_t=-\D^2 u$, which turned out to be a good approximation
of the TFE for small $n>0$; \cite{Gl4}).

Thus,  the goal of the paper was to describe some leading key
ideas concerning (i) formation of similarity solutions for the
stable TFEs and (ii) extra relations between the CP and the FBP.
Some of the most difficult
conclusions remain formal, the related
mathematics turns out to be very difficult, and we have posed several
open problems for future research.




\end{document}